\documentclass{svmult}
\usepackage{amsmath,amssymb,amsfonts}
\usepackage{verbatim}

\spnewtheorem{theorem*}{Theorem}[section]{\bf}{\it}
\spnewtheorem{lemma*}[theorem*]{Lemma}{\bf}{\it}
\spnewtheorem{proposition*}[theorem*]{Proposition}{\bf}{\it}
\spnewtheorem{corollary*}[theorem*]{Corollary}{\bf}{\it}
\spnewtheorem{definition*}[theorem*]{Definition}{\bf}{\it}

\def\a{\mathfrak{t}}
\def\A{{\rm A}}
\def\ab{\bar{a}}
\def\al{\alpha}
\def\At{\widetilde A}

\def\b{\mathfrak{b}}
\def\B{{\rm B}}
\def\Bb{\ts\overline{\ns\rm B\!}\,}
\def\be{\beta}
\def\bet{\widetilde{\be}}
\def\bg{\tilde g}
\def\bi{\tilde\imath}
\def\bj{\tilde\jmath}
\def\bh{\tilde h}
\def\bk{\tilde k}
\def\bl{\tilde l}

\def\bt{\ts\,\raise-0.5pt\hbox{\small$\boxtimes$}\,\,}
\def\Bt{\widetilde B}
\def\C{{\rm C}}
\def\cb{\bar{c}}
\def\CC{\mathbb{C}}
\def\com{\ts,\hskip-.5pt}
\def\Ct{\widetilde C}

\def\d{\partial}
\def\db{\bar{d}}
\def\D{\mathfrak{A}}
\def\de{\delta}
\def\De{\Delta}
\def\Dt{\widetilde D}

\def\E{\mathcal{E}}
\def\Eb{\widetilde{\mathcal{E}}}
\def\End{\operatorname{End}\ts}
\def\ep{\varepsilon}
\def\Ep{E^{\ts\prime}}
\def\Et{\widetilde E}

\def\f{\mathfrak{f}}
\def\F{\mathcal{F}}
\def\Fp{F^{\,\prime}}
\def\Ft{\widetilde{F}}
\def\Ftp{\widetilde{F}^{\,\prime}}

\def\g{\mathfrak{g}}
\def\ga{\gamma}
\def\ge{\geqslant}
\def\gl{\mathfrak{gl}}

\def\h{\mathfrak{h}}
\def\H{\mathfrak{H}}
\def\Hh{\mathfrak{B}}
\def\Hom{\operatorname{Hom}}
\def\Hp{\H^{\ts\prime}}

\def\I{{\rm I}}
\def\Ib{\tts\overline{\nns\rm I\nns}\tts}
\def\io{\iota}

\def\J{{\rm J}}
\def\Jb{\,\overline{\!\rm J\ns}\ts}
\def\Jp{{\rm J}^{\ts\prime}}
\def\Jpb{\,\overline{\!\rm J\ns}\ts^{\ts\prime}}

\def\ka{\kappa}

\def\la{\lambda}
\def\las{\la^{\ts\prime}}
\def\lat{\bar\la}
\def\lcd{\ts,\ldots,}
\def\le{\leqslant}

\def\mus{\mu^{\ts\prime}}
\def\mut{\bar\mu}

\def\n{\mathfrak{n}}
\def\nns{\hskip-.5pt}
\def\np{\n^{\ts\prime}}
\def\ns{\hskip-1pt}
\def\nus{\nu^{\ts\prime}}
\def\nut{\bar{\nu}}

\def\om{\omega}
\def\op{\oplus}
\def\ot{\otimes}

\def\p{\mathfrak{p}}
\def\P{\mathcal{G}}
\def\Pb{\bar P}
\def\PD{\mathcal{GD}}
\def\ph{\chi}
\def\Pp{P^{\ts\prime}}

\def\q{\mathfrak{q}}
\def\qp{\q^{\ts\prime}}
\def\Qb{\bar Q}

\def\r{\mathfrak{r}}
\def\rp{\r^{\ts\prime}}
\def\Rp{R^{\ts\prime}}

\def\s{\mathfrak{s}}
\def\S{\operatorname{S}}
\def\si{\sigma}
\def\sib{\ts\overline{\ns\si\ns}\ts}
\def\sih{\widetilde{\si}}
\def\sip{\si^{\ts\prime}}
\def\so{\mathfrak{so}}
\def\sp{\mathfrak{sp}}
\def\Sp{S^{\ts\prime}}
\def\St{\widetilde{S}}
\def\Sym{\mathfrak{S}}

\def\th{\theta}
\def\Tp{T^{\ts\prime}}
\def\ts{\hskip1pt}
\def\tts{\hskip.5pt}

\def\uo{\ts\overline{\ns u\nns}\,}
\def\U{\operatorname{U}}
\def\Uhb{\,\overline{\!\U(\h)\!\!\!}\,\,\,}

\def\vuo{\ts\overline{\ns v\ot u\ns}\,}

\def\Wb{\ts\overline{\ns W^{\phantom{!}}\!\!}}
\def\Wt{\widetilde{W}}

\def\X{\operatorname{X}}
\def\xic{\check\xi}

\def\Y{\operatorname{Y}}

\def\Z{\operatorname{Z}}
\def\ZZ{\mathbb{Z}}

\begin{document}


\title*{Twisted Yangians and Mickelsson Algebras II}
\author{Sergey Khoroshkin\inst{1}
\and Maxim Nazarov\inst{2}}
\institute{
Institute for Theoretical and Experimental Physics,
Moscow 117259,
Russia;
\\
\texttt{khor@itep.ru}
\and
Department of Mathematics,
University of York,
York YO10 5DD,
England;
\\
\texttt{mln1@york.ac.uk}
}

\maketitle


\smartqed
\renewcommand{\theequation}{\thesection.\arabic{equation}}
\makeatletter
\@addtoreset{equation}{section}
\makeatother


\section*{\normalsize 0.\ Introduction}
\setcounter{section}{0}
\setcounter{equation}{0}

This article is a continuation of our work \cite{KN2} which
concerned two known functors. The definition of one of these two
functors belongs to V.\,Drinfeld \cite{D2}. Let $\D_N$ be the
\textit{degenerate affine Hecke algebra\/} corresponding to the
general linear group $GL_N$ over a non-Archimedean local field. This
is an associative algebra over the 
field $\CC$
which contains the symmetric group ring $\CC\,\Sym_N$ as a subalgebra.
Let $\Y(\gl_n)$ be the
\textit{Yangian\/} of the general linear Lie algebra $\gl_n\ts$.
This is a deformation of the universal enveloping algebra
of the polynomial current Lie algebra
$\gl_n[u]$ in the class of Hopf algebras \cite{D1}.
It contains the universal enveloping algebra $\U(\gl_n)$ as a subalgebra.
There is also a homomorphism of associative algebras
$\Y(\gl_n)\to\U(\gl_n)$ identical on the subalgebra
$\U(\gl_n)\subset\Y(\gl_n)\ts$.
In \cite{D2} for any $\D_N\ts$-module $M\ts$, an action of the algebra
$\Y(\gl_n)$ was defined on the vector space
$(\ts M\ot(\CC^{\ts n})^{\ot N}){}^{\ts\Sym_N}_{\,-}$
of the diagonal skew $\Sym_N$-invariants in the tensor
product of the vector spaces $M$ and $(\CC^{\ts n})^{\ot N}\ts$.
Thus one gets a functor from the category of all
$\D_N\ts$-modules to the category of $\Y(\gl_n)\ts$-modules,
the \textit{Drinfeld functor}.

In \cite{KN1} we studied the composition of the Drinfeld functor
with another functor, introduced by I.\,Cherednik \cite{C}.
This second functor was also studied by
T.\,Arakawa, T.\,Suzuki and A.\,Tsuchiya \cite{A,AS,AST}.
For any module $U$ over the Lie algebra $\gl_{\ts l}\ts$,
an action of the algebra $\D_N$ can be defined on the
tensor product $U\ot(\CC^{\ts l})^{\ot N}$
of $\gl_{\ts l}\ts$-modules. This action of $\D_N$
commutes with the diagonal action of $\gl_{\ts l}$ on the tensor product.
Thus one gets a functor from the category of all
$\gl_{\ts l}\ts$-modules to the category of bimodules over $\gl_{\ts l}$
and $\D_N\ts$, the \textit{Cherednik functor}.
By applying the Drinfeld functor to the $\D_N\ts$-module
$M=U\ot(\CC^{\ts l})^{\ot N}\ts$, one turns to an $\Y(\gl_n)\ts$-module
the vector space
$$
(\ts U\ot(\CC^{\ts l})^{\ot N}\ot(\CC^{\ts n})^{\ot N}){}^{\ts\Sym_N}_{\,-}
=
U\ns\ot\ts\operatorname{\mathrm{\Lambda}}^N(\CC^{\ts l}\ot\CC^{\ts n})\ts.
$$
The action of the associative algebra
$\Y(\gl_n)$ on this vector space
commutes with the action of 
$\gl_{\ts l}\ts$.
By taking the direct sum of these $\Y(\gl_n)\ts$-modules
over $N=0,1\lcd n$
we turn to an $\Y(\gl_n)\ts$-module
the space $U\ot\ts\mathrm{\Lambda}\ts(\CC^{\ts l}\ot\CC^{\ts n})\ts$.
It is also a $\gl_{\ts l}\ts$-module; denote this bimodule by
$\E_{\ts l}\ts(U)\ts$. We identify the exterior algebra
$\mathrm{\Lambda}\ts(\CC^{\ts l}\ot\CC^{\ts n})$ with the
Grassmann algebra $\P\ts(\CC^{\ts l}\ot\CC^{\ts n})\ts$,
and denote by $\PD\ts(\CC^{\ts l}\ot\CC^{\ts n})$ the ring of
$\CC\ts$-endomorphisms of $\P\ts(\CC^{\ts l}\ot\CC^{\ts n})\ts$.
The action of the Yangian $\Y(\gl_n)$ on its module
$\E_{\ts l}\ts(U)$ is then determined by a homomorphism
$\al_{\ts l}:\Y(\gl_n)\to
\U(\gl_{\ts l})\ot\ts\PD\ts(\CC^{\ts l}\ot\CC^{\ts n})\ts$,
see Proposition~\ref{dast} below.

Now let $\f_m$ be either the orthogonal Lie algebra $\so_{2m}$
or the symplectic Lie algebra $\sp_{2m}\ts$.
The first objective of the present article
is to define analogues of the functor $\E_{\ts l}$
and of the homomorphism $\al_{\ts l}$
for the Lie algebra $\f_m$ instead of $\gl_{\ts l}\ts$.
The role of the Yangian $\Y(\gl_n)$ is played here by the
\textit{twisted Yangian\/} $\Y(\g_n)\ts$, which is
a right coideal subalgebra of the Hopf algebra $\Y(\gl_n)\ts$.
Here $\g_n$ is a Lie subalgebra of $\gl_n\ts$, orthogonal
in the case $\f_m=\so_{2m}$ and symplectic in the case $\f_m=\sp_{2m}\ts$;
in the latter case $n$ has to be even.
Let the superscript ${}^{\ts\prime}$ indicate the transposition
in $\gl_n$ relative to the bilinear form on $\CC^n$
preserved by the subalgebra $\g_n\subset\gl_n\ts$, so that
$
\g_n=\{\ts A\in\gl_n\,\ts|\,A^{\ts\prime}=-A\ts\}\ts.
$
As an associative algebra, 
$\Y(\g_n)$ is a deformation of the
universal enveloping algebra of the
\textit{twisted polynomial current Lie algebra\/}
$$
\{\ts A(u)\in\gl_n[u]\,\ts|\,A^{\ts\prime}(u)=-A(-u)\ts\}\ts.
$$

Twisted Yangians were introduced
by G.\,Olshanski \cite{O2},
their structure has been studied in \cite{MNO}.
In Section 2 of the the present article we introduce a homomorphism
$\Y(\g_n)\to\U(\ts\f_m)\ot\ts\PD\ts(\CC^{\ts m}\ot\CC^{\ts n})\ts$,
see our Propositions~\ref{xb} and \ref{tb}.
The image of $\Y(\g_n)$ under this homomorphism commutes with the image
of the algebra $\U(\ts\f_m)$ under its diagonal embedding
\eqref{xact}
to the tensor product $\U(\ts\f_m)\ot\ts\PD\ts(\CC^{\ts m}\ot\CC^{\ts n})\,$;
here we use the homomorphism
$\zeta_{\ts n}:\U(\ts\f_m)\to\PD\ts(\CC^{\ts m}\ot\CC^{\ts n})$
defined by \eqref{gan}.
The twisted Yangian $\Y(\g_n)$
contains the universal enveloping algebra $\U(\g_n)$ as a subalgebra.
There is also a homomorphism
$\pi_n:\Y(\g_n)\to\U(\g_n)$ identical on the subalgebra
$\U(\g_n)\subset\Y(\g_n)\ts$.
Our results extend the classical theorem \cite{H}
stating that the image of $\U(\ts\f_m)$ in
$\PD\ts(\CC^{\ts m}\ot\CC^{\ts n})$
under the homomorphism
$\zeta_n$ consists of all $G_n$-invariant elements.
Here $G_n$ is either the orthogonal or the
symplectic group, so that $\g_n$ is its Lie algebra;
the group $G_n$ acts on 
$\PD\ts(\CC^{\ts m}\ot\CC^{\ts n})$
via its natural action on $\CC^{\ts n}\ts$.

In the present article
we prefer to work with a certain central extension
$\X(\g_n)$ of the algebra
$\Y(\g_n)\ts$, called the
\textit{extended twisted Yangian}. Central elements
$\,O^{(1)}\ns\com O^{(2)}\ns\com\,\ldots\,$ of the algebra $\X(\g_n)$
generating the kernel of the canonical homomorphism $\X(\g_n)\to\Y(\g_n)$
are given in Section 1, together with the definitions of
$\X(\g_n)$ and $\Y(\g_n)\ts$. There is also a homomorphism
$\X(\g_n)\to\X(\g_n)\ot\Y(\gl_n)\ts$.
Using it, the tensor product
of any modules over the algebras $\X(\g_n)$ and $\Y(\gl_n)$
becomes another module over $\X(\g_n)\ts$.
Moreover, this homomorphism is a coaction
of the Hopf algebra $\Y(\gl_n)$ on the algebra $\X(\g_n)\ts$.
We define a homomorphism
$\be_m:\X(\g_n)\to\U(\ts\f_m)\ot\ts\PD\ts(\CC^{\ts m}\ot\CC^{\ts n})$
which is our analogue of the homomorphism $\al_{\ts l}\ts$, see
Proposition~\ref{xb}.
The image of $\X(\g_n)$ under $\be_m$ commutes with the image
of the algebra $\U(\ts\f_m)$ under its embedding \eqref{xact} to
$\U(\ts\f_m)\ot\ts\PD\ts(\CC^{\ts m}\ot\CC^{\ts n})\ts$.
The reason why we work with $\X(\g_n)$ rather than with
$\Y(\g_n)$ is explained in Section~2.

The generators of the algebra $\X(\g_n)$ appear as coefficients
of certain series $S_{ij}(u)$ in the variable $u$
where $i\com j=1\lcd n\ts$.
We define the homomorphism $\be_m$ by applying it to the
coefficients, and by giving the resulting series
$\be_m(S_{ij}(u))$ with coefficients in
$\U(\ts\f_m)\ot\ts\PD\ts(\CC^{\ts m}\ot\CC^{\ts n})$ explicitly.
Then we define another homomorphism
$$
\bet_m:\X(\g_n)\to\U(\ts\f_m)\ot\ts\PD\ts(\CC^{\ts m}\ot\CC^{\ts n})
$$
which factors through the canonical homomorphism
$\X(\g_n)\to\Y(\g_n)\ts$. 
Thus we obtain the homomorphism
$\Y(\g_n)\to\U(\ts\f_m)\ot\ts\PD\ts(\CC^{\ts m}\ot\CC^{\ts n})$
mentioned above.
Every series $\bet_m(S_{ij}(u))$ is the product
of $\be_m(S_{ij}(u))$ with a certain series
with coefficients from
$\Z(\ts\f_m)\ot1\ts$,
where $\Z(\ts\f_m)$ is the centre of the algebra $\U(\ts\f_m)\ts$.

The defining relations of the algebra $\X(\g_n)$ can be
written as the \textit{reflection equation} \eqref{rsrs}
on the $n\times n$ matrix
$S(u)$ whose $i\com j$ entry is the series $S_{ij}(u)$.
This terminology was introduced by physicists;
see for instance \cite{KS} and references therein.

Now let $V$ be any $\f_m\ts$-module. Using the homomorphism $\be_m\ts$,
we turn the vector space $V\ns\ot\ts\P\ts(\CC^{\ts m}\ot\CC^{\ts n})$
into a bimodule over $\f_m$ and $\X(\g_n)\ts$.
We denote this bimodule by $\F_m(V)\ts$.
The functor $\F_m$ is our analogue of the functor $\E_{\ts l}$ for
$\f_m$ instead of $\gl_{\ts l}\ts$.
When $m=0\ts$, we set $\F_{\ts0}(V)=\CC$ so that
$\be_{\ts0}$ is the composition of the canonical
homomorphism $\X(\g_n)\to\Y(\g_n)$ with the restriction of
the counit homomorphism $\Y(\gl_n)\to\CC$ to $\Y(\g_n)\ts$.

Here we show that the functor $\F_m$
shares the three fundamental properties of the functor $\E_{\ts l}$
considered in \cite{KN2}. The first of these properties of
$\E_{\ts l}$ concerns \textit{parabolic induction\/}
from the direct sum of Lie algebras
$\gl_m\op\gl_{\ts l}$ to $\gl_{m+l}\ts$.
Let $\p$ be the maximal parabolic subalgebra of $\gl_{m+l}$
containing the direct sum $\gl_m\op\gl_{\ts l}\ts$.
Let $\q\subset\gl_{m+l}$ be the Abelian subalgebra with
$
\gl_{m+l}=\q\op\p\ts.
$
For any $\gl_m\ts$-module $W$ let $W\bt U$ be the
$\gl_{m+l}\ts$-module
parabolically induced from the
$\gl_m\op\gl_{\ts l}\ts$-module $W\ot U\ts$.
This is a module induced from the subalgebra $\p\ts$.
Consider the space $\E_{\ts m+l}\ts(\ts W\ns\bt U\ts)_{\ts\q}$
of $\q\ts$-coinvariants of the $\gl_{m+l}\ts$-module
$\E_{\ts m+l}\ts(\ts W\ns\bt U\ts)\ts$.
This space is an $\Y(\gl_n)\ts$-module, which also
inherits the action of the Lie algebra $\gl_m\op\gl_{\ts l}\ts$.
The additive group $\CC$ acts on
the Hopf algebra $\Y(\gl_n)$ by automorphisms.
Let $\E_{\ts l}^{\,-z}\ts(U)$ be the $\Y(\gl_n)\ts$-module obtained
from $\E_{\ts l}\ts(U)$ by pulling it back through the automorphism
of $\Y(\gl_n)$ corresponding to 
$z\in\CC\ts$.
The automorphism itself is denoted by $\tau_{z}\ts$, see \eqref{tauz}.
Thus the underlying vector space of the $\Y(\gl_n)\ts$-module
$\E_{\ts l}^{\,-z}\ts(U)$ is $U\ot\ts\P\ts(\CC^{\ts l}\ot\CC^{\ts n})\ts$,
whereon the action of $\Y(\gl_n)$ is
defined by the composition of two homomorphisms,
\begin{equation}
\label{diag1}
\Y(\gl_n)
\underset{\tau_z}\longrightarrow
\Y(\gl_n)
\underset{\al_{\ts l}}\longrightarrow
\U(\gl_{\ts l})\ot\ts\PD\ts(\CC^{\ts l}\ot\CC^{\ts n})\ts.
\end{equation}
The target algebra here acts on $U\ot\ts\P\ts(\CC^{\ts l}\ot\CC^{\ts n})$
by definition. As a $\gl_{\ts l}\ts$-module $\E_{\ts l}^{\,-z}\ts(U)$
coincides with $\E_{\ts l}\ts(U)\ts$.
In \cite{KN2} we proved that
the bimodule $\E_{\ts m+l}\ts(\ts W\ns\bt U\ts)_{\ts\q}$
of $\Y(\gl_n)$ and $\gl_m\op\gl_{\ts l}$
is equivalent to 
$\E_{\ts m}(W)\ot\E_{\ts l}^{\ts m}\ts(U)\ts$.
We use the comultiplication on $\Y(\gl_n)\ts$.

Our Theorem~\ref{parind} is an analogue of this
comultiplicative property of 
$\E_{\ts l}\ts$. Take the maximal parabolic subalgebra
of the Lie algebra $\f_{m+l}$ containing the direct sum
$\f_m\op\gl_{\ts l}\,$; we do not exclude the case $m=0$ here.
Using that subalgebra, determine
the $\f_{m+l}\ts$-module $V\bt U$ parabolically induced from
the $\f_m\op\gl_{\ts l}\ts$-module $V\ot U\ts$.
Consider the space of coinvariants of the
$\f_{m+l}\ts$-module $\F_{m+l}\ts(\ts V\ns\bt U\ts)$
relative to the nilpotent subalgebra of $\f_{m+l}$
complementary to our parabolic subalgebra.
This space is a bimodule over $\f_m\op\gl_{\ts l}$
and $\X(\g_n)\ts$. We prove that this bimodule
is essentially equivalent to the tensor product
$\F_m(V)\ot\E_{\ts l}^{\ts z}\ts(U)$ with
$z=m-\frac12$ for $\f_m=\so_{2m}\ts$, and
$z=m+\frac12$ for $\f_m=\sp_{2m}\ts$.
More precisely, the underlying vector space of the $\X(\g_n)\ts$-module
$\F_m(V)\ot\E_{\ts l}^{\ts z}\ts(U)$ is
\begin{equation}
\label{diag2}
V\ot\ts\P\ts(\CC^{\ts m}\ot\CC^{\ts n})\ot
U\ot\ts\P\ts(\CC^{\ts l}\ot\CC^{\ts n})\ts,
\end{equation}
whereon the action of $\X(\g_n)$ is
defined by the composition of two homomorphisms,
$$
\X(\g_n)\to
\X(\g_n)\ot\Y(\gl_n)\to
\U(\ts\f_m)\ot\ts\PD\ts(\CC^{\ts m}\ot\CC^{\ts n})\ot
\U(\ts\gl_{\ts l})\ot\ts\PD\ts(\CC^{\ts l}\ot\CC^{\ts n})\ts.
$$
Here the first homomorphism is the
coaction $\Y(\gl_n)$ on $\X(\g_n)\ts$,
while the second one is the tensor product of the homomorphisms
$
\be_m:\X(\g_n)\to
\U(\ts\f_m)\ot\ts\PD\ts(\CC^{\ts m}\ot\CC^{\ts n})
$
and
$$
\al_{\ts l}\,\tau_{-z}:\Y(\gl_n)\to
\U(\ts\gl_{\ts l})\ot\ts\PD\ts(\CC^{\ts l}\ot\CC^{\ts n})\ts;
$$
see \eqref{diag1}.
By multiplying the image of 
$S_{ij}(u)\in\X(\g_n)[[u^{-1}]]$
under this composition by a certain series with the
coefficients from the subalgebra
$$
1\ot1\ot\Z(\gl_{\ts l})\ot1
\ts\subset\ts
\U(\ts\f_m)\ot\ts\PD\ts(\CC^{\ts m}\ot\CC^{\ts n})\ot
\U(\ts\gl_{\ts l})\ot\ts\PD\ts(\CC^{\ts l}\ot\CC^{\ts n})\,,
$$
we get another homomorphism
$
\X(\g_n)\to\U(\ts\gl_{\ts l})\ot\ts\PD\ts(\CC^{\ts l}\ot\CC^{\ts n})\ts.
$
The latter homorphism defines
another action of $\X(\g_n)$ on the
vector space \eqref{diag2}.
Theorem~\ref{parind} states that this action
is equivalent to the action of $\X(\g_n)$
on the space of coinvariants of $\F_{m+l}\ts(\ts V\ns\bt U\ts)\ts$.
Moreover, the actions of the direct summand $\f_m$ of
$\ts\f_m\op\gl_{\ts l}$
on $\F_m(V)\ot\E_{\ts l}^{\ts z}\ts(U)$
and on the space of coinvariants of
$\F_{m+l}\ts(\ts V\ns\bt U\ts)$ are also equivalent,
while the actions of the direct summand $\gl_{\ts l}$
differ only by the automorphism \eqref{autol} of
the Lie algebra $\gl_{\ts l}\ts$. Hence Theorem~\ref{parind}
describes the first fundamental property of the functor $\F_m\ts$.

Let us now discuss the second fundamental property of $\F_m\ts$.
In \cite{TV} V.\,Tarasov and A.\,Varchenko established a
correspondence between canonical intertwining operators
on $l\ts$-fold tensor products of certain $\Y(\gl_n)\ts$-modules,
and the \textit{extremal cocycle\/}
on the Weyl group $\Sym_{\ts l}$ of the reductive Lie algebra $\gl_{\ts l}$
defined by D.\,Zhelobenko \cite{Z1}. In \cite{TV}
each of the $l$ tensor factors
is obtained from one of $\gl_n\ts$-modules
$\operatorname{S\ts}^N(\CC^{\ts n})$
by pulling back through the homomorphism $\Y(\gl_n)\to\U(\gl_n)$
and then back through one of the automorphisms
$\tau_z:\Y(\gl_n)\to\Y(\gl_n)$.
Here $\operatorname{S\ts}^N(\CC^{\ts n})$
is the $N$-th symmetric power of the vector space $\CC^{\ts n}\ts$,
while the homomorphism $\Y(\gl_n)\to\U(\gl_n)$ is defined by \eqref{eval}.
In \cite{KN1} we gave a representation theoretic explanation
of that correspondence from \cite{TV}, by employing the theory of
Mickelsson algebras \cite{M1,M2} as developed in~\cite{KO}.

For any $N\in\{1\lcd n\}$ and any $z\in\CC\ts$
we denote by $P_z^{\ts N}$
the $\Y(\gl_n)\ts$-module obtained by pulling back the action
of $\U(\gl_n)$ on the subspace of $\P\ts(\CC^{\ts n})$
of degree $N$ through the homomorphism $\Y(\gl_n)\to\U(\gl_n)$
and then through the automorphism $\tau_{-\ts z}$ of $\Y(\gl_n)\ts$.
The action of the algebra $\Y(\gl_n)$ on $P_z^{\ts N}$
is defined by the composition of homomorphisms
\begin{equation}
\label{diag3}
\Y(\gl_n)
\underset{\tau_{-z}}\longrightarrow
\Y(\gl_n)
\to
\U(\gl_n)
\to
\PD\ts(\CC^{\ts n})\ts.
\end{equation}
Here the second homomorphism is the one defined by \eqref{eval}, while
the algebra $\PD\ts(\CC^{\ts n})$ acts on $\P\ts(\CC^{\ts n})$ naturally.
Using the functor $\E_{\ts l}\ts$, in \cite{KN2} we established
a correspondence between intertwining operators
on the $l\ts$-fold
tensor products of modules of the form $P_z^{\ts N}\ts$,
and the same extremal cocycle on $\Sym_{\ts l}$
as considered in \cite{KN1}.
This is an \lq\lq\ts antisymmetric\ts\rq\rq\
version of the
correspondence first established in \cite{TV}.
The parameters $z$ corresponding to the $l$ tensor factors
are in general position, that is their
differences do not belong to $\ZZ\ts$. Then each of
the tensor products is irreducible as $\Y(\gl_n)$-module \cite{NT}.
Hence the intertwining operators between them
are unique up to multipliers from $\CC\ts$.


In the present article
we show that the functor $\F_m$ plays a role
similar to that of $\E_{\ts l}\ts$,
when the Lie algebra $\gl_{\ts l}$ is replaced by $\f_m\ts$.
Namely,
we establish a correspondence between intertwining operators
of certain $\X(\g_n)\ts$-modules, and the
extremal cocycle
on the hyperoctahedral group $\H_m$ corresponding
to the reductive Lie algebra $\f_m\ts$.
Here $\H_m$ is regarded as the Weyl group of $\f_m=\sp_{2m}\ts$,
and as an extension of the Weyl group of $\f_m=\so_{2m}$
by a Dynkin diagram automorphism. In both cases, the definition
of the extremal cocycle
is essentially due to D.\,Zhelobenko \cite{Z1}.
However, the original extremal cocycle
has been defined on the Weyl group of $\f_m\ts$,
which in the case $\f_m=\so_{2m}$ is only a subgroup of $\H_m$ of index~2.
An extension of the original definition to the
whole group $\H_m$ was given in \cite{KN3}.
All necessary details on the extremal cocycle corresponding to $\f_m$
are also reviewed in Section~4 of the present article.

The twisted Yangian $\Y(\g_n)$ is determined 
by a distinguished involutive automorphism \eqref{transauto}
of the algebra $\Y(\gl_n)\ts$. The automorphism \eqref{transauto}
corresponds to the automorphism
$$
A(u)\mapsto-\ts A^{\ts\prime}(-u)
$$
of the Lie algebra $\gl_n[u]\ts$,
when the algebra $\Y(\gl_n)$ is regarded as a deformation of
the universal enveloping algebra of $\gl_n[u]\ts$.
By pulling the $\Y(\gl_n)\ts$-module $P_z^{\ts N}$ back through the
automorphism \eqref{transauto} we get another
$\Y(\gl_n)\ts$-module, which we denote by $P_z^{\ts-N}\ts$.
The underlying vector space of $P_z^{\ts-N}$ consists
of elements of $\P\ts(\CC^{\ts n})$ of degree $N$,
whereon the action of $\Y(\gl_n)$
is defined by the composition of four homomorphisms
$$
\Y(\gl_n)
\to
\Y(\gl_n)
\underset{\tau_{-z}}\longrightarrow
\Y(\gl_n)
\to
\U(\gl_n)
\to
\PD\ts(\CC^{\ts n})\ts.
$$
Here the first map is the automorphism \eqref{transauto},
the other three are the same as in~\eqref{diag3}.

Now take any $\nu_1\lcd\nu_m\in\{1\lcd n\}$
and any $z_1\lcd z_m\in\CC$ such that
$z_a-z_b\notin\ZZ$
and
$z_a+z_b\notin\ZZ$
when
$a\neq b\ts$.
In the case
$\f_m=\sp_{2m}$ we also assume that $2\ts z_a\notin\ZZ$
for any $a\ts$.
The hyperoctahedral group $\H_m$ can be realized
as the group of all permutations $\si$ of
$-\ts m\lcd-1\com1\lcd m\ts$
such that $\si\ts(-c)=-\ts\si\ts(c)$ for any $c\ts$.
In Section~5 of the present article,
we show how the value of the extremal cocycle
for the Lie algebra $\f_m$ at an element $\si\in\H_m$
determines an intertwining operator of $\X(\g_n)\ts$-modules
\begin{equation}
\label{nuz}
P_{z_m}^{\,\nu_m}
\ot\ldots\ot
P_{z_1}^{\,\nu_1}
\,\to\,
P_{\widetilde{z}_m}^{\,\de_m\ts\widetilde{\nu}_m}
\ot\ldots\ot
P_{\widetilde{z}_1}^{\,\de_1\ts\widetilde{\nu}_1}
\end{equation}
where
\begin{equation}
\label{defdelta}
\widetilde\nu_a=\nu_{\ts|\si^{-1}(a)|}\ts,
\quad
\widetilde{z}_a=z_{\ts|\si^{-1}(a)|}
\quad\textrm{and}\quad
\de_a=\operatorname{sign}\si^{-1}(a)
\end{equation}
for each $a=1\lcd m\ts$.
The tensor products in \eqref{nuz} are those of
$\Y(\gl_n)\ts$-modules. By restricting both
tensor products to the subalgebra
$\Y(\g_n)\subset\Y(\gl_n)$ and by pulling the
restrictions back through the canonical homomorphism
$\X(\g_n)\to\Y(\g_n)$,
both tensor products in \eqref{nuz} become $\X(\g_n)\ts$-modules.
Thus the actions of the algebra $\X(\g_n)$ on both tensor products
in \eqref{nuz} are obtained by using the composition
$$
\X(\g_n)\to\Y(\g_n)\to\Y(\gl_n)\to\Y(\gl_n)^{\ot\ts n}.
$$
Here the first map is the canonical homomorphism,
the second is the embedding defining $\Y(\g_n)$, while the third
is $m$ fold comultiplication.
It was proved in \cite{MN} that
under our assumptions on $z_1\lcd z_m$
both 
tensor products in \eqref{nuz}
are irreducible $\X(\g_n)\ts$-modules, equivalent to each other.
Hence an intertwining operator between them
is unique up to a multiplier from $\CC\ts$.
For our operator,
this multiplier is determined by Proposition~\ref{isis}.

To obtain our intertwining
operator \eqref{nuz}
we use the theory of Mickelsson algebras, like we did in \cite{KN1,KN2}.
Our particular Mickelsson algebra is determined by the pair
formed by the tensor product
$\U(\ts\f_m)\ot\ts\PD\ts(\CC^{\ts m}\ot\CC^{\ts n})$
and by its subalgebra
$\U(\ts\f_m)$ relative to the 
embedding \eqref{xact}.
The extended twisted Yangian $\X(\g_n)$ appears naturally here,
because its image relative to 
$\be_m$ commutes with the image of $\U(\ts\f_m)$ in the tensor product.
Another expression for an intertwining operator
\eqref{nuz} was given in \cite{N2}.

In Section 2 we choose a triangular decomposition \eqref{tridec}
of the Lie algebra $\f_m$ into a direct sum
of a Cartan subalgebra $\h$ and of two maximal nilpotent subalgebras
$\n\,,\np\ts$.
For any formal power series $f(u)$ in $u^{-1}$ with coefficients
from $\CC$ and leading term $1$, the assignments \eqref{fus}
define an automorphism of the algebra $\X(\g_n)\ts$.
Up to pulling it back through such an automorphism,
the source $\X(\g_n)\ts$-module in \eqref{nuz}
arises as the space of $\n\ts$-coinvariants of weight $\la$
for the $\f_m\ts$-module $\F_m(M_\mu)\ts$,
where $M_\mu$ is the Verma module over $\f_m$ with the
highest vector of weight $\mu$ annihilated by the
action of the subalgebra $\np\subset\f_m\ts$.
The weights $\la$ and $\mu$
relative to the Cartan subalgebra $\h$
are determined here by the parameters $\nu_1\lcd\nu_m$
and $z_1\lcd z_m$ from \eqref{nuz}.
We denote the space of $\n\ts$-coinvariants of weight $\la$
by $\F_m(M_\mu)_{\ts\n}^{\ts\la}\ts$.
The algebra $\X(\g_n)$ acts on the latter space, because
the action of $\X(\g_n)$ on $\F_m(M_\mu)$ commutes with that of $\f_m\ts$.
We prove that
the above defined action of the algebra $\X(\g_n)$ on the source
tensor product in \eqref{nuz} is equivalent to the action
on the vector space of $\F_m(M_\mu)_{\ts\n}^{\ts\la}\ts$,
defined by the composition 
\begin{equation}
\label{diag4}
\X(\g_n)\to\X(\g_n)\to\End(\ts\F_m(M_\mu))\ts.
\end{equation}
Here the first map is the automorphism \eqref{fus} where
$f(u)^{-1}$ equals the product \eqref{muprod}.
The second map here is the defining homomorphism of
the $\X(\g_n)\ts$-module $\F_m(M_\mu)\ts$.

To get the target $\X(\g_n)\ts$-module in \eqref{nuz},
we generalize our definition of the functor $\F_m\ts$.
In the beginning of Section 5,
for any sequence $\de=(\ts\de_1\lcd\de_m)$ of $m$
elements of the set $\{1\com-1\}$ we define a functor
$\F_\de\,$, with the same source and target categories as
the functor $\F_m$ has. Moreover, for any $\f_m\ts$-module
$V$ the underlying vector spaces
of the bimodules $\F_\de(V)$ and $\F_m(V)$ are the same,
that is $V\ot\P\ts(\CC^{\ts m}\ot\CC^{\ts n})\ts$.
The actions of $\f_m$ and $\X(\g_n)$ on $\F_\de(V)$
are obtained by pushing forward the defining homomorphisms
$$
\zeta_n:\U(\ts\f_m)\to\PD\ts(\CC^{\ts m}\ot\CC^{\ts n})
\quad\text{and}\quad
\be_m:\X(\g_n)\to\U(\ts\f_m)\ot\ts\PD\ts(\CC^{\ts m}\ot\CC^{\ts n})
$$
through a certain authomorphism $\varpi$ of the ring
$\PD\ts(\CC^{\ts m}\ot\CC^{\ts n})$ depending on $\de\ts$.
Namely, the automorphism $\varpi$ is defined by the assignments
\eqref{compfour}. Thus to define the functor $\F_\de\,$, we use
the compositions $\varpi\,\zeta_n$ and $(1\ot\varpi)\,\be_m$
instead of the homomorphisms $\zeta_n$ and $\be_m$ respectively.
In particular, we have $\F_\de\ts(V)=\F_m(V)$ for
the sequence $\de=(1\lcd1)\ts$.

Up to pulling it back through an automorphism of the form \eqref{fus},
the target $\X(\g_n)\ts$-module in \eqref{nuz}
arises as the space of $\n\ts$-coinvariants of weight $\si\ts\circ\la$
for the $\f_m\ts$-module $\F_\de\ts(M_{\ts\si\ts\circ\ts\mu})\ts$.
The sequence $\de=(\ts\de_1\lcd\de_m)$ is as defined in
\eqref{defdelta}, and
the symbol $\circ$ here indicates the \textit{shifted action\/}
of the group $\H_m$ on the weights of $\h\ts$.
Our Proposition \ref{siverma} states that 
action of the algebra $\X(\g_n)$ on the target
tensor product in \eqref{nuz} is equivalent to the action
on the vector space of
$\F_\de\ts(M_{\ts\si\ts\circ\ts\mu})_{\ts\n}^{\ts\si\ts\circ\ts\la}\ts$,
defined by the composition 
\begin{equation}
\label{diag5}
\X(\g_n)\to\X(\g_n)\to\End(\ts\F_m(M_{\ts\si\ts\circ\ts\mu}))\ts.
\end{equation}
Here the first map is the automorphism \eqref{fus} where
$f(u)^{-1}$ equals the product \eqref{muprod}.
The second map here is the defining homomorphism of
the $\X(\g_n)\ts$-module $\F_m(M_{\ts\si\ts\circ\ts\mu}\ts)\ts$.

In Section 5 we show that value of the extremal cocycle
for the Lie algebra $\f_m$ at the element $\si\in\H_m$
determines an intertwining operator of $\X(\g_n)\ts$-modules
\begin{equation}
\label{diag6}
\F_m(M_\mu)_{\ts\n}^{\ts\la}\,\to\,
\F_\de\ts(M_{\ts\si\ts\circ\ts\mu})_{\ts\n}^{\ts\si\ts\circ\ts\la}\ts.
\end{equation}
The product \eqref{muprod} does not depend on the element
$\si\in\H_m\ts$, so that the authomorphisms \eqref{fus}
of the algebra $\X(\g_n)$ in \eqref{diag4} and \eqref{diag5} are the same.
Hence by replacing the source and the target $\X(\g_n)\ts$-modules
by their equivalent modules, we obtain our intertwining
operator \eqref{nuz}. The role played by the functor
$\F_m$ in this construction of the operator \eqref{nuz} is the
second fundamental property of that functor.

The third fundamental property of the functor $\E_{\ts l}$ considered
in \cite{KN2} is its connection with the
\textit{centralizer construction\/}
of the Yangian $\Y(\gl_n)$ proposed by G.\,Olshanski \cite{O1}.
For any two irreducible polynomial modules
$U$ and $U^{\ts\prime}$ over the Lie algebra $\gl_{\ts l}\ts$,
the results of \cite{O1} provide an action
of $\Y(\gl_n)$ on the vector space
\begin{equation}
\label{diag7}
\Hom_{\,\gl_{\ts l}}(\ts U^{\ts\prime}\ts,
U\ot\P\ts(\ts\CC^{\ts l}\ot\CC^{\ts n}\ts))\ts.
\end{equation}
Moreover, this action is irreducible. In \cite{KN2}
we proved that the same action of $\Y(\gl_n)$ on the vector space
\eqref{diag7} is obtained when the target $\gl_{\ts l}\ts$-module
$U\ot\P\ts(\ts\CC^{\ts l}\ot\CC^{\ts n}\ts)$ in \eqref{diag7}
is regarded as the bimodule $\E_{\ts l}\ts(U)$
over $\Y(\gl_n)$ and $\gl_{\ts l}\ts$.

There is a centralizer construction of 
$\Y(\g_n)$ again due
to G.\,Olshanski \cite{O2}, see also \cite{MO} and Section 6 here.
That construction
served as a motivation for introducing the twisted Yangians.
For any irreducible finite-dimensional modules
$V$ and $V^{\ts\prime}$ of the Lie algebra $\f_m\ts$,
the results of \cite{O2} provide an action
of the algebra $\X(\g_n)$ on the vector~space
\begin{equation}
\label{vvp}
\Hom_{\,\f_m}(\ts V^{\ts\prime}\ts,
V\ot\P\ts(\ts\CC^{\ts m}\ot\CC^{\ts n}\ts))\,.
\end{equation}
The group $G_n$ also acts on this vector space,
via its natural action on $\CC^{\ts n}$.

When $\g_n$ is an
orthogonal Lie algebra, the space \eqref{vvp} is irreducible
under the joint action of $\X(\g_n)$ and $G_n\ts$.
When $\g_n$ is symplectic, 
\eqref{vvp} is irreducible under the action of the $\X(\g_n)$ alone.
Our Theorem~\ref{5.1} states that the action
of $\X(\g_n)$ on \eqref{vvp} is essentially the same as
the action obtained from the bimodule
$\F_m(V)=V\ns\ot\ts\P\ts(\CC^{\ts m}\ot\CC^{\ts n})$
of $\X(\g_n)$ and $\f_m\ts$. More precisely, the action
of $\X(\g_n)$ on the vector space \eqref{vvp}
provided by \cite{O2} can also be obtained
from an action of $\X(\g_n)$ on the target
$\f_m\ts$-module $V\ot\P\ts(\ts\CC^{\ts m}\ot\CC^{\ts n}\ts)$
in \eqref{vvp}. The latter action is not exactly that on
$\F_m(V)\ts$, but is defined by the composition
$$
\X(\g_n)
\to
\X(\g_n)
\underset{\be_m}\longrightarrow
\U(\ts\f_m)\ot\ts\PD\ts(\CC^{\ts m}\ot\CC^{\ts n})
$$
where the first map is the automorphism \eqref{fus} with 
$f(u)$ given by \eqref{fuss}. The second map
is the defining homomorphism of the $\X(\g_n)\ts$-module $\F_m(V)$.
This third property of $\F_m$ was the origin
of our definition of this functor.
Thus we have two different descriptions of
the same action of $\X(\g_n)$ on \eqref{vvp}.
Another two, still different descriptions of the same action of $\X(\g_n)$
on the vector space \eqref{vvp} were provided in
\cite{M} and \cite{N2}~respectively.

The functor $\F_m$ here is an
\lq\lq\ts antisymmetric\ts\rq\rq\
version of a functor introduced in \cite{KN3}.
The exterior algebra $\mathrm{\Lambda}\ts(\CC^{\ts m}\ot\CC^{\ts n})$ here
replaces the symmetric algebra
$\S\ts(\CC^{\ts m}\ot\CC^{\ts n})$ in \cite{KN3}.
Analogues of the three fundamental properties of
$\F_m$ were also given in \cite{KN3}.


\section*{\normalsize 1.\ Twisted Yangians}
\setcounter{section}{1}
\setcounter{equation}{0}
\setcounter{theorem*}{0}

Let $G_n$ be one of the complex Lie groups $O_n$ and $Sp_n\ts$.
We regard $G_n$ as the subgroup  of the general linear Lie group
$GL_n\ts$, preserving a non-degenerate bilinear form $\langle\ ,\,\rangle$
on the vector space $\CC^{\ts n}\ts$. This form is symmetric in the case
$G_n=O_n\ts$, and alternating in the case $G_n=Sp_n\ts$. In the latter case
$n$ has to be even. We always assume that the integer $n$ is positive.
Throughout this article, we will use the following convention. Whenever
the double sign $\ts\pm\ts$ or $\ts\mp\ts$ appears, the upper sign corresponds
to the case
$G_n=O_n$ while the lower sign
corresponds to the case
$G_n=Sp_n\ts$.

Let $i$ be any of the indices $1\lcd n\ts$. If $i$ is even, put
$\bi=i-1\ts$. If $i$ is odd and $i<n$, put $\bi=i+1\ts$. Finally,
if $i=n$ and $n$ is odd, put $\bi=i\ts$.
Let $e_1\lcd e_n$ be the vectors of the standard basis in $\CC^{\ts n}$.
Choose the bilinear form on $\CC^{\ts n}$ so that for any two basis vectors
$e_i$ and $e_j$ we have
$\langle\ts e_i\com e_j\ts\rangle=\th_i\,\de_{\ts\bi j}$ where
$\th_i=1$ or $\th_i=(-1)^{\ts i-1}$
in the case of the symmetric or alternating form.

Let $E_{\ts ij}\in\End(\CC^{\ts n})$ be the standard matrix units.
We will also regard these matrix units as basis elements of
the general linear Lie algebra $\gl_n\ts$.
Let $\g_n$ be the Lie algebra of the group $G_n\ts$,
so that $\g_n=\so_n$ or $\g_n=\sp_n$
in the case of the symmetric or alternating form on $\CC^{\ts n}$.
The Lie subalgebra $\g_n\subset\gl_n$ is spanned by the elements
$
E_{\ts ij}-\th_i\ts\th_j\ts E_{\ts\bj\ts\bi}\,.
$

Take the \textit{Yangian\/} $\Y(\gl_n)$ of the Lie algebra $\gl_n\ts$.
The unital associative algebra $\Y(\gl_n)$ over $\CC$
has a family of generators
$
T_{ij}^{\ts(1)},T_{ij}^{\ts(2)},\ts\ldots
$
where $i\com j=1\lcd n\ts$.
Defining relations for these generators
can be written using the series
\begin{equation*}
T_{ij}(u)=
\de_{ij}+T_{ij}^{\ts(1)}u^{-\ns1}+T_{ij}^{\ts(2)}u^{-\ns2}+\,\ldots
\end{equation*}
where $u$ is a formal parameter. Let $v$ be another formal parameter.
Then the defining relations in the associative algebra $\Y(\gl_n)$
can be written as
\begin{equation}
\label{yrel}
(u-v)\,[\ts T_{ij}(u)\ts,T_{kl}(v)\ts]\ts=\;
T_{kj}(u)\ts T_{il}(v)-T_{kj}(v)\ts T_{il}(u)\,.
\end{equation}

The algebra $\Y(\gl_n)$ is commutative if $n=1\ts$.
By \eqref{yrel}, for any $z\in\CC$
the assignments
\begin{equation}
\label{tauz}
\tau_z\ts:\,T_{ij}(u)\ts\mapsto\,T_{ij}(u-z)
\end{equation}
define an automorphism $\tau_z$ of the algebra $\Y(\gl_n)\ts$.
Here each of the formal
power series $T_{ij}(u-z)$ in $(u-z)^{-1}$ should be re-expanded
in $u^{-1}$, and every assignment \eqref{tauz} is a correspondence
between the respective coefficients of series in $u^{-1}$.
Relations \eqref{yrel} also show that for any
formal power series $g(u)$ in $u^{-1}$ with coefficients from
$\CC$ and leading term $1$, the assignments
\begin{equation}
\label{fut}
T_{ij}(u)\ts\mapsto\,g(u)\,T_{ij}(u)
\end{equation}
define an automorphism of the algebra $\Y(\gl_n)\ts$.
Using \eqref{yrel}, one can directly verify that the assignments
\begin{equation}
\label{eval}
T_{ij}(u)\ts\mapsto\,\de_{ij}+E_{ij}\ts u^{-1}
\end{equation}
define a homomorphism of unital associative algebras
$\Y(\gl_n)\to\U(\gl_n)\ts$.

There is an embedding $\U(\gl_n)\to\Y(\gl_n)\ts$,
defined by mapping $E_{ij}\mapsto T_{ij}^{\ts(1)}$. So
$\Y(\gl_n)$ contains the universal enveloping
algebra $\U(\gl_n)$ as a subalgebra. The homomorphism \eqref{eval} is
identical on the subalgebra $\U(\gl_n)\subset\Y(\gl_n)\ts$.

Let $T(u)$ be the $n\times n$ matrix
whose $i\com j$ entry is the series $T_{ij}(u)\ts$.
The relations \eqref{yrel} can be rewritten
by using the \textit{Yang R-matrix\/}.
This is the $n^2\times n^2$ matrix
\begin{equation}
\label{ru}
R(u)\,=\,u-\sum_{i,j=1}^n\,E_{ij}\ot E_{ji}
\end{equation}
where the tensor factors $E_{ij}$ and $E_{ji}$
are regarded as $n\times n$ matrices.
Note that
\begin{equation}
\label{rur}
R(u)\,R(-u)=1-u^2\ts.
\end{equation}
Take $n^2\times n^2$ matrices whose entries are
series with coefficients from $\Y(\gl_n)\ts$,
$$
T_1(u)=T(u)\ot1
\ \quad\text{and}\ \quad
T_2(v)=1\ot T(v)\,.
$$
The collection of relations \eqref{yrel} for all possible indices
$i\com j\com k\com l$ can be written as
\begin{equation}
\label{rtt}
R(u-v)\,T_1(u)\,T_2(v)\,=\,T_2(v)\,T_1(u)\,R(u-v)\,.
\end{equation}

Using this form of the defining relations together with \eqref{rur},
one shows that 
\begin{equation}
\label{tin}
T(u)\mapsto T(-u)^{-1}
\end{equation}
defines an involutive automorphism of the algebra $\Y(\gl_n)\ts$.
Here each entry of the inverse matrix $T(-u)^{-1}$
is a formal power series in $u^{-1}$ with coefficients
from the algebra $\Y(\gl_n)\ts$, and the assignment \eqref{tin}
is as a correspondence between the respective matrix entries.

The Yangian $\Y(\gl_n)$ is a Hopf algebra over the field $\CC\ts$.
The comultiplication $\De:\Y(\gl_n)\to\Y(\gl_n)\ot\Y(\gl_n)$ is defined by
the assignment
\begin{equation}\label{1.33}
\De:\,T_{ij}(u)\ts\mapsto\ts\sum_{k=1}^n\ T_{ik}(u)\ot T_{kj}(u)\,.
\end{equation}
When taking tensor products of $\Y(\gl_n)\ts$-modules,
we use the comultiplication \eqref{1.33}.
The counit homomorphism
$\Y(\gl_n)\to\CC$ is defined by the assignment
$
T_{ij}(u)\ts\mapsto\ts\de_{ij}\ts.
$
The antipodal map $\Y(\gl_n)\to\Y(\gl_n)$ is defined by the assignment
$
T(u)\mapsto T(u)^{-1}.
$
This map is an anti-automorphism of the associative algebra $\Y(\gl_n)\ts$.
For further details on the Hopf algebra
structure on $\Y(\gl_n)$ see \cite[Chapter 1]{MNO}.

Let $\Tp(u)$ be the transpose to the matrix $T(u)$
relative to the form $\langle\ ,\,\rangle$ on $\CC^{\ts n}\ts$.
The $i\com j$ entry of the matrix $\Tp(u)$ is
$\ts\th_i\ts\th_j\ts T_{\ts\bj\ts\bi\ts}(u)\ts$.
Define the $n^2\times n^2$ matrices
$$
\Tp_1(u)=\Tp(u)\ot1
\ \quad\text{and}\ \quad
\Tp_2(v)=1\ot\Tp(v)\,.
$$
Note that the Yang $R$-matrix \eqref{ru} is invariant under applying
the transposition relative to $\langle\ ,\,\rangle$ to both tensor factors.
Hence the relation \eqref{rtt} implies that
$$
\Tp_1(u)\,\Tp_2(v)\,R(u-v)
\,=\,
R(u-v)\,\Tp_2(v)\,\Tp_1(u)\,,
$$
\begin{equation}
\label{rttp}
R(u-v)\,\Tp_1(-u)\,\Tp_2(-v)
\,=\,
\Tp_2(-v)\,\Tp_1(-u)\,R(u-v)\,.
\vspace{4pt}
\end{equation}
To obtain the latter relation, we used \eqref{rur}.
By comparing \eqref{rtt} and \eqref{rttp},
an involutive automorphism of the algebra $\Y(\gl_n)$
can be defined by the assignment
\begin{equation}
\label{transauto}
T(u)\mapsto\Tp(-u)\ts.
\end{equation}
This assignments is
understood as a correspondence between respective matrix entries.

Now take the product $\Tp(-u)\,T(u)\ts$.
The $i\com j$ entry of this 
matrix is the series
\begin{equation}
\label{yser}
\sum_{k=1}^n\,\th_i\ts\th_k\,T_{\,\bk\ts\bi\ts}(-u)\,T_{kj}(u)\,.
\end{equation}
The \textit{twisted Yangian\/}
corresponding to the form $\langle\ ,\,\rangle$
is the subalgebra of $\Y(\gl_n)$ generated by
coefficients of all series \eqref{yser}.
We denote this subalgebra by $\Y(\g_n)\ts$.

To give defining relations for these generators of
$\Y(\g_n)\ts$, let us introduce the \textit{extended twisted Yangian}
$\X(\ts\g_n)\ts$. The unital associative algebra $\X(\ts\g_n)$ has a
family of generators
$
S_{ij}^{\ts(1)},S_{ij}^{\ts(2)},\ts\ldots
$
where $i\com j=1\lcd n\ts$.
Put
\begin{equation*}
S_{ij}(u)=
\de_{ij}+S_{ij}^{\ts(1)}u^{-\ns1}+S_{ij}^{\ts(2)}u^{-\ns2}+\,\ldots
\end{equation*}
and let $S(u)$ be the $n\times n$ matrix
whose $i\com j$ entry is the series $S_{ij}(u)\ts$.
Also introduce the $n^2\times n^2$ matrix
\begin{equation}
\label{rup}
\Rp(u)\,=\,u-\sum_{i,j=1}^n\,\th_i\,\th_j\,E_{ij}\ot E_{\ts\bi\ts\bj}
\end{equation}
which is obtained from the Yang $R$-matrix \eqref{ru} by applying
to any of the two tensor factors
the transposition relative to the form $\langle\ ,\,\rangle$ on $\CC^n\ts$.
Note the relation
\begin{equation}
\label{rurp}
\Rp(u)\,\Rp(\ts n-u)=u\ts(\ts n-u)\ts.
\end{equation}
Take $n^2\times n^2$ matrices whose entries are
series with coefficients from the algebra $\X(\g_n)\ts$,
$$
S_1(u)=S(u)\ot1
\ \quad\text{and}\ \quad
S_2(v)=1\ot S(v)\,.
$$
Defining relations in the algebra $\X(\g_n)$ can then be written
as a single matrix relation
\begin{equation}
\label{rsrs}
R(u-v)\,S_1(u)\,\Rp(-u-v)\,S_2(v)\,=\,S_2(v)\,\Rp(-u-v)\,S_1(u)\,R(u-v)\,.
\end{equation}
It is equivalent to the collection of relations
$$
(u^2-v^2)\,[\ts S_{ij}(u)\ts,S_{kl}(v)\ts]\ts=\ts
(u+v)(\ts S_{kj}(u)\,S_{il}(v)-S_{kj}(v)\,S_{il}(u))
$$
$$
\mp\,(u-v)\,(\,
\th_k\ts\th_j\,S_{i\ts\bk}(u)\,S_{\ts\bj\ts l}(v)-
\th_i\ts\th_l\,S_{k\ts\bi}(v)\,S_{\ts\bl\ts j}(u))
$$
\vspace{-8pt}
\begin{equation}
\label{xrel}
\pm\,\ts\th_i\ts\th_j\,
(\ts S_{k\ts\bi\ts}(u)\,S_{\ts\bj\ts l}(v)-
S_{k\ts\bi\ts}(v)\,S_{\ts\bj\ts l}(u))\,.
\vspace{4pt}
\end{equation}
Similarly to \eqref{fut}, this collection of
relations shows that for any
formal power series $f(u)$ in $u^{-1}$ with the coefficients from
$\CC$ and leading term $1$, the assignments
\begin{equation}
\label{fus}
S_{ij}(u)\ts\mapsto\,f(u)\,S_{ij}(u)
\end{equation}
define an automorphism of the algebra
$\X(\g_n)\ts$.
See \cite[Section 1]{KN3}
for the proof of 

\begin{proposition*}
\label{xyp}
One can define a homomorphism $\X(\g_n)\to\Y(\g_n)$ by assigning
\begin{equation}
\label{xy}
S(u)\,\mapsto\,\Tp(-u)\,T(u)\ts.
\end{equation}
\end{proposition*}

By definition, the homomorphism \eqref{xy} is surjective.
Further, the algebra $\X(\g_n)$ has a distinguished family
of central elements. Indeed,
by dividing each side of the equality \eqref{rsrs}
by $S_2(v)$ on the left and right and then setting $v=-u\ts$,
we get
$$
\Rp(0)\,S_1(u)\,R(2u)\,S_2(-u)^{-1}=\ts
S_2(-u)^{-1}\ts R(2u)\,S_1(u)\,\Rp(0)\,.
$$
The rank of the matrix $\Rp(0)$ equals $1$. So
the last displayed equality implies existence of a formal power series
$O(u)$ in $u^{-1}$ with the coefficients in $\X(\g_n)$ and leading term $1$,
such that
\begin{equation}
\label{5.2}
\Rp(0)\,S_1(u)\,R(2u)\,S_2(-u)^{-1}=\ts
(\ts2u\mp1\ts)\,O(u)\,\Rp(0)\,.
\end{equation}
By \cite[Theorem 6.3]{MNO} all coefficients of the series $O(u)$
belong to the centre of $\X(\g_n)\ts$.
Let us write
$$
O(u)\,=\,1+O^{(1)}u^{\ns-1}+O^{(2)}u^{-\ns2}+\,\ldots\,.
$$
By \cite[Theorem 6.4]{MNO}
the kernel of the homomorphism \eqref{xy}
coincides with the (two-sided) ideal generated by the central elements
$O^{\ts(1)}\ns\com O^{\ts(2)}\ns\com\,\ldots$
defined as coefficients of the series $O(u)\ts$. Using \eqref{rur},
one derives from \eqref{5.2} the relation $O(u)\ts O(-u)=1\ts$.

Thus the twisted Yangian $\Y(\g_n)$ can be defined as the
associative algebra with the generators
$
S_{ij}^{\ts(1)},S_{ij}^{\ts(2)},\ts\ldots
$
which satisfy the relation $O(u)=1$ and the
\textit{reflection equation\/} \eqref{rsrs}.
For more details on the definition of the algebra $\Y(\g_n)$
see \cite[Chapter~3]{MNO}.

In the present article we need the algebra $\X(\g_n)$ which is
determined by \eqref{rsrs} alone, because this algebra admits
an analogue of the automorphism \eqref{tin} of the Yangian $\Y(\gl_n)\ts$.
Indeed, using \eqref{rsrs} together with the relations
\eqref{rur} and \eqref{rurp}, one shows that the assignment
\begin{equation}
\label{sin}
\om_n:\ts S(u)\ts\mapsto{S(-\ts u-{n}/2\ts)}^{-1}
\end{equation}
defines an involutive automorphism $\om_n$ of $\X(\g_n)\ts$.
However, $\om_n$ does not determine an automorphism of the algebra
$\Y(\g_n)\ts$, because the map $\om_n$ does
not preserve the ideal of $\X(\g_n)$ generated by the elements
$O^{(1)}\ns\com O^{(2)}\ns\com\,\ldots\ts\,$;
see \cite[Section 6.6]{MNO}.
Note that by multiplying
\eqref{5.2} on the right by $S_2(-u)\ts$, the relation $O(u)=1$
can be rewritten as
\begin{equation}
\label{spu}
\Sp(u)\ts=\ts S(-u)\ts\pm\ts\frac{S(u)-S(-u)}{2u}
\end{equation}
where $\Sp(u)$ is the transpose to the matrix $S(u)$
relative to the form $\langle\ ,\,\rangle$ on $\CC^{\ts n}\ts$.

The definition \eqref{5.2} of the series $O(u)$
implies that the assignment \eqref{fus} determines an automorphism
of the quotient algebra $\Y(\g_n)$ of $\X(\g_n)\ts$,
if and only if $f(u)=f(-u)$.
If $z\neq0$, the automorphism $\tau_z$ of
$\Y(\gl_n)$ does not preserve the subalgebra
$\Y(\g_n)\subset\Y(\gl_n)\ts$. 
There is no analogue of the automorphism $\tau_z$
for the algebra $\X(\g_n)$.

However, there is an analogue of the
homomorphism $\Y(\gl_n)\to\U(\gl_n)$ defined by \eqref{eval}.
Namely, one can define a homomorphism $\pi_n:\X(\g_n)\to\U(\g_n)$
by the assignments
\begin{equation}
\label{pin}
\pi_n:\,S_{ij}(u)\,\mapsto\,\de_{ij}+
\frac{E_{ij}-\th_i\ts\th_j\ts E_{\ts\bj\ts\bi}}
{\textstyle u\pm\frac12}
\end{equation}
This can be proved by using the defining relations \eqref{xrel},
see \cite[Proposition~3.11]{MNO}.
Furthermore, the central elements
$O^{\ts(1)}\ns\com O^{\ts(2)}\ns\com\,\ldots$ of $\X(\g_n)$
belong to the kernel of $\pi_n\ts$.
Thus $\pi_n$ factors through the homomorphism $\X(\g_n)\to\Y(\g_n)$
defined by \eqref{xy}.

Further, there is an embedding $\U(\g_n)\to\Y(\g_n)$ defined by mapping
each element $E_{ij}-\,\th_i\ts\th_j\ts E_{\ts\bj\ts\bi}\in\g_n$
to the coefficient at $u^{-1}$ of the series \eqref{yser}. Hence
$\Y(\g_n)$ contains the universal enveloping
algebra $\U(\g_n)$ as a subalgebra.
The homomorphism $\Y(\g_n)\to\U(\g_n)$ corresponding to $\pi_n$
is evidently identical on the subalgebra $\U(\g_n)\subset\Y(\g_n)\ts$.

For any positive integer $l\ts$, consider
the vector space $\CC^{\ts l}$ and the corresponding Lie algebra
$\gl_{\ts l}\ts$. Let $E_{ab}\in\End(\CC^{\ts l})$ with
$a\com b=1\lcd l$ be the standard matrix units.
When regarding these matrix units as generators of
the universal enveloping algebra $\U(\gl_{\ts l})\ts$, introduce the
$l\times l$ matrix $E$ whose $a\com b$ entry is the
generator $E_{ab}\ts$. Denote by $\Ep$ the
$l\times l$ matrix whose $a\com b$ entry is the
generator $E_{\ts ba}\ts$. Then consider the matrix inverse
$(\ts u-\Ep\ts)^{-1}\ts$. The $a\com b$ entry
$(\ts u-\Ep\ts)^{-1}_{\ts ab}$ of the inverse matrix
is a formal power series in $u^{-1}$
with the leading term $\de_{ab}\,u^{-1}$ and the
coefficients from the algebra $\U(\gl_{\ts l})\ts$.

Take the tensor product of the vector spaces
$\CC^{\ts l}\ot\CC^{\ts n}$. Let $x_{ai}$ with
$a=1\lcd l$ and $i=1\lcd n$ be the standard coordinate functions
on $\CC^{\ts l}\ot\CC^{\ts n}\ts$.
Consider the \textit{Grassmann algebra\,}
$\P\ts(\CC^{\ts l}\ot\CC^{\ts n})\ts$.
It is generated by the elements $x_{ai}$
subject to the anticommutation relations
$
x_{ai}\,x_{bj}=-\ts x_{bj}\,x_{ai}
$
for all indices $a\com b=1\lcd l$ and $i\com j=1\lcd n\ts$.
We will denote the operator of the left multiplication
by $x_{ai}$ on $\P\ts(\CC^{\ts l}\ot\CC^{\ts n})$ by the same symbol.
Let $\d_{ai}$ be the operator of left derivation on
$\P\ts(\CC^{\ts l}\ot\CC^{\ts n})$
corresponding to the variable $x_{ai}\ts$, 
it is also called the \textit{inner multiplication}
in $\P\ts(\CC^{\ts l}\ot\CC^{\ts n})$ corresponding to 
$x_{ai}\ts$.

The ring of $\CC\ts$-endomorphisms of $\P\ts(\CC^{\ts l}\ot\CC^{\ts n})$
is generated by all operators
$x_{ai}$ and $\d_{ai}\ts$,
see for instance \cite[Appendix~2.3]{H}.
This ring will be denoted by
$\PD\ts(\CC^{\ts l}\ot\CC^{\ts n})$.
In this ring, we have 
\begin{equation}
\label{cliff}
x_{ai}\,\d_{bj}\ts+\,\d_{bj}\,x_{ai}\,=\,\de_{ab}\,\de_{ij}\,.
\end{equation}
Hence the ring $\PD\ts(\CC^{\ts l}\ot\CC^{\ts n})$
is isomorphic to the Clifford algebra
corresponding to the direct sum of the vector space
$\CC^{\ts l}\ot\CC^{\ts n}$ with its dual.

The Lie algebra $\gl_{\ts l}$ acts on the vector space
$\P\ts(\CC^{\ts l}\ot\CC^{\ts n})$ so that the generator
$E_{ab}$ acts as the operator
\begin{equation}
\label{glmact}
\sum_{k=1}^n\,
x_{ak}\,\d_{\ts bk}\ts.
\end{equation}
Denote by $\A_{\ts l}$ the tensor product of associative algebras
$\U(\gl_{\ts l})\ot\ts\PD\ts(\CC^{\ts l}\ot\CC^{\ts n})\ts$. We have an
embedding $\U(\gl_{\ts l})\to\A_{\ts l}$ defined for $a\com b=1\lcd l$
by the mappings
\begin{equation}
\label{eabact}
E_{ab}\mapsto
E_{ab}\ot1\,+\,
\sum_{k=1}^n\,
1\ot x_{ak}\,\d_{\ts bk}\ts.
\end{equation}
The following proposition was proved in \cite[Section 1]{KN2},
see also \cite[Section 3]{A}.

\begin{proposition*}
\label{dast}
{\rm\,\,(i)}
One can define a homomorphism $\al_{\ts l}:\Y(\gl_n)\to\A_{\ts l}$
by mapping
\begin{equation}
\label{ehom}
\al_{\ts l}:\
T_{ij}(u)\ts\mapsto\,\de_{ij}+
\sum_{a,b=1}^l(\ts u-\Ep\ts)^{-1}_{\ts ab}\ot x_{ai}\,\d_{\ts bj}
\end{equation}
{\rm(ii)}
The image of\/ $\Y(\gl_n)$ in\/ $\A_{\ts l}$ relative to this
homomorphism commutes with the image of\/ $\U(\gl_{\ts l})$ in\/ $\A_{\ts l}$
relative to the embedding \eqref{eabact}.
\end{proposition*}

Note that
$$
\al_{\ts l}:\
T_{ij}^{\ts(1)}\,\mapsto\,\sum_{c=1}^l\,1\ot x_{ci}\,\d_{cj}\ts.
$$
Hence the restriction of the homomorphism $\al_{\ts l}$
to the subalgebra $\U(\gl_n)\subset\Y(\gl_n)$ corresponds
to the natural action of the Lie algebra $\gl_n$ on
$\P\ts(\CC^{\ts l}\ot\CC^{\ts n})$.

Denote by $Z(u)$ the trace of the inverse matrix
$(u+E)^{-1}\ts$, so that
\begin{equation}
\label{zu}
Z(u)\,=\,\sum_{c=1}^l\,(\ts u+E\ts)^{-1}_{\ts cc}\,.
\end{equation}
Then $Z(u)$ is a formal power series in $u^{-1}$
with the coefficients from the algebra  $\U(\gl_{\ts l})\ts$.
It is well known that these coefficients actually belong to the
centre $\Z\ts(\gl_{\ts l})$ of $\U(\gl_{\ts l})\ts$.
Note that the leading term of this series is $l\ts u^{-1}$.

Let us choose the Borel subalgebra $\b$ of the Lie algebra $\gl_{\ts l}$
spanned by the elements $E_{ab}$ where $a\le b\ts$.
Let $\a\subset\b$ be the Cartan subalgebra of $\gl_{\ts l}$
with the basis $(E_{11}\lcd E_{ll})\ts$. Consider
the corresponding \textit{Harish-Chandra homomorphism\/}
$\varphi_{\ts l}:\ts\U(\gl_{\ts l})^\a\to\U(\a)\ts$. By definition,
for any $\a\ts$-invariant element $X\in\U(\gl_{\ts l})$ the difference
$X-\ts\varphi_{\ts l}\ts(X)$ belongs to the left ideal of
$\U(\gl_{\ts l})$ generated by the elements $E_{ab}$ where $a<b\ts$.
Restriction of the homomorphism $\varphi_{\ts l}$ to
$\Z\ts(\gl_{\ts l})\subset\U(\gl_{\ts l})^\a$ is injective.
It is well known that
\begin{equation}
\label{hczu}
1+\varphi_{\ts l}\ts(Z(u))\,=\,\prod_{a=1}^l\,\ts
\Bigl(1+\frac1{u+l-a+E_{aa}}\ts\Bigr)\,,
\end{equation}
see for instance \cite[Theorem 3]{PP}.
For the proof of the next lemma see \cite[Section~1]{KN3}
where the parameter $u$ should be now replaced by $-u\ts$.

\begin{lemma*}
\label{eep}
For any indices\/ $a\com d=1\lcd l$ we have the equality
$$
(u+E\ts)^{-1}_{\ts da}=(1+Z(u))\,(\ts u+l+\Ep\ts)^{-1}_{\ts ad}\,.
$$
\end{lemma*}

Now let $U$ be a module of the Lie algebra $\gl_{\ts l}\ts$.
Using the homomorphism \eqref{ehom} we can turn the tensor
product of $\gl_{\ts l}\ts$-modules $U\ot\P\ts(\CC^{\ts l}\ot\CC^{\ts n})$
to a bimodule over $\gl_{\ts l}$ and $\Y(\gl_n)\ts$.
This bimodule is denoted by $\E_{\ts l}\ts(U)\ts$.
More generally, for $z\in\CC$ denote by
$\E_{\ts l}^{\ts z}\ts(U)$ the $\Y(\gl_n)\ts$-module obtained
from  $\E_{\ts l}\ts(U)$ via pull-back through the automorphism
$\tau_{-z}$ of $\Y(\gl_n)\ts$, see \eqref{tauz}. It is determined by
the homomorphism $\Y(\gl_n)\to\A_{\ts l}$ such that
$$
T_{ij}(u)\ts\mapsto\,\de_{ij}\,+
\sum_{a,b=1}^l(\ts u+z-\Ep\ts)^{-1}_{\ts ab}\ot x_{ai}\,\d_{\ts bj}
$$
for any $i\com j=1\lcd n\ts$.
As a $\gl_{\ts l}\ts$-module $\E_{\ts l}^{\ts z}\ts(U)$
coincides with $\E_{\ts l}\ts(U)$ by definition.
In the next section we will introduce analogues
of the homomorphism \eqref{eabact}
and of the correspondence $U\mapsto\E_{\ts l}\ts(U)$ for
the twisted Yangian $\Y(\g_n)$ instead of $\Y(\gl_n)\ts$.


\section*{\normalsize 2. Howe duality}
\setcounter{section}{2}
\setcounter{equation}{0}
\setcounter{theorem*}{0}

We will work with one of
the pairs $(\so_{2m}\com O_n)$ and $(\sp_{2m}\com Sp_n)\ts$.
The second member of the pair will be the Lie group $G_n\ts$.
The first member will be the Lie algebra $\f_m$ defined below.
These pairs appear in the context of the skew Howe duality,
see \cite[Section 4.3]{H}.

Take the even-dimensional vector space $\CC^{\ts 2m}$.
Equip $\CC^{\ts 2m}$ with a non-degenerate bilinear form,
symmetric in the case $G_n=O_n$ and alternating in the case $G_n=Sp_n\ts$.
Let $\f_m$ be the subalgebra of the general Lie algebra
$\gl_{\ts 2m}$ preserving our bilinear form on $\CC^{\ts 2m}$.
We have $\f_m=\so_{2m}$ or $\f_m=\sp_{2m}$ respectively
in the case of a symmetric or an alternating form on $\CC^{\ts 2m}$.

Let us label the standard basis vectors of $\CC^{\ts 2m}$ by the
numbers $-\ts m\lcd-1\com1\lcd m\ts$.
Let $E_{ab}\in\End(\CC^{\ts2m})$ be the standard matrix units,
where the indices $a\com b$ run through these numbers.
We will also regard these matrix units as basis elements of 
$\gl_{\ts2m}\ts$. Put
\begin{equation}
\label{eab}
\ep_{ab}=1
\ \quad\text{or}\ \quad
\ep_{ab}=\operatorname{sign}a\cdot\operatorname{sign}b
\end{equation}
respectively in the case of a symmetric or an alternating form on
$\CC^{\ts 2m}$.
Then choose the form on $\CC^{\ts2m}$ so that
the Lie subalgebra $\f_m\subset\gl_{\ts2m}$ is spanned
by the elements
\begin{equation}
\label{fab}
F_{ab}=E_{ab}-\ep_{ab}\,E_{-b,-a}\,.
\end{equation}
In the universal enveloping algebra $\U(\ts\f_m)$
we have the commutation relations
\begin{equation}
\label{ufmrel}
[\ts F_{ab}\com F_{cd}\ts]=
\de_{cb}\,F_{ad}-
\de_{ad}\,F_{cb}-
\ep_{ab}\,\de_{c,-a}\,F_{-b,d}+
\ep_{ab}\,\de_{-b,d}\,F_{c,-a}\,.
\end{equation}

Let $F$ be the $2\ts m\times2\ts m$ matrix
whose $a\com b$ entry is the element $F_{ab}\ts$.
Denote by $F(u)$ the inverse to the matrix $u+F\ts$.
Let $F_{ab}(u)$ be the $a\com b$ entry of the inverse matrix.
Any of these entries may be regarded as a formal power series in $u^{-1}$
with the coefficients from the algebra $\U(\ts\f_m)\ts$. Then
\begin{equation}
\label{fabu}
F_{ab}(u)=\ts\de_{ab}\,u^{-1}\,+\,
\sum_{s=0}^\infty\,
\sum_{\ |c_1|,\ldots,|c_s|=1}^m
(-1)^{\ts s+1}\,
F_{ac_1}\ts F_{c_1c_2}\ldots\ts F_{c_{s-1}c_s}\ts F_{c_sb}\,\,
u^{-s-2}\,.
\end{equation}
When $s=0\ts$, the sum over
$c_1\lcd c_s$ in \eqref{fabu} is understood as 
$-\ts F_{ab}\,u^{-2}\ts$.
Let us denote by $W(u)$ the trace of the matrix $F(u)\ts$, that is
\begin{equation}
\label{wu}
W(u)=\sum_{|c|=1}^m F_{cc}(u)\,.
\end{equation}
The coefficients of the series $W(u)$ belong to the centre $\Z\ts(\f_m)$
of the algebra $\U(\ts\f_m)\ts$.

In what follows, the upper signs in $\ts\pm\ts$ and $\ts\mp\ts$
correspond to the case of a symmetric form on $\CC^{\ts 2m}$
while the lower signs correspond to the case of an alternating form
on $\CC^{\ts 2m}$. In these cases
we also have respectively a symmetric or alternating form
on $\CC^{\ts n}$. Thus the choice of signs in
$\ts\pm\ts$ and $\ts\mp\ts$ here agrees with our general convention
on double signs.
Let $\Fp(u)$ be the transpose to the matrix $F(u)$
relative to our 
bilinear form on $\CC^{\ts2m}$ so that
the $a\com b$ entry $\Fp_{ab}(u)$ of the matrix $\Fp(u)$ equals
$\ts\ep_{ab}\,F_{-b,-a}(u)\ts$. For the proof of next proposition
see \cite[Section 2]{KN3} where $u$ should be now replaced by $-\ts u\ts$.

\begin{proposition*}
\label{1.2}
We have the equality of\/ $2m\times2m$ matrices
$$
-\ts\Fp(u)\,=\,\bigl(\ts W(u)\ts\mp\ts\frac1{2u+2m\mp1}+1\tts\bigr)
\,F(\tts-\ts u-2m\pm1)\,\pm\,\frac{F(u)}{2u+2m\mp1}\,\,.
\vspace{4pt}
$$
\end{proposition*}

\begin{corollary*}
\label{1.3}
We have the equality
$$
\bigl(\ts W(u)\ts\mp\ts\frac1{2u+2m\mp1}+1\tts\bigr)
\bigl(\ts W(\tts -\ts u-2m\pm1\tts)\ts\pm\ts\frac1{2u+2m\mp1}+1\tts\bigr)
$$
$$
=1-\frac1{(2u+2m\mp1)^2}\,\,.
$$
\end{corollary*}

On the space $\CC^{\ts m}\ot\CC^{\ts n}\ts$,
we have the coordinate functions
$x_{ai}$ where $a=1\lcd m$ and $i=1\lcd n\ts$.
Consider the Grassmann algebra $\P\ts(\CC^{\ts m}\ot\CC^{\ts n})$
corresponding to this vector space.
We will denote the operator of the left multiplication
by $x_{ai}$ on $\P\ts(\CC^{\ts m}\ot\CC^{\ts n})$ by the same symbol.
Let $\d_{ai}$ be the left derivation on
$\P\ts(\CC^{\ts m}\ot\CC^{\ts n})$ relative to $x_{ai}\ts$.
There is an action of 
$\f_m$ on $\P\ts(\CC^{\ts m}\ot\CC^{\ts n})\ts$,
commuting with the natural action of the group
$G_n\ts$. The corresponding homomorphism
$\zeta_{\ts n}:\U(\ts\f_m)\to\PD\ts(\CC^{\ts m}\ot\CC^{\ts n})$
is defined by the following mappings for $a\com b=1\lcd m\ts$:
$$
\zeta_{\ts n}:\ F_{ab}\,\mapsto\,
-\,\de_{ab}\,n/2\,\,+\,
\sum_{k=1}^n\,x_{ak}\,\d_{\ts bk}\,,
$$
\begin{equation}
\label{gan}
F_{a,-b}\,\mapsto\,\sum_{k=1}^n\,
\th_k\,x_{a\bk}\,x_{\ts bk}\,,
\ \quad
F_{-a,b}\,\mapsto\,\sum_{k=1}^n\,
\th_k\,\d_{ak}\,\d_{\ts b\bk}\,.
\end{equation}
The homomorphism property here
can be verified by using the relations \eqref{ufmrel}.
Moreover, the image of the homomorphism $\zeta_{\ts n}$ coincides with the
subring of all $G_n$-invariants in $\PD\ts(\CC^{\ts m}\ot\CC^{\ts n})\ts$;
see \cite[Subsections 3.8.7 and 4.3.3]{H}.
Let $\B_m$ be the tensor product of associative algebras
$\U(\ts\f_m)\ot\ts\PD\ts(\CC^{\ts m}\ot\CC^{\ts n})\ts$. Take the
embedding $\U(\ts\f_m)\to\B_m\ts$ defined by mapping
\begin{equation}
\label{xact}
X\,\mapsto\,
X\ot1\,+\,
1\ot\zeta_{\ts n}\ts(X)
\ \quad\text{for each}\ \quad
X\in\f_m\,.
\end{equation}

\begin{proposition*}
\label{xb}
{\rm\,\,(i)}
One can define a homomorphism\/ $\be_m:\X(\g_n)\to\B_m$
so that the series $S_{ij}(u)$ is mapped to
the series with coefficients in the algebra $\B_m\ts$,
\begin{gather}
\nonumber
\de_{ij}\,+
\sum_{a,b=1}^m\,(\,
F_{-a,-b}\ts(u\pm{\textstyle\frac12}-m)
\ot\ts x_{ai}\,\d_{\ts bj}
\ts+\ts
F_{-a,b}\ts(u\pm{\textstyle\frac12}-m)
\ot\ts\th_j\,x_{ai}\,x_{\ts b\bj}
\\
\label{fhom}
+\ts\,
F_{a,-b}\ts(u\pm{\textstyle\frac12}-m)
\ot\ts\th_i\,\d_{a\bi}\,\d_{\ts bj}\ts
+\ts
F_{ab}\ts(u\pm\textstyle{\frac12}-m)
\ot\ts\th_i\,\th_j\,\d_{a\bi}\,x_{\ts b\bj}\,)\,.
\\
\nonumber
\end{gather}
{\rm(ii)}
The image of\/ $\X(\g_n)$ in\/ $\B_m$ relative to this homomorphism commutes
with the image of\/ $\U(\ts\f_m)$ in\/ $\B_m$ relative to the embedding
\eqref{xact}.
\end{proposition*}

Proposition \ref{xb} can be proved by direct calculation
using the defining relations \eqref{xrel}. That calculation
is omitted here. In Section 6
we will give a more conceptual proof of the proposition.
Now let the indices $c$ and $d$ run through the sequence
$-\ts m\lcd-1\com1\lcd m\ts$. For $c<0$ put
$p_{\ts ci}=x_{\ts-c,i}$ and $q_{\ts ci}=\d_{\ts-c,i}\ts$.
For $c>0$ put
$p_{\ts ci}=\th_i\,\d_{c\ts\bi}$
and
$q_{\ts ci}=\th_i\,x_{c\ts\bi}\,$.
Then our definition of the homomorphism $\be_m$ can be written as
\begin{equation}
\label{fhompq}
\be_m:\,S_{ij}(u)\,\mapsto\,
\de_{ij}\,+
\sum_{|c|,|d|=1}^m\,
F_{cd}\ts(u\pm{\textstyle\frac12}-m)\ot p_{\ts ci}\,q_{\ts dj}\,,
\end{equation}
similarly to \eqref{ehom}. Moreover, then by the definition \eqref{gan}
\begin{equation}
\label{ganpq}
\zeta_{\ts n}:\ F_{cd}\,\mapsto\,
-\,\de_{cd}\,{n}/2\,\,+\,
\sum_{k=1}^n\,q_{ck}\,p_{\ts dk}\,.
\end{equation}

Using \eqref{wu}, let us define a formal power series $\Wb(u)$ in $u^{-1}$
with coefficiens in the centre $\Z\ts(\f_m)$ of the algebra $\U(\ts\f_m)$
by the equation
$$
\bigl(\ts1\mp\frac1{2u}\ts\bigr)\,\,\Wb(u)\,=\,
\textstyle
W(u\pm\frac12-m)\,.
$$
By Corollary \ref{1.3},
$$
(\ts\Wb(u)+1\ts))\,(\ts\Wb(-u)+1\ts)=1\,.
$$
Hence there is a formal power series $\Wt(u)$ in $u^{-1}$
with coefficiens in $\Z\ts(\f_m)$ and leading term $1\ts$, such that
\begin{equation}
\label{wtu}
\Wt(-u)\,\ts\Wt(u)^{-1}=1+\Wb(u)\,.
\end{equation}
The series $\Wt(u)$ is not unique. But its coefficient
at $u^{-1}$ is always $-\ts m\ts$, because the leading term
of the series $W(u)$ is $2m\ts u^{-1}$.
Let $\bet_m$ be the homomorphism
$\X(\g_n)\to\B_m$ defined by assigning to
$S_{ij}(u)$ the series \eqref{fhom} multiplied by
\begin{equation}
\label{wub}
\Wt(u)\ot1\ts\in\ts\B_m\ts[[u^{-1}]]\ts.
\end{equation}
The homomorphism property of $\bet_m$ follows from Part~(i) of
Proposition \ref{xb}, see also the defining relations \eqref{xrel}.
Part (ii) implies that
the image of $\bet_m$ commutes with the image of\/ $\U(\ts\f_m)$ in
the algebra\/ $\B_m$ relative to the embedding
\eqref{xact}.

\begin{proposition*}
\label{tb}
The elements
$\,O^{(1)}\ns\com O^{(2)}\ns\com\,\ldots\,$ of\/ $\X(\g_n)$
belong to the kernel of\/ $\bet_m\ts$.
\end{proposition*}

\begin{proof}
Let us denote by $\St_{ij}(u)$ the product of the
series \eqref{fhom} and \eqref{wub}. Using the equivalent presentation
\eqref{spu} of the relation $O(u)=1$, we have to prove the equality
\begin{equation}
\label{stu}
\th_i\,\th_j\,\St_{\ts\bj\ts\bi\ts}(u)
\ts=\ts
\St_{ij}(-u)\ts\pm\ts\frac{\St_{ij}(u)-\St_{ij}(-u)}{2u}
\end{equation}
for any $i\com j=1\lcd n\ts$.
By the definition of the series $\Wt(u)\ts$, we have the relation
\begin{equation}
\label{wtw}
\Wt(u)\,(\ts1+W(u\pm{\textstyle\frac12}-m)\ts)
\ts=\ts
\Wt(-u)\ts\pm\ts\frac{\Wt(u)-\Wt(-u)}{2u}\ .
\end{equation}
Further, let us introduce the $2m\times 2m$ matrix
\begin{equation}
\label{ftud}
\Ft(u)\ts=\ts\Wt(u)\,F(u\pm{\textstyle\frac12}-m)
\end{equation}
and its transpose $\Ftp(u)$ relative to our bilinear form
on $\CC^{\ts2m}$. By Proposition~\ref{1.2},
\begin{equation}
\label{ftu}
\Ftp(u)
\ts=\ts
-\ts\Ft(-u)\ts\mp\ts\frac{\Ft(u)-\Ft(-u)}{2u}\ .
\end{equation}

By changing the indices
$i\com j$ in \eqref{fhom} respectively
to $\bj\ts\com\bi$ and multiplying
the resulting series by $\th_i\,\th_j$ we get
$$
\de_{ij}\,+
\sum_{a,b=1}^m\,(\,
F_{-a,-b}\ts(u\pm{\textstyle\frac12}-m)
\ot\ts\th_i\,\th_j\,x_{a\bj}\,\d_{\ts b\ts\bi}
\,\pm\,F_{-a,b}\ts(u\pm{\textstyle\frac12}-m)
\ot\ts\th_j\,x_{a\bj}\,x_{\ts bi}
$$
$$
\pm\,\,
F_{a,-b}\ts(u\pm{\textstyle\frac12}-m)
\ot\ts\th_i\,\d_{\ts aj}\,\d_{\ts b\ts\bi}
\ts+\ts
F_{ab}\ts(u\pm{\textstyle\frac12}-m)
\ot\ts\d_{\ts aj}\,x_{\ts bi}\,)
\vspace{12pt}
$$
$$
=\ (\ts1+W(u\pm{\textstyle\frac12}+m)\ts)\ot\de_{ij}\ +
$$
$$
\sum_{a,b=1}^m\,(\,
-\ts F_{\ts -b,-a}\ts(u\pm{\textstyle\frac12}-m)
\ot\ts\th_i\,\th_j\,\d_{\ts a\bi}\,x_{\ts b\bj}
\,\mp\,F_{\ts -b,a}\ts(u\pm{\textstyle\frac12}-m)
\ot\ts\th_j\,x_{ai}\,x_{\ts b\bj}
$$
$$
\mp\,\,
F_{\ts b,-a}\ts(u\pm{\textstyle\frac12}-m)
\ot\ts\th_i\,\d_{\ts a\bi}\,\d_{\ts bj}
\ts-\ts
F_{\ts ba}\ts(u\pm{\textstyle\frac12}-m)
\ot\ts x_{ai}\,\d_{\ts bj}\,)
\vspace{12pt}
$$
$$
=\ (\ts1+W(u\pm{\textstyle\frac12}-m)\ts)\ot\de_{ij}
$$
$$
-\,\sum_{a,b=1}^m\,(\,
\Fp_{ab}\ts(u\pm{\textstyle\frac12}-m)
\ot\ts\th_i\,\th_j\,\d_{\ts a\bi}\,\ts x_{\ts b\bj}
\ts+\ts
\Fp_{-a,b}\ts(u\pm{\textstyle\frac12}-m)
\ot\ts\th_j\,x_{ai}\,x_{\ts b\bj}\ +
$$
$$
\Fp_{a,-b}\ts(u\pm{\textstyle\frac12}-m)
\ot\ts\th_i\,\d_{\ts a\bi}\,\d_{\ts bj}
\ts+\ts
\Fp_{-a,-b}\ts(u\pm{\textstyle\frac12}-m)
\ot\ts x_{ai}\,\d_{\ts bj}\,)\,.
\vspace{8pt}
$$
Multiplying the expession in the last three lines by
$\Wt(u)\ot1$ and using the definition \eqref{ftud}, we get
$$
\Wt(u)\,(\ts1+W(u\pm{\textstyle\frac12}-m)\ts)\ot\de_{ij}
\vspace{4pt}
$$
$$
-\,\sum_{a,b=1}^m\,(\,
\Ftp_{ab}\ts(u)
\ot\ts\th_i\,\th_j\,\d_{\ts a\bi}\,\ts x_{\ts b\bj}
\ts+\ts
\Ftp_{-a,b}\ts(u)
\ot\ts\th_j\,x_{ai}\,x_{\ts b\bj}\ +
$$
$$
\Ftp_{a,-b}\ts(u)
\ot\ts\th_i\,\d_{\ts a\bi}\,\d_{\ts bj}
\ts+\ts
\Ftp_{-a,-b}\ts(u)
\ot\ts x_{ai}\,\d_{\ts bj}\,)\,.
\vspace{8pt}
$$
The required equality \eqref{stu} now follows from \eqref{wtw} and \eqref{ftu}.
\qed
\end{proof}

So the homomorphism $\bet_m:\X(\g_n)\to\B_m$
factors to a homomorphism $\Y(\g_n)\to\B_m\ts$.
This is an analogue of the homomorphism \eqref{ehom}
for the twisted Yangian $\Y(\g_n)$ instead of $\Y(\gl_n)$.
Recall that
$$
\Wt(u)=1-m\ts u^{-1}+\ldots
$$
so that
$$
\bet_m:\,S_{ij}^{\ts(1)}\mapsto\,-\ts m\,\de_{ij}+\,\sum_{c=1}^m\,\,
(\ts
1\ot\ts x_{ci}\,\d_{\ts cj}
+
1\ot\ts\th_i\,\th_j\,\d_{\ts c\ts\bi}\,x_{c\ts\bj}\ts)
$$
$$
=\,\sum_{c=1}^m\,\,
(\ts
1\ot\ts x_{ci}\,\d_{\ts cj}
-
1\ot\ts\th_i\,\th_j\,x_{c\ts\bj}\,\d_{\ts c\ts\bi}\ts)\,.
$$
Thus for any formal power series $\Wt(u)$
in $u^{-1}$ which has its coefficiens in $\Z\ts(\f_m)$,
has the leading term $1\ts$ and satisfies the equation \eqref{wtu},
the restriction of the homomorphism $\Y(\g_n)\to\B_m\ts$
to the subalgebra $\U(\g_n)\subset\Y(\g_n)$ corresponds
to the natural action of the Lie algebra $\g_n$ on
the vector space $\P\ts(\CC^{\ts m}\ot\CC^{\ts n})$.

The series $\Wt(u)$ is not unique, and it will be more convenient
for us to work with the homomorphism $\be_m:\X(\g_n)\to\B_m$
defined in Proposition \ref{xb}. Using this homomorphism
and the action of the Lie algebra $\f_m$ on $\P\ts(\CC^{\ts m}\ot\CC^{\ts n})$
as defined by \eqref{gan}, for arbitrary $\f_m$-module $V$
we can turn the tensor product $V\ot\P\ts(\CC^{\ts m}\ot\CC^{\ts n})$
to a bimodule over $\f_m$ and $\X(\g_n)\ts$.
This bimodule will be denoted by $\F_m(V)\ts$.

Consider the \textit{triangular decomposition}
of the Lie algebra $\f_m\ts$,
\begin{equation}
\label{tridec}
\f_m=\n\op\h\op\np
\end{equation}
where $\h$ is the Cartan subalgebra of $\f_m$ with the basis
$(\ts F_{-m,-m}\lcd F_{-1,-1})\ts$. Further,
$\n$ and $\np$ are the nilpotent
subalgebras of $\f_m$ spanned by elements $F_{ab}$ where
$a>b$ and $a<b$ respectively; the indices $a\com b$ here
can be positive or negative.
For each $\f_m$-module $V$, we will denote by $V_{\ts\n}$ the vector space
$V\tts/\ts\n\cdot V$
of coinvariants of the action of the subalgebra $\n\subset\f_m$ on $V$.
The Cartan subalgebra $\h\subset\f_m$ acts on
the vector space $V_{\ts\n}\ts$.

Now consider
the bimodule $\F_m(V)\ts$.
The action of $\X(\g_n)$ on this bimodule commutes with the
action of the Lie algebra $\f_m\ts$,
and hence with the action of the subalgebra $\n\subset\f_m\ts$.
Therefore the space $\F_m(V)_{\ts\n}$
of coinvariants of the action of $\n$
is a quotient of the $\X(\g_n)$-module $\F_m(V)\ts$.
Thus we get a functor from the category of all
$\f_m$-modules to the category of bimodules over $\h$ and $\X(\g_n)\ts$,
\begin{equation}
\label{zelefun}
V\mapsto\ts\F_m(V)_{\ts\n}=
(\ts V\ot\P\ts(\CC^{\ts m}\ot\CC^{\ts n})\ts)_{\ts\n}\ts.
\end{equation}

The assignments $E_{ab}\mapsto F_{ab}$ for all $a\com b=1\lcd m$
define a Lie algebra embedding $\gl_m\ns\to\f_m\ts$; see relations
\eqref{ufmrel}. Using this embedding, consider the decomposition
\begin{equation}
\label{pardec}
\f_m=\ts\r\ts\op\ts\gl_m\nns\op\rp
\end{equation}
where $\r$ and $\rp$ are the Abelian subalgebras of $\f_m$
spanned respectively by the elements $F_{a,-b}$ and
$F_{-a,b}$ for all $a\com b=1\lcd m\ts$.
For any $\gl_m$-module $U$, let $V$ be the $\f_m$-module
\textit{parabolically induced\/} from the $\gl_m$-module $U$.
To define $V$, one first extends the action of the Lie algebra
$\gl_m$ on $U$ to the maximal parabolic
subalgebra $\gl_m\op\ts\rp\subset\f_m\ts$, so that
every element of the summand $\rp$ acts on $U$ as zero.
By definition, $V$ is the $\f_m$-module induced from the
$\gl_m\op\ts\rp\ts$-module $U\ts$. Note that here we have
a canonical embedding $U\to V$ of
$\gl_m\op\ts\rp\ts$-modules\ts;
we will denote by $\uo$ the image of an element $u\in U$ under this
embedding. The $\f_m\ts$-module $V$ determines
the bimodule $\F_m(V)$ over $\f_m$ and $\X(\g_n)\ts$.
The space of $\r\ts$-coinvariants $\F_m(V)_{\ts\r}$ is
then a bimodule over $\gl_m$ and $\X(\g_n)\ts$.

On the other hand, for any $z\in\CC$ consider
the bimodule $\E_m^{\ts z}\ts(U)$ over the Lie algebra
$\gl_m$ and over the Yangian $\Y(\gl_n)\ts$.
By restricting the module $\E_m^{\ts z}\ts(U)$ from the algebra
$\Y(\gl_n)$ to its subalgebra $\Y(\g_n)$ and then by using the homomorphism
$\X(\g_n)\to\Y(\g_n)$ defined by \eqref{xy}, we can regard
$\E_m^{\ts z}\ts(U)$ as a module over the algebra $\X(\g_n)$
instead of $\Y(\gl_n)\ts$. This module is
determined by the homomorphism $\X(\g_n)\to\A_m$ such that
for any  $i\com j=1\lcd n$ the series $S_{ij}(u)$ is mapped to 
\begin{equation}
\label{rezhom}
\sum_{k=1}^n\,\th_i\,\th_k\,
\al_m\bigl(\ts T_{\ts\bk\ts\bi}(-\,u+z)\,T_{kj}(u+z)\bigr)\,;
\end{equation}
see \eqref{yser} and \eqref{ehom}.
Let us now map 
$S_{ij}(u)$ to the series \eqref{rezhom} multiplied~by
\begin{equation}
\label{difz}
(1+Z\ts(\ts u-z-m\ts))\ot1\ts\in\ts\A_m\ts[[u^{-1}]]\ts;
\end{equation}
see \eqref{zu} where the positive integer $l$
has to be now replaced by $m\ts$.
The latter mapping determines another homomorphism
$\X(\g_n)\to\A_m\ts$. Using it, we turn the
vector space $U\ot\P\ts(\CC^{\ts m}\ot\CC^{\ts n})$
of the $\X(\g_n)\ts$-module $\E_m^{\ts z}\ts(U)$ to
another $\X(\g_n)\ts$-module,
to be denoted by $\Eb_m^{\ts z}\ts(U)\ts$.
Further, define an action of the Lie algebra
$\gl_m$ on $\,\Eb_m^{\ts z}\ts(U)$ by pulling its action on
$\E_m^{\ts z}(U)$ back through the automorphism
\begin{equation}
\label{autom}
E_{ab}\,\mapsto\,-\,\de_{ab}\,{n}/2\,+\,E_{ab}
\quad\text{for}\quad
a\com b=1\lcd m\ts.
\end{equation}
Thus the action of $\gl_m$ on $\,\Eb_m^{\ts z}\ts(U)$
is determined by the composition of homomorphisms
$$
\U(\gl_m)\to\U(\gl_m)\to
\End(\ts U\ot\P\ts(\CC^{\ts m}\ot\CC^{\ts n}))\,,
$$
where the first map is the automorphism \eqref{autom}
while the second map corresponds to the natural action of $\gl_m$
on $\E_m^{\ts z}(U)\ts$. The following proposition is a
particular case of Theorem \ref{parind} from the next section.

\begin{proposition*}
\label{parres}
For the\/ $\f_m$-module $V$ parabolically induced from any\/
$\gl_m$-module $U$, the bimodule $\F_m(V)_{\ts\r}$
of\/ $\gl_m$ and\/ $\X(\g_n)$ is equivalent to\/ $\Eb_m^{\ts z}\ts(U)$
where $z=\mp\tts\frac12\ts$.
\end{proposition*}

Further, let $u$ and $f$ range over the vector spaces
$U$ and $\P\ts(\CC^{\ts m}\ot\CC^{\ts n})$ respectively.
In the next section, we will show that the linear map
$$
U\ot\P\ts(\CC^{\ts m}\ot\CC^{\ts n})
\to
(\ts V\ot\P\ts(\CC^{\ts m}\ot\CC^{\ts n}))_{\ts\r}
$$
defined by mapping $u\ot f$ to the class of $\,\uo\ot f$ in
the space of $\r\ts$-coinvariants, is an equivalence of bimodules
$\Eb_m^{\ts z}\ts(U)\to\F_m(V)_{\ts\r}$ over $\gl_m$ and~$\X(\g_n)$.

An element $\mu$ of the vector space $\h^\ast$ dual to
$\h$ is called a \textit{weight\/}.
A weight $\mu$ can be identified with
the sequence $(\ts\mu_1\lcd\mu_m\ts)$ of its \textit{labels},
where $$
\mu_a=\mu\ts(F_{\ts a-m-1,a-m-1})=-\ts\mu\ts(F_{\ts m-a+1,m-a+1})
\ \quad\textrm{for}\ \quad
a=1\lcd m\ts.
$$
The \textit{Verma module\/} $M_{\ts\mu}$ of the Lie algebra
$\f_m$ is the quotient of the algebra
$\U(\ts\f_m)$ by the left ideal generated by all elements
$X\in\np$, and by all elements $X-\mu\ts(X)$ where $X\in\h\ts$.
The elements of the Lie algebra $\f_m$ act on this quotient
via left multiplication. The image of the identity element $1\in\U(\ts\f_m)$
in this quotient is denoted by $1_\mu\ts$. Then
$X\cdot 1_\mu=0$ for all $X\in\np$, and
$X\cdot 1_\mu=\mu\ts(X)\cdot 1_\mu$
for all $X\in\h\ts$.
Denote by $L_\mu$ be the quotient of the Verma module $M_\mu$
relative to the maximal proper submodule. This quotient is a
simple $\f_m$-module of the highest weight $\mu\ts$.

For $z\in\CC$ denote by $P_z$
the $\Y(\gl_n)\ts$-module obtained by pulling the standard action
of $\U(\gl_n)$ on $\P(\CC^{\ts n})$ back through
the homomorphism $\Y(\gl_n)\to\U(\gl_n)$ defined by \eqref{eval},
and then back through the automorphism $\tau_{-z}$ of $\Y(\gl_n)\ts$.
Let $x_1\lcd x_n$ be the standard generators of $\P(\CC^{\ts n})$
and let $\d_1\lcd\d_n$ be the corresponding left derivations.
Using \eqref{diag3}, the action of $\Y(\gl_n)$ on $P_z$
is determined by the homomorphism $\Y(\gl_n)\to\PD\ts(\CC^{\ts n})\ts$,
\begin{equation}
\label{pz}
T_{ij}(u)\mapsto\de_{ij}+\frac{\,x_i\ts\d_j}{u+z}\ .
\end{equation}

Using the comultiplication \eqref{1.33},
for any $z_1\lcd z_m\in\CC$ define the tensor product of
$\Y(\gl_n)\ts$-modules
\begin{equation}
\label{pzm}
P_{z_m}\ns\ot\ldots\ot P_{z_1}\ts.
\end{equation}
For $a=1\lcd m$ let $\deg{\nns}_a$ be the linear operator
on this tensor product, corresponding to evaluation of the total
degree in $x_1\lcd x_n$ in the tensor factor $P_{z_a}\ts$;
that is the $a\ts$-th tensor factor when counting from right to left.
By restricting this tensor
product of $\Y(\gl_n)\ts$-modules to the subalgebra
$\Y(\g_n)\subset\Y(\gl_n)$ and then using the homomorphism
$\X(\g_n)\to\Y(\g_n)$ defined by \eqref{xy}, we can regard
the tensor product \eqref{pzm}
as a module over the extended twisted Yangian $\X(\g_n)\ts$.

\begin{corollary*}
\label{verma}
The bimodule\/ $\F_m(M_\mu)_{\ts\n}$ over\/ $\h$ and\/ $\X(\g_n)$ is
equivalent to the tensor product
\begin{equation}
\label{promuz}
P_{\mu_m+z}\ot P_{\mu_{m-1}+z+1}\ot\ldots\ot P_{\mu_1+z+m-1}
\vspace{6pt}
\end{equation}
pulled back through the automorphism of\/ $\X(\g_n)$ defined by \eqref{fus},
where $f(u)$ equals
\begin{equation}
\label{fuprod}
\prod_{a=1}^m\,\ts
\Bigl(1+\frac1{u-z-m+a-1-\mu_a}\ts\Bigr)\,;
\end{equation}
here $z=\mp\tts\frac12\ts$.
The element\/ $F_{\ts m-a+1,m-a+1}\in\h$ acts on 
\eqref{promuz} as the operator
\begin{equation}
\label{faa}
-\,n/2\ts+\ts\deg{\nns}_a-\mu_{\tts a}\,.
\end{equation}
\end{corollary*}

\begin{proof}
We have an embedding of $\gl_m$ to $\f_m$ such that
$E_{aa}\mapsto F_{aa}$ for $a=1\lcd m\ts$.
Then the Cartan subalgebra $\a$ of $\gl_m$ is identified
with the Cartan subalgebra $\h$ of $\f_m\ts$. Put
$\ab=m-a+1$ for short. If we regard the
weight $\mu$ as an element of $\a^\ast$, then
$$
\mu\ts(E_{\ts\ab\ts\ab\ts})=-\,\mu_{\tts a}
\quad\text{for}\quad
a=1\lcd m\ts.
$$
Let $U$ be the Verma module of the Lie algebra
$\gl_m$ corresponding to $\mu\in\a^\ast\ts$.
It is defined as the quotient of the algebra
$\U(\gl_m)$ by the left ideal generated by all elements
$E_{ab}$ with $a<b\ts$,
and by all elements $E_{aa}-\mu\ts(E_{aa})\ts$.
The Verma module $M_\mu$ of the Lie algebra $\f_m$ is then
equivalent to the module $V$ parabolically induced from
the $\gl_m$-module $U$.
Here we use the decomposition \eqref{pardec}.

Let $\s$ denote the subalgebra of the Lie algebra $\gl_m$
spanned by all elements $E_{ab}$ with $a>b\,$.
Using our embedding of $\gl_m$ to $\f_m\ts$, we can also
regard $\s$ as a subalgebra of $\f_m\ts$.
The Lie algebra $\n$ of $\f_m$ is then spanned by
$\r$ and $\s\ts$. By Proposition \ref{parres},
the bimodule $\F_m(M_\mu)_{\ts\n}$ over $\h$ and $\X(\g_n)$
is equivalent to $\Eb_m^{\ts z}\ts(U)_{\ts\s}$ where
$z=\mp\tts\frac12\ts$.
To describe the latter bimodule, let us firstly
consider the bimodule $\E_m^{\ts z}\ts(U)_{\ts\s}$
over $\a$ and $\Y(\gl_n)\ts$.
Using \cite[Corollary 2.4]{KN2}, the bimodule
$\E_m^{\ts z}\ts(U)_{\ts\s}$ is equivalent to the
tensor product of $\Y(\gl_n)\ts$-modules \eqref{promuz}
where the element $E_{\ts\ab\ts\ab\ts}\in\a$ acts as 
$\,\deg{\nns}_a-\mu_{\tts a}\ts$. After pulling the action
of Lie algebra $\gl_m$ on $\E_m^{\ts z}(U)$ back
through the automorphism \eqref{autom},
the element $E_{\ts\ab\ab\ts}\in\a$ will act on the tensor product
of vector spaces \eqref{promuz} as \eqref{faa}.

To complete the proof of Corollary \ref{verma},
recall that the action of $\X(\gl_n)$ on
$\Eb_m^{\ts z}\ts(U)$ differs from that on
$\E_m^{\ts z}\ts(U)$ by multiplying the series \eqref{rezhom} by 
\eqref{difz}. Using \eqref{hczu}, the series
$1+Z\ts(\ts u-z-m\ts)$ in $u^{-1}$ with the coefficients in $\Z(\gl_m)$
acts on the Verma module $U$ via scalar multiplication by
the series \eqref{fuprod}.
\qed
\end{proof}

By definition, the vector spaces of the two equivalent
bimodules in Corollary~\ref{verma} are
$(M_\mu\ot\P\ts(\CC^{\ts m}\ot\CC^{\ts n}))_{\ts\n}$
and $\P\ts(\CC^{\ts n})^{\ts\ot m}$ respectively.
We can determine a linear map from the latter vector space
to the former one,
by mapping $f_1\ot\ldots\ot f_m$ to
the class of the element $1_\mu\ot f$ in
the space of $\n\ts$-coinvariants. Here
for any $m$ polynomials $f_1\lcd f_m$ in the
$n$ anticommuting variables $x_1\lcd x_n$
the polynomial $f$ in the $m\ts n$
anticommuting variables $x_{11}\lcd x_{mn}$ is defined by setting
\begin{equation}
\label{fff}
f(x_{11}\lcd
x_{mn})\ts=
f_1(x_{11}\lcd x_{1n})\ts\ldots\ts f_m(x_{m1}\lcd x_{mn})\ts.
\end{equation}
This provides the bimodule equivalence in
Corollary \ref{verma},
see \cite[Corollary 2.4]{KN2}
and the remarks made immediately after stating Proposition \ref{parres} here.

For any $z\in\CC$ denote by $\Pp_z$
the $\Y(\gl_n)\ts$-module obtained by pulling $P_z$
back through the automorphism \eqref{transauto}
of $\Y(\gl_n)\ts$.
According to \eqref{pz}, the action of $\Y(\gl_n)$ on $\Pp_z$
is determined by the homomorphism $\Y(\gl_n)\to\PD\ts(\CC^{\ts n})\ts$,
\begin{equation}
\label{ppz}
T_{ij}(u)\mapsto\de_{ij}-
\frac{\,\th_i\ts\th_j\ts x_{\ts\bj}\,\d_{\ts\bi}}{u-z}\ .
\end{equation}

\begin{lemma*}
\label{ppl}
The\/ $\Y(\gl_n)\ts$-module $\Pp_z$
can also be obtained by pushing the action of\/
$\Y(\gl_n)$ on $P_{\ts-z-1}$ forward through the automorphism
of\/ $\PD\ts(\CC^{\ts n})$ such that for each $i=1\lcd n$
\begin{equation}
\label{onefour}
x_i\mapsto\th_i\,\d_{\ts\bi}
\quad\text{and}\quad
\d_i\mapsto\th_i\,x_{\ts\bi}\,,
\end{equation}
and by pulling the resulting action back through the automorphism
\eqref{fut} of\/ $\Y(\gl_n)$ where
\begin{equation}
\label{fuz}
g(u)=1-\frac1{u-z}\ .
\end{equation}
\end{lemma*}

Thus the action of $\Y(\gl_n)$ on $\Pp_z$ can also be determined by
the composition 
$$
\Y(\gl_n)\to
\Y(\gl_n)\underset{\tau_{z+1}}\longrightarrow\Y(\gl_n)\to
\U(\gl_n)\to\PD\ts(\CC^{\ts n})\to\PD\ts(\CC^{\ts n})\ts.
$$
Here the first map is the automorphism \eqref{fut} of $\Y(\gl_n)$
where the series $g(u)$ is given by \eqref{fuz}, the last map
is the automorphism \eqref{onefour} of $\PD\ts(\CC^{\ts n})$,
while the other three maps are defined like in \eqref{diag3}.

\begin{proof}
Applying the automorphism \eqref{onefour} to the right hand
side of \eqref{pz} and replacing the parameter $z$ there by $-z-1$ we get
$$
\de_{ij}+\frac{\,\th_i\ts\th_j\ts\d_{\ts\bi}\,x_{\ts\bj}}{u-z-1}
\,=\,
\de_{ij}+\frac{\,\de_{ij}-\th_i\ts\th_j\ts x_{\ts\bj}\,\d_{\ts\bi}}{u-z-1}
\,=\,
\frac{u-z}{u-z-1}
\,\Bigl(\de_{ij}-
\frac{\,\th_i\ts\th_j\ts x_{\ts\bj}\,\d_{\ts\bi}}{u-z}\ts\Bigr)
$$
which, after multiplying it by \eqref{fuz},
becomes the right hand side of \eqref{ppz}.
\qed
\end{proof}


\section*{\normalsize 3. Parabolic induction}
\setcounter{section}{3}
\setcounter{equation}{0}
\setcounter{theorem*}{0}

The twisted Yangian $\Y(\g_n)$ is not just a subalgebra of
$\Y(\gl_n)\ts$, it is also a right coideal of the coalgebra
$\Y(\gl_n)$ relative to the comultiplication \eqref{1.33}. Indeed,
let us apply this comultiplication to the
$i\com j$ entry of the $n\times n$ matrix $\Tp(-u)\,T(u)\ts$.
We get the sum
$$
\sum_{k=1}^n\,\th_i\ts\th_k\,
\De\ts(\,T_{\,\bk\ts\bi\ts}(-u)\,T_{kj}(u))\,=
\vspace{-2pt}
$$
$$
\sum_{g,h,k=1}^n\,\th_i\ts\th_j\,
(\,T_{\ts\bk\ts\bg\ts}(-u)\ot T_{\ts\bg\ts\bi\ts}(-u))\,
(\,T_{kh}(u)\ot T_{hj}(u))\,=
$$
$$
\sum_{g,h,k=1}^n\,\th_g\ts\th_k\,
T_{\ts\bk\ts\bg\ts}(-u)\,T_{kh}(u)
\ts\ot\ts
\th_i\ts\th_g\,T_{\ts\bg\ts\bi\ts}(-u)\,T_{hj}(u)\ts.\
\vspace{2pt}
$$
In the last displayed line, by performing the summation
over $k=1\lcd n$ in the first tensor factor we get
the $g\com h$ entry of the matrix $\Tp(-u)\,T(u)\ts$. Therefore
$$
\De\ts(\ts\Y(\g_n))\subset\Y(\g_n)\ot\Y(\gl_n)\ts.
$$
For the extended twisted Yangian $\X(\g_n)\ts$,
one can define a homomorphism of associative algebras
$$
\X(\g_n)\to\X(\g_n)\ot\Y(\gl_n)
$$
by assigning
\begin{equation}
\label{comod}
S_{ij}(u)\,\mapsto
\sum_{g,h=1}^n\,S_{gh}(u)
\ts\ot\ts
\th_i\ts\th_g\,T_{\ts\bg\ts\bi\ts}(-u)\,T_{hj}(u)\ts.
\vspace{2pt}
\end{equation}
The homomorphism property can be verified directly, see \cite[Section 3]{KN3}.
Using the homomorphism \eqref{comod}, the tensor product
of any modules over the algebras $\X(\g_n)$ and $\Y(\gl_n)$
becomes another module over $\X(\g_n)\ts$.

Furthermore, the homomorphism \eqref{comod} is a \textit{coaction}
of the Hopf algebra $\Y(\gl_n)$ on the algebra $\X(\g_n)\ts$.
Formally, one can define a homomorphism of associative algebras
$$
\X(\g_n)\to\X(\g_n)\ot\Y(\gl_n)\ot\Y(\gl_n)
$$
in two different ways: either by using the assignment \eqref{comod} twice,
or by using \eqref{comod} and then \eqref{1.33}. Both ways however
lead to the same result, see again \cite[Section 3]{KN3}.

Now for any positive integer $l$ consider the general linear Lie algebra
$\gl_{\ts2m+2l}$ and its subalgebra $\f_{m+l}\ts$. This subalgebra
is spanned by the elements $F_{ab}$ where
\begin{equation}
\label{mliml}
a\com b=-\ts m-l\lcd-1\com1\lcd m+l\ts.
\end{equation}
Extend the notation \eqref{eab} and \eqref{fab}
to all these indices $a\com b\ts$.
Now identify $\f_m$ with the subalgebra
of $\f_{m+l}$ spanned by the elements $F_{ab}$ where
$a\com b=-\ts m\lcd-1\com1\lcd m\ts$.
Choose the embedding of the Lie algebra $\gl_{\ts l}$ to $\f_{m+l}$
determined by the mappings
\begin{equation}
\label{gll}
E_{ab}\mapsto F_{m+a,m+b}
\ \quad\text{for}\ \quad
a\com b=1\lcd l\ts.
\end{equation}
Let $\q\com\qp$ be the subalgebras of $\f_{m+l}$
spanned respectively by the elements $F_{\ts ab}\com F_{\ts ba}$ where
$$
a=m+1\lcd m+l
\ \quad\text{and}\ \quad
b=-\ts m-l\lcd-1\com1\lcd m\,;
$$
these two subalgebras of $\f_{m+l}$ are nilpotent.
Put $\p=\f_m\op\ts\gl_{\ts l}\op\qp\ts$.
Then $\p$ is a maximal parabolic subalgebra of the reductive
Lie algebra $\f_{m+l}\ts$, and
$
\f_{m+l}=\q\op\p\ts.
$
We do not exclude the case $m=0$ here.
In this case the nilpotent subalgebras
$\q$ and $\qp$ of $\f_{m+l}$ become the Abelian subalgebras
$\r$ and $\rp$ of the Lie algebra $\f_{\ts l}\ts$;
see the decomposition \eqref{pardec}
where the positive integer $m$ is to be replaced by $l\ts$.
Note that the meaning of the symbols $\p$ and $\q$ here is
different from that in Section 0.

Let $V$ and $U$ be any modules of the Lie algebras
$\f_m$ and $\gl_{\ts l\ts}$ respectively.
Denote by $V\ns\bt U$ the $\f_{m+l}\ts$-module
\textit{parabolically induced\/}
from the $\f_m\op\ts\gl_{\ts l}\ts$-module $V\ns\ot U$.
To define $V\ns\bt U$, one first extends the action of the Lie algebra
$\f_m\op\ts\gl_{\ts l}$ on $V\ns\ot U$ to the Lie algebra $\p\ts$, so that
every element of the subalgebra $\qp\subset\p$ acts on $V\ns\ot U$ as zero.
By definition, $V\ns\bt U$ is the $\f_{m+l}\ts$-module induced from the
$\p\ts$-module $V\ns\ot U\ts$.
Note that here we have
a canonical embedding $V\ot U\to V\ns\bt U$ of
$\p\ts$-modules\ts;
we will denote by $\vuo$ the image of an element $v\ot u\in V\ot U$
under this embedding.

Consider the bimodule $\F_{m+l}\ts(\ts V\ns\bt U\ts)$
over $\f_{m+l}$ and $\X(\g_n)\ts$.
Here the action of $\,\X(\gl_n)$ commutes with the
action of the Lie algebra $\f_{m+l}\ts$,
and hence with the action of the subalgebra $\q\subset\f_{m+l}\ts$.
Therefore the vector space
$\F_{m+l}\ts(\ts V\ns\bt U\ts)_{\ts\q}$
of coinvariants of the action of the subalgebra $\q$
is a quotient of the $\X(\g_n)$-module
$\F_{m+l}\ts(\ts V\ns\bt U\ts)\ts$.
Note that the subalgebra
$\f_m\op\ts\gl_{\ts l}\subset\f_{m+l}$ also acts on this quotient.

For any $z\in\CC$ consider the bimodule $\E_{\ts l}^{\ts z}\ts(U)$
over $\gl_{\ts l}$ and $\Y(\gl_n)$ defined as in the end of Section~1.
Also consider the bimodule $\F_m(V)$ over $\f_m$ and $\X(\g_n)\ts$.
Using the homomorphism \eqref{comod}, the tensor product
of vector spaces $\F_m(V)\ot\ts\E_{\ts l}^{\ts z}\ts(U)$
becomes a module over $\X(\g_n)\ts$. This module is determined
by the homomorphism $\X(\g_n)\to\B_m\ot\A_l$ such that
for any $i\com j=1\lcd n$ the series $S_{ij}(u)$ is mapped to
\begin{equation}
\label{baser}
\sum_{g,h=1}^n\,\be_m\bigl(\ts S_{gh}(u)\bigr)
\ts\ot\ts
\th_i\ts\th_g\,
\al_{\ts l\ts}\bigr(\ts T_{\ts\bg\ts\bi\ts}(-\,u+z)\,T_{hj}(u+z)\bigl)\ts.
\end{equation}
Let us now map the series $S_{ij}(u)$ to
the series \eqref{baser} multiplied by
\begin{equation}
\label{zuser}
\bigl(\ts1\ot1\tts\bigr)
\ot
\bigl(\ts(1+Z\ts(\ts u-z-l\ts))\ot1\ts\bigr)
\ts\in\ts\B_m\ot\A_l\,[[u^{-1}]]\ts,
\end{equation}
see \eqref{zu}.
This mapping determines another homomorphism
$\X(\g_n)\to\B_m\ot\A_l\ts$.
Using it, we turn the vector space
of the $\X(\g_n)\ts$-module $\F_m(V)\ot\ts\E_{\ts l}^{\ts z}\ts(U)$
to yet another $\X(\g_n)\ts$-module, which will be denoted by
$\F_m(V)\ts\,\widetilde{\ot}\,\ts\E_{\ts l}^{\ts z}\ts(U)\,$.
Define an action of the Lie algebra $\gl_{\ts l}$ on
the latter $\X(\g_n)\ts$-module by pulling the action of $\gl_{\ts l}$
on $\E_{\ts l}\ts(U)$ back through the automorphism
\begin{equation}
\label{autol}
E_{ab}\,\mapsto\,-\,\de_{ab}\,n/2\,+\,E_{ab}
\quad\text{for}\quad
a\com b=1\lcd l\ts.
\end{equation}
The Lie algebra $\f_m$ acts on the $\X(\g_n)\ts$-module
$\F_m(V)\ts\,\widetilde{\ot}\,\ts\E_{\ts l}^{\ts z}\ts(U)$
via the tensor factor $\F_m(V)\ts$. Thus
$\F_m(V)\ts\,\widetilde{\ot}\,\ts\E_{\ts l}^{\ts z}\ts(U)$
becomes a bimodule
over the direct sum of Lie algebras $\f_m\op\ts\gl_{\ts l}$
and over the extended twisted Yangian $\X(\g_n)\ts$.
In the case $m=0$ the next theorem becomes Proposition~\ref{parres},
where the positive integer $m$ has to be replaced by $l\ts$.
Here we assume that $\F_{\ts0}(V)=\CC$ so that
$\be_{\ts0}(S_{ij}(u))=\de_{ij}\ts$.

\begin{theorem*}
\label{parind}
The bimodule $\F_{m+l}\ts(\ts V\ns\bt U\ts)_{\ts\q}$
over\/ $\f_m\op\ts\gl_{\ts l}$ and\/ $\X(\g_n)$ is
equivalent to
$\F_m(V)\,\ts\widetilde{\ot}\,\ts\E_{\ts l}^{\ts z}\ts(U)$
where $z=m\mp\tts\frac12\ts$.
\end{theorem*}

\begin{proof}
In the remainder of this section we shall prove Theorem \ref{parind}.
As vector spaces,
\begin{gather*}
\F_{m+l}\ts(\ts V\ns\bt U\ts)_{\ts\q}=
(\ts V\ns\bt U\ot\P\ts(\CC^{\ts m+l}\ot\CC^{\ts n}))_{\ts\q}\,,
\\
\F_m(V)\,\ts\widetilde{\ot}\,\ts\E_{\ts l}^{\ts z}\ts(U)=
V\ot\P\ts(\CC^{\ts m}\ot\CC^{\ts n})\ot U\ot\P\ts(\CC^{\ts l}\ot\CC^{\ts n})
\ts.
\end{gather*}

We can determine a linear map from the latter vector space to the
former one by mapping any element
$v\ot f\ot u\ot g$
to the class of $\vuo\ot f\ot g$ in the space of $\q\ts$-coinvarians.
Here $v\in V$, $f\in\P\ts(\CC^{\ts m}\ot\CC^{\ts n})$ and
$u\in U$, $g\in\P\ts(\CC^{\ts l}\ot\CC^{\ts n})$ whereas
the tensor product $f\ot g$
is identified with an element of $\P\ts(\CC^{\ts m+l}\ot\CC^{\ts n})$
in a natural way, which corresponds to the decomposition
\begin{equation}
\label{mnl}
\CC^{\ts m+l}\ot\CC^{\ts n}=
\CC^{\ts m}\ot\CC^{\ts n}\op\ts\CC^{\ts l}\ot\CC^{\ts n}\ts.
\end{equation}
We will show
that this map establishes an equivalence of bimodules
in Theorem~\ref{parind}.

The 
vector space of the $\f_{m+l}\ts$-module $V\ns\bt U$
can be identified with the tensor product $\U(\q)\ot V\ns\ot U$
where the Lie subalgebra $\q\subset\f_{m+l}$ acts via
left multiplication on the first tensor factor.
Then $\vuo=1\ot v\ot u\ts$, so that
the tensor product $V\ns\ot U$ gets identified with the subspace
\begin{equation}
\label{onesub}
1\ot V\ns\ot U\subset\ts\U(\q)\ot V\ns\ot U\ts.
\end{equation}
On this subspace,
every element of the subalgebra $\qp\subset\f_{m+l}$ acts as zero,
while the two direct summands
of subalgebra $\f_m\op\ts\gl_{\ts l}\subset\f_{m+l}$
act non-trivially only on the tensor factors $V$ and $U$ respectively.
All this determines the action of Lie algebra $\f_{m+l}$ on
$\U(\q)\ot V\ns\ot U$. Now consider $\F_{m+l}\ts(\ts V\ns\bt U\ts)$
as a $\f_{m+l}\ts$-module, we will denote it by $M$ for short.
Then $M$ is the tensor product of two $\f_{m+l}\ts$-modules,
\begin{equation}
\label{mmod}
M=(\ts V\ns\bt U)\ot\P\ts(\CC^{\ts m+l}\ot\CC^{\ts n})=
\U(\q)\ot V\ns\ot U\ot\P\ts(\CC^{\ts m+l}\ot\CC^{\ts n})\ts.
\end{equation}

The vector spaces of the $\X(\g_n)\ts$-module $\F_m\ts(V)$ and of the
$\Y(\gl_n)\ts$-module $\E_{\ts l}^{\ts z}\ts(U)$ are respectively
$V\ns\ot\ts\P\,(\CC^{\ts m}\ot\CC^{\ts n})$ and
$U\ns\ot\ts\P\,(\CC^{\ts l}\ot\CC^{\ts n})\ts$.
The action of the Lie algebra $\f_m$ on the first
vector space is defined by \eqref{gan}. By pulling back through the
automorphism \eqref{autol}, the action
of the Lie algebra $\gl_{\ts l}$ on the second
vector space is defined by mapping
$$
E_{ab}\,\mapsto-\,\de_{ab}\,{n}/2\,+\,
E_{ab}\ot1\,+\,
\sum_{k=1}^n\,
1\ot x_{ak}\,\d_{\ts bk}
\ \quad\text{for}\ \quad
a\com b=1\lcd l\ts.
$$
Identify the tensor product of these two vector spaces with the vector space
\begin{equation}
\label{vuprod}
V\ns\ot U\ot
\P\,(\CC^{\ts m}\ot\CC^{\ts n})\ot\P\,(\CC^{\ts l}\ot\CC^{\ts n})=
V\ns\ot U\ot
\P\,(\CC^{\ts m+l}\ot\CC^{\ts n})
\end{equation}
where we use the direct sum decomposition \eqref{mnl}.
We get an action of
the direct sum of Lie algebras $\f_m\op\ts\gl_{\ts l}$ on
the vector space \eqref{vuprod}.

Let us now define a linear map
$$
\ph:\,V\ns\ot U\ot\P\,(\CC^{\ts m+l}\ot\CC^{\ts n})\ts\to M\ts/\,\q\cdot M
$$
by the assignment
$$
\ph:\,y\ot x\ot t\,\mapsto\,1\ot y\ot x\ot t\,+\,\q\cdot M
$$
for any vectors $y\in V$, $x\in U$ and
$t\in\P\,(\CC^{\ts m+l}\ot\CC^{\ts n})\ts$.
The operator $\ph$ intertwines the actions of the Lie
algebra $\f_m\op\ts\gl_{\ts l}\ts$; see the definition \eqref{gan}
where $m$ is to be replaced by $m+l\ts$.
Let us demonstrate that the operator $\ph$ is bijective.

Firstly consider the
action of the Lie subalgebra $\q\subset\f_{m+l}$ on the vector space
$$
\P\ts(\CC^{\ts m+l})\ts=\ts
\P\ts(\CC^{\ts m})\ot\P\ts(\CC^{\ts l})\ts;
$$
the action is defined by \eqref{gan} where $n=1$,
and the integer $m$ is replaced by $m+l\ts$.
This vector space admits a descending filtration by the subspaces
$$
\mathop{\op}\limits_{K=N}^{l}\,
\P\,(\CC^{\ts m})\ot
\P^{\ts K}(\CC^{\ts l})
\ \quad\textrm{where}\ \quad
N=0\com1\lcd l\,.
$$
Here $\P^{\ts K}(\CC^{\ts l})$ stands for
the homogeneous subspace of $\P\,(\CC^{\ts l})$ of degree $K\ts$.
The action of the Lie algebra $\q$ on $\P\ts(\CC^{\ts m+l})$
preserves each of the filtration subspaces, and becomes trivial
on the associated graded space.

Similarly, for any $n=1\com2\com\ts\ldots$
the vector space $\P\,(\CC^{\ts m+l}\ot\CC^{\ts n})$ admits an descending
filtration by $\q\ts$-submodules such that $\q$ acts trivially on
each of the corresponding graded subspaces.
The latter filtration induces a filtration of $M$
by $\q\ts$-submodules such that on
the corresponding graded quotient
$\operatorname{gr}M\ts$, the Lie algebra $\q$
acts via left multiplication
on the first tensor factor $\U(\q)$ in \eqref{mmod}.
The space $V\ns\ot U\ot\P\,(\CC^{\ts m+l}\ot\CC^{\ts n})$
is therefore isomorphic to the space
of coinvariants $(\ts\operatorname{gr}M)_{\ts\q}$
via the bijective linear map
$$
y\ot x\ot t\,\mapsto\,
1\ot y\ot x\ot t\,+\,\q\cdot(\ts\operatorname{gr}M\ts)\ts.
$$
Therefore the linear map $\ph$ is bijective as well.
It now remains to show that the map $\ph$ intertwines the actions of
the algebra $\X(\g_n)\ts$.

In  this section we will use the symbol $\,\equiv\,$ to indicate
equalities in the algebra $\U(\ts\f_{m+l})$ modulo the left ideal
generated by the elements of the subalgebra $\qp\subset\f_{m+l}\ts$.
Any two elements of $\U(\ts\f_{m+l})$ related by $\,\equiv\,$
act on the subspace \eqref{onesub} in the same way.
We will extend the relation $\,\equiv\,$ to
formal power series in $u^{-1}$ with
coefficients in $\U(\ts\f_{m+l})\ts$, and then
to matrices whose entries are these series. Put
\begin{equation}
\label{v}
\textstyle
v=u\pm\frac12-m-l
\quad\text{and}\quad
w=-\ts u\pm\frac12-m-l\ts.
\end{equation}

The definition of the $\X(\g_n)\ts$-module
$M$ involves the $(2m+2l)\times(2m+2l)$ matrix whose
$a\com b$ entry is $\de_{ab}\,v+F_{ab}\ts$.
The rows and columns of this matrix are labelled by the indices
\eqref{mliml}. We proved in \cite[Section~3]{KN3} that the
inverse to this matrix
is related by $\,\equiv\,$ to the block matrix
\begin{equation}
\label{compound}
\begin{bmatrix}
\,H\,&0\,&0\,
\\
\,I\,&J\,&0\,
\\
\,P\,&Q\,&R\,
\end{bmatrix}
\end{equation}
where the blocks $H\com P\com R$ are certain
matrices of size $l\times l$
while the blocks $I\com J\com Q$
are certain matrices of sizes
$2m\times l\ts$,
$2m\times 2m\ts$,
$l\times 2m$ respectively.
Let us label the rows and columns of the
blocks by the same indices as in the compound matrix \eqref{compound}.
For instance, the rows and columns of the $l\times l$ matrix
$R$ are labelled by $m+1\lcd m+l\ts$.

Keeping to the notation of Section 2, let $F$ be the $2m\times2m$
matrix whose $c\com d$ entry is $F_{cd}$ for
$c\com d=-\ts m\lcd-1\com1\lcd m\ts$. Let $F(u)$ be the inverse
to the matrix $u+F$. The entries of the matrix $F(u)$ are
formal power series in $u^{-1}$ with coefficients in the
algebra $\U(\ts\f_m)\ts$, see \eqref{fabu}.
But now the algebra $\U(\ts\f_m)$ is regarded as a as subalgebra of
$\U(\ts\f_{m+l})\ts$. Let us denote by $W(u)$ the trace of the matrix $F(u)$,
as we did in Section 2.

Denote by $E$ the $l\times l$ matrix whose $a\com b$ entry is $F_{ab}$
for $a\com b=m+1\lcd m+l\ts$.
Using our embedding \eqref{gll}
of the Lie algebra $\gl_{\ts l}$ to $\f_{m+l}\ts$,
this notation agrees with the notation of Section 1.
But now we use the indices $a\com b=m+1\lcd m+l$
to label the rows and columns of the matrix $E\ts$.
Let $E(v)$ the inverse to the matrix $v+E\ts$.
Let $E_{ab}(u)$ be the $a\com b$ entry of the inverse matrix.
Let $Z(v)$ be the trace of the inverse matrix.
The coefficients of the formal power series $Z(v)$ in $v^{-1}$ belong to the
centre of the algebra $\U(\gl_{\ts l})\ts$, which is
now regarded as a subalgebra of $\U(\ts\f_{m+l})\ts$.
Further, for any indices $a\com b=m+1\lcd m+l\ts$ denote
$\Et_{ab}(v)=(v+l+\Ep\ts)_{\ts ba}^{-1}\ts$.
Then by Lemma \ref{eep}
\begin{equation}
\label{eet}
(1+Z(v))\,\Et_{ab}(v)=E_{ab}(v)\ts.
\end{equation}

Let $a\com b=m+1\lcd m+l$ and $c,d=-\ts m\lcd-1\com1\lcd m\ts$.
By \cite[Section 3]{KN3}
\begin{align*}
-\,H_{-b,-a}
&\,=\,
(1+Z(v))\,\Bigl(
\bigl(\ts W(v+l\ts)\ts\mp\ts\frac1{2u}+1\ts\bigr)\,
\Et_{ab}(w)\,\pm\,\frac1{2u}\ts\,\Et_{ab}(v)\Bigr)\,,
\\[10pt]
-\,I_{-d,-a}
&\,=\,
\sum_{b>m\ge c\ge-m}\,
\ep_{ad}\,F_{\ts bc}\,(1+Z(v))\,\Bigl(\ts
\ep_{cd}\,\Et_{ab}(w)\,F_{-d,-c\ts}(v+l\ts)
\\
&\hspace{92pt}
\pm\,\ts\frac{\Et_{ab}(w)-\Et_{ab}(v)}{2u}\,
F_{cd\ts}(v+l\ts)
\Bigr)\ts,
\end{align*}
\begin{gather*}
J_{cd}
\,=\,
(1+Z(v))\,F_{cd}\ts(v+l\ts)\ts,
\\[16pt]
P_{\ts b,-a}
\,=\,
\sum_{e,f>m}
F_{\ts f,-e}\,E_{\ts bf}(v)\,\Et_{ae}(w)\,\pm\hspace{-6pt}
\sum\limits_{\substack{e,f>m\\m\ge c,d\ge-m}}\hspace{-6pt}
\ep_{ad}\,\ts
F_{\ts f,-d}\,F_{\ts ec}\,E_{\ts be}(v)\,
\Et_{af}(w)\,F_{cd\ts}(v+l\ts)\,,
\\[2pt]
-\,Q_{ad}\,=
\sum_{e>m\ge c\ge-m}
F_{\ts ec}\,E_{ae}(v)\,F_{cd\ts}(v+l\ts)\ts,
\qquad
R_{\ts ab}\ts=\ts E_{ab}(v)\ts.
\end{gather*}
By definition
of $\X(\g_n)\ts$-module $M\ts$,
the action of $\X(\g_n)$ on the elements of the~subspace
\begin{equation}
\label{1vu}
1\ot V\ns\ot U\ot\P\ts(\CC^{\ts m+l}\ot\CC^{\ts n})\subset M
\end{equation}
can be now described by assigning to every series $S_{ij}(u)$
the following sum of series with coefficients in the algebra
$\B_{m+l}=\U(\ts\f_{m+l})\ot\ts\PD\ts(\CC^{\ts m+l}\ot\CC^{\ts n})\ts$:
\begin{gather}
\nonumber
\de_{ij}\,+\,\sum_{a,b>m}\,
R_{ab}\ts\ot\ts\th_i\,\th_j\,\d_{a\bi}\,x_{\ts b\bj}\ \,+
\\[2pt]
\nonumber
\sum_{a,b>m}\,
H_{-b,-a}\ts\ot\ts x_{\ts bi}\,\d_{aj}\ \,+
\\[2pt]
\nonumber
\sum_{a>m\ge d>0}\,(\,
I_{-d,-a}\ts\ot\ts x_{di}\,\d_{aj}
\ts+\ts
I_{d,-a}\ts\ot\ts\th_i\,\d_{d\bi}\,\d_{aj}\,)\ \,+
\\[2pt]
\nonumber
\sum_{m\ge c,d>0}\,(\,
J_{-c,-d}\ts\ot\ts x_{ci}\,\d_{dj}
+
J_{-c,d}\ts\ot\ts\th_j\,x_{ci}\,x_{d\bj}
+
J_{c,-d}\ts\ot\ts\th_i\,\d_{c\bi}\,\d_{dj}
+
J_{cd}\ts\ot\ts\th_i\,\th_j\,\d_{c\bi}\,x_{d\bj}\,)
\\[2pt]
\nonumber
+\,
\sum_{a,e>m}\,
P_{e,-a}\ts\ot\ts\th_i\,\d_{e\bi}\,\d_{aj}\ \,+
\\[2pt]
\label{fivlin}
\sum_{a>m\ge d>0}\,(\ts
Q_{a,-d}\ts\ot\ts\th_i\,\d_{a\bi}\,\d_{dj}
\ts+\ts
Q_{ad}\ts\ot\ts\th_i\,\th_j\,\d_{a\bi}\,x_{d\bj}\,)\,.
\end{gather}
Here for $a=1\lcd m+l$ and $i=1\lcd n$
we use the standard generators
$x_{ai}$ of the Grassmann algebra $\P\ts(\CC^{\ts m+l}\ot\CC^{\ts n})\ts$.
Then $\d_{ai}$ is the left derivation on
$\P\ts(\CC^{\ts m+l}\ot\CC^{\ts n})$ relative to $x_{ai}\ts$.
The generators $x_{ai}$ with $a\leqslant m$ and $a>m$
correspond to the first and the second
direct summands in \eqref{mnl}.

Consider the action of $\X(\g_n)$
on the elements of the subspace \eqref{1vu} modulo $\q\cdot M\ts$, using
the definition \eqref{gan} where $m$ is to be replaced by $m+l\ts$.
>From now till the end of this section,
we will be assuming that
$a\com b\com e\com f=m+1\lcd m+l$ while $c\com d=1\lcd m\ts$.
The indices $g\com h$ and $k$ will run through $1\lcd n\ts$.

By our description of the block $R\ts$,
the sum displayed in the first of six lines \eqref{fivlin} acts
on the elements of the subspace \eqref{1vu} as the sum
\begin{gather}
\nonumber
\de_{ij}\,\ts+\ts\sum_{a,b}
E_{ab}(v)\ts\ot\ts\th_i\,\th_j\,\d_{a\bi}\,x_{\ts b\bj}\ =
\\
\label{term0}
\de_{ij}\,(1+Z(v))\ts-\ts\sum_{a,b}
E_{ab}(v)\ts\ot\ts\th_i\,\th_j\,x_{\ts b\bj}\,\d_{a\bi}\,.
\end{gather}

By our description of the block $H\ts$,
the sum displayed in the second of six lines \eqref{fivlin}
acts on the elements of \eqref{1vu}
as the sum over the indices
$a\com b$ of
\begin{equation}
\label{term1}
-\,\,(1+Z(v))\,\Bigl(
\bigl(\ts W(v+l\ts)\ts\mp\ts\frac1{2u}+1\ts\bigr)\,
\Et_{ab}(w)\,\pm\,\frac1{2u}\ts\,\Et_{ab}(v)\Bigr)
\ts\ot\ts x_{\ts bi}\,\d_{aj}\,.
\vspace{4pt}
\end{equation}

By our description of the block $I\ts$,
the sum in the third of six lines \eqref{fivlin} acts
on the elements of \eqref{1vu}
as the sum over the indices $a\com b\com c\com d$ of
\begin{align*}
&\mp
F_{\ts b,-c}\,(1+Z(v))\,
\Bigl(\ts
\Et_{ab}(w)\,F_{-d,c\ts}(v+l\ts)
+
\frac{\Et_{ab}(w)-\Et_{ab}(v)}{2u}\,
F_{-c,d\ts}(v+l\ts)
\Bigr)
\,\ot\,x_{di}\,\d_{aj}
\\[6pt]
&-
F_{\ts bc}\,(1+Z(v))\,
\Bigl(\ts
\Et_{ab}(w)\,F_{-d,-c\ts}(v+l\ts)
\pm
\frac{\Et_{ab}(w)-\Et_{ab}(v)}{2u}\,
F_{cd\ts}(v+l\ts)
\Bigr)
\,\ot\,x_{di}\,\d_{aj}
\\[6pt]
&\mp
F_{\ts b-c}\,(1+Z(v))\,\Bigl(\ts
\Et_{ab}(w)\,
F_{dc\ts}(v+l\ts)
\pm
\frac{\Et_{ab}(w)-\Et_{ab}(v)}{2u}\,
F_{-c,-d\ts}(v+l\ts)
\Bigr)
\,\ot\,\th_i\,\d_{d\bi}\,\d_{aj}
\\[6pt]
&-
F_{\ts bc}\,(1+Z(v))\,\Bigl(\ts
\Et_{ab}(w)\,
F_{d,-c\ts}(v+l\ts)
+
\frac{\Et_{ab}(w)-\Et_{ab}(v)}{2u}\,
F_{c,-d\ts}(v+l\ts)
\Bigr)
\,\ot\,\th_i\,\d_{d\bi}\,\d_{aj}\,.
\end{align*}
Here $F_{\ts b,-c}\in\q$ and $F_{\ts bc}\in\q\ts$. Hence
modulo $\q\cdot M\ts$, the expression displayed in the latter four
lines acts on the elements of \eqref{1vu} as the sum over the index $k$ of
\begin{align*}
\Bigl(
&\pm
\Bigl(\ts
\Et_{ab}(w)\,F_{-d,c\ts}(v-l\ts)
+
\frac{\Et_{ab}(w)-\Et_{ab}(v)}{2u}\,
F_{-c,d\ts}(v+l\ts)
\Bigr)
\,\ot\,
\th_k\,x_{\ts b\bk}\,x_{ck}\,
x_{di}\,\d_{aj}
\\[4pt]
&+
\Bigl(\ts
\Et_{ab}(w)\,F_{-d,-c\ts}(v+l\ts)
\pm
\frac{\Et_{ab}(w)-\Et_{ab}(v)}{2u}\,
F_{cd\ts}(v+l\ts)
\Bigr)
\,\ot\,x_{\ts bk}\,\d_{ck}\,x_{di}\,\d_{aj}
\\[4pt]
&\pm
\Bigl(\ts
\Et_{ab}(w)\,
F_{dc\ts}(v+l\ts)
\pm
\ts\frac{\Et_{ab}(w)-\Et_{ab}(v)}{2u}\,
F_{-c,-d\ts}(v+l\ts)
\Bigr)
\,\ot\,\th_i\,
\th_k\,x_{\ts b\bk}\,x_{ck}\,
\d_{d\bi}\,\d_{aj}
\\[4pt]
&+
\Bigl(\ts
\Et_{ab}(w)\,
F_{d,-c\ts}(v+l\ts)
+
\frac{\Et_{ab}(w)-\Et_{ab}(v)}{2u}\,
F_{c,-d\ts}(v+l\ts)
\Bigr)
\,\ot\,
\th_i\,
x_{\ts bk}\,\d_{ck}\,
\d_{d\bi}\,\d_{aj}\,\Bigr)
\end{align*}
\begin{equation}
\label{term2}
\times\ \bigl(\ts(1+Z(v))\ot1\ts\bigr)\,.
\vspace{4pt}
\end{equation}

By our description of the block $J\ts$,
the sum displayed in the fourth of six lines \eqref{fivlin}
acts on the elements of \eqref{1vu} as the sum over $c\com d$ of
\begin{gather}
\nonumber
\bigl(\ts(1+Z(v))\ot1\ts\bigr)\,\bigl(\,
F_{-c,-d\ts}(v+l\ts)\ts\ot\ts x_{ci}\,\d_{dj}
\ts+\ts
F_{-c,d\ts}(v+l\ts)\ts\ot\ts\th_j\,x_{ci}\,x_{d\bj}
\\[6pt]
\label{term3}
+\,\ts
F_{c,-d\ts}(v+l\ts)\ts\ot\ts\th_i\,\d_{c\bi}\,\d_{dj}
\ts+\ts
F_{cd\ts}(v+l\ts)\ts\ot\ts\th_i\,\th_j\,\d_{c\bi}\,x_{d\bj}\,\bigr)\,.
\end{gather}

By our description of the block $P\ts$,
the sum displayed in the fifth of six lines \eqref{fivlin}
acts on the elements of the subspace \eqref{1vu} as
the sum over the indices $a\com b\com e\com f$~of
$$
F_{\ts f,-b}\,E_{ef}(v)\,\Et_{ab}(w)
\ts\ot\ts\th_i\,\d_{e\bi}\,\d_{aj}
$$
plus the action of the sum over the indices
$a\com b\com c\com d\com e\com f$ of
\begin{align*}
&
F_{fd}\,
F_{\ts b,-c}\,
E_{eb}(v)\,\Et_{af}(w)\,F_{-c,-d\ts}(v+l)
\ts\ot\ts\th_i\,\d_{e\bi}\,\d_{aj}
\\[6pt]
\pm\,&
F_{f,-d}\,
F_{\ts b,-c}\,
E_{eb}(v)\,\Et_{af}(w)\,F_{-c,d}(v+l)
\ts\ot\ts\th_i\,\d_{e\bi}\,\d_{aj}
\\[6pt]
+\,&
F_{fd}\,
F_{\ts bc}\,
E_{eb}(v)\,\Et_{af}(w)\,F_{c,-d\ts}(v+l)
\ts\ot\ts\th_i\,\d_{e\bi}\,\d_{aj}
\\[6pt]
\pm\,&
F_{f,-d}\,
F_{\ts bc}\,
E_{eb}(v)\,\Et_{af}(w)\,F_{cd\ts}(v+l)
\ts\ot\ts\th_i\,\d_{e\bi}\,\d_{aj}\,.
\end{align*}
Modulo $\q\cdot M\ts$, here the expression to be summed over
the indices $a\com b\com e\com f$ acts
on the elements of the subspace \eqref{1vu} as the sum over the index $k$ of
$$
-\,E_{ef}(v)\,\Et_{ab}(w)
\ts\ot\ts
\th_i\,\th_k\,
x_{f\bk}\,x_{\ts bk}\,
\d_{e\bi}\,\d_{aj}
$$
while the expression to be summed over $a\com b\com c\com d\com e\com f$ acts
as the sum over $g\com h$ of
\begin{align*}
\bigl(\ts E_{eb}(v)\,\Et_{af}(w)\ot1\ts\bigr)\,
\bigl(\,
&F_{-c,-d\ts}(v+l)
\ts\ot\ts
\th_i\,
\th_g\,
x_{\ts b\bg}\,x_{cg}\,
x_{fh}\,\d_{dh}\,
\d_{e\bi}\,\d_{aj}
\\[6pt]
\pm\,&
F_{-c,d}(v+l)
\ts\ot\ts
\th_i\,\th_g\,\th_h\,
x_{\ts b\bg}\,x_{cg}\,
x_{f\bh}\,x_{dh}\,
\d_{e\bi}\,\d_{aj}
\\[6pt]
+\,&
F_{c,-d\ts}(v+l)
\ts\ot\ts
\th_i\,
x_{\ts bg}\,\d_{cg}\,
x_{\ts fh}\,\d_{dh}\,
\d_{e\bi}\,\d_{aj}
\\[6pt]
\pm\,&
F_{cd\ts}(v+l)
\ts\ot\ts
\th_i\,\th_h\,
x_{\ts bg}\,\d_{cg}\,
x_{f\bh}\,x_{dh}\,
\d_{e\bi}\,\d_{aj}\ts\bigr)\,.
\end{align*}

We have $\th_{\bk}=\pm\,\th_k$ for $k=1\lcd n\ts$.
Using the commutation relations in the ring
$\PD\ts(\CC^{\ts m+l}\ot\CC^{\ts n})\ts$,
the sum over the index $k$ above equals the sum over $k$ of
\begin{equation}
\label{term41}
E_{ef}(v)\,\Et_{ab}(w)
\ts\ot\ts
\th_i\,\th_k\,
x_{f\bk}\,\d_{e\bi}\,
x_{\ts bk}\,\d_{aj}
\end{equation}
plus
\begin{equation}
\label{term410}
\mp\,\ts\de_{\ts be}\,
E_{ef}(v)\,\Et_{ab}(w)
\ts\ot\ts
x_{f\bi}\,\d_{aj}\ts.
\vspace{6pt}
\end{equation}
Similarly, the sum over the indices
$g\com h$ equals the sum over $g\com h$ of
\begin{align*}
\nonumber
\bigl(\,
F_{-c,-d\ts}(v+l\ts)\ts\ot\ts x_{cg}\,\d_{dh}
&\ts+\ts
F_{-c,d\ts}(v+l\ts)\ts\ot\ts\th_h\,x_{cg}\,x_{d\bh}\,\,+
\\[6pt]
F_{c,-d\ts}(v+l\ts)\ts\ot\ts\th_g\,\d_{c\bg}\,\d_{dh}
&\ts+\ts
F_{cd\ts}(v+l\ts)\ts\ot\ts\th_g\,\th_h\,\d_{c\bg}\,x_{d\bh}\,\bigr)\ \times
\end{align*}
\begin{equation}
\label{term421}
\bigl(\ts E_{eb}(v)\ot\th_i\,\th_g\,x_{\ts b\bg}\,\d_{e\bi}\ts\bigr)\,
\bigl(\ts\Et_{af}(w)\ot x_{\ts fh}\,\d_{aj}\ts\bigr)
\end{equation}
plus the sum over $k$ of
\begin{align}
\nonumber
\bigl(\,\de_{ef}\,E_{eb}(v)\,&\Et_{af}(w)\ot1\ts\bigr)\ \times
\\[6pt]
\nonumber
\bigl(\ts
-\,
F_{-c,-d\ts}(v+l\ts)
\ts\ot\ts\th_i\,
\th_k\,x_{\ts b\bk}\,x_{ck}\,
\d_{d\bi}\,\d_{aj}
&\,\mp\,
F_{-c,d\ts}(v+l\ts)
\ts\ot\ts
\th_k\,x_{\ts b\bk}\,x_{ck}\,
x_{di}\,\d_{aj}
\\[6pt]
\label{term422}
-\,
F_{c,-d\ts}(v+l\ts)
\ts\ot\ts
\th_i\,
x_{\ts bk}\,\d_{ck}\,
\d_{d\bi}\,\d_{aj}
&\,\mp\,
F_{cd\ts}(v+l\ts)
\ts\ot\ts
x_{\ts bk}\,\d_{ck}\,x_{di}\,\d_{aj}\ts\bigl)\,.
\end{align}

By our description of the block $Q\ts$,
the sum displayed in the last of the six lines \eqref{fivlin}
acts on the elements of \eqref{1vu}
as the sum over $a\com b\com c\com d$ of
$$
\begin{aligned}
&-\,(\ts
F_{\ts b,-c}\,E_{ab}(v)\,F_{-c,-d\ts}(v+l\ts)
\ts+\ts
F_{\ts bc}\,E_{ab}(v)\,F_{c,-d\ts}(v+l\ts))
\ts\ot\ts\th_i\,\d_{a\bi}\,\d_{dj}
\\[4pt]
&-\,(\ts
F_{\ts b,-c}\,E_{ab}(v)\,F_{-c,d\ts}(v+l\ts)
\ts+\ts
F_{\ts bc}\,E_{ab}(v)\,F_{cd\ts}(v+l\ts))
\ts\ot\ts\th_i\,\th_j\,\d_{a\bi}\,x_{d\bj}\,.
\end{aligned}
$$
Modulo $\q\cdot M\ts$, the expression in the above two lines acts
on the elements of the subspace \eqref{1vu} as the sum over $k$ of
\begin{gather*}
\bigl(\ts E_{ab}(v)\ot1\ts\bigr)\ \times
\\[4pt]
\bigl(\,
F_{-c,-d\ts}(v+l\ts)
\ts\ot\ts
\th_i\,\th_k\,x_{\ts b\bk}\,x_{ck}\,
\d_{a\bi}\,\d_{dj}
\,+\,
F_{c,-d\ts}(v+l\ts)
\ts\ot\ts
\th_i\,
x_{\ts bk}\,\d_{ck}\,
\d_{a\bi}\,\d_{dj}\,\,+
\\[4pt]
\
F_{-c,d\ts}(v+l\ts)
\ts\ot\ts
\th_k\,x_{\ts b\bk}\,x_{ck}\,
\th_i\,\th_j\,\d_{a\bi}\,x_{d\bj}
\,+\,
F_{cd\ts}(v+l\ts)
\ts\ot\ts
x_{\ts bk}\,\d_{ck}\,\th_i\,\th_j\,\d_{a\bi}\,x_{d\bj}\ts\bigr)\ts.
\end{gather*}
Note that this sum over the index $k$
can be rewritten as the sum over $k$ of
\begin{gather}
\nonumber
\bigl(\ts
F_{-c,-d\ts}(v+l\ts)
\ts\ot\ts
x_{ck}\,\d_{dj}
\,+\,
F_{-c,d\ts}(v+l\ts)
\ts\ot\ts
\th_j\,
x_{ck}\,\,x_{d\bj}\,\,+
\\[4pt]
\nonumber
\ \ts
F_{c,-d\ts}(v+l\ts)
\ts\ot\ts
\th_k\,
\d_{c\bk}\,\,\d_{dj}
\,+\,
F_{cd\ts}(v+l\ts)
\ts\ot\ts
\th_k\,\th_j
\d_{c\bk}\,x_{d\bj}\ts\bigr)\ \times
\\[4pt]
\label{term5}
\bigl(\ts-\,E_{ab}(v)
\ts\ot\ts
\th_i\,\th_k\,
x_{\ts b\bk}\,\d_{a\bi}\ts\bigr)\,.
\end{gather}

Consider the sum of the expressions \eqref{term422}
over the running indices $e\com f\ts$. Add this sum to
the expression displayed in the five lines \eqref{term2}. Using the relation
\begin{equation}
\label{can3}
\sum_{e}\,\Et_{\ts eb}(v)\,\Et_{ae}(w)\,=\,
\frac{\ts\Et_{ab}(w)-\Et_{ab}(v)\ts}{2u}
\end{equation}
together with \eqref{eet} and performing cancellations,
we get the expression
\begin{gather*}
\bigl(\,
\pm\,\ts
F_{-d,c\ts}(v+l\ts)
\ts\ot\ts
\th_k\,x_{\ts b\bk}\,x_{ck}\,
x_{di}\,\d_{aj}
\ts+\ts
F_{-d,-c\ts}(v+l\ts)
\ts\ot\ts x_{\ts bk}\,\d_{ck}\,x_{di}\,\d_{aj}\,\,+
\\[4pt]
F_{d,-c\ts}(v+l\ts)
\ts\ot\ts
\th_i\,
x_{\ts bk}\,\d_{ck}\,
\d_{d\bi}\,\d_{aj}
\,\pm\,
F_{dc\ts}(v+l\ts)
\ts\ot\ts\th_i\,
\th_k\,x_{\ts b\bk}\,x_{ck}\,
\d_{d\bi}\,\d_{aj}\ts\bigr)\ \times
\\[4pt]
\bigl(\ts(1+Z(v))\,\Et_{ab}(w)\ot1\ts\bigr)\,.
\end{gather*}
After exchanging the running indices $c$ and $d\ts$,
the sum over the index $k$ of the expressions in the last three
displayed lines
can be rewritten as
\begin{equation}
\label{term24221}
\de_{cd}\,(1+Z(v))\,
\bigl(\ts F_{-c,-d\ts}(v+l\ts)+F_{cd\ts}(v+l\ts)\bigr)\,
\Et_{ab}(w)
\ts\ot\ts
x_{\ts bi}\,\d_{aj}
\end{equation}
plus the sum over $k$ of
\begin{gather}
\nonumber
\bigl(\ts
F_{-c,-d\ts}(v+l\ts)
\ts\ot\ts
x_{ci}\,\d_{dk}
\,+\,
F_{-c,d\ts}(v+l\ts)
\ts\ot\ts
\th_k\,
x_{ci}\,x_{d\bk}\,\,+
\\[4pt]
\nonumber
\ \ \ts
F_{c,-d\ts}(v+l\ts)
\ts\ot\ts
\th_i\,
\d_{c\bi}\,\d_{dk}
\,+\,
F_{cd\ts}(v+l\ts)
\ts\ot\ts
\th_i\,\th_k
\d_{c\bi}\,x_{d\bk}
\ts\bigr)\ \times
\\[4pt]
\label{term24222}
\bigl(\ts-\,(1+Z(v))\,\Et_{ab}(w)
\ts\ot\ts
x_{\ts bk}\,\d_{aj}\ts\bigr)\,.
\end{gather}
Here we again used the commutation relations in the ring
$\PD\ts(\CC^{\ts m+l}\ot\CC^{\ts n})\ts$.

Let us now perform the summation over all running indices
in the four expressions
\eqref{term3},\eqref{term421},\eqref{term5},\eqref{term24222}
and then take their total. By exchanging the running indices
$b$ and $f$ in \eqref{term421}, and by replacing the running index $k$
in \eqref{term5},\eqref{term24222} by $g\com h$ respectively,
the total can be written as the sum over indices $c\com d$ and
$g\com h$ of
\begin{gather}
\nonumber
\bigl(\ts(1+Z(v))\ot1\ts\bigr)\ \times
\\[4pt]
\nonumber
\bigl(\,
F_{-c,-d\ts}(v+l\ts)\ts\ot\ts x_{cg}\,\d_{dh}
\,+\,
F_{-c,d\ts}(v+l\ts)\ts\ot\ts\th_h\,x_{cg}\,x_{d\bh}\,\,+
\\[4pt]
\nonumber
\ \ \,\ts
F_{c,-d\ts}(v+l\ts)\ts\ot\ts\th_g\,\d_{c\bg}\,\d_{dh}
\,+\,
F_{cd\ts}(v+l\ts)\ts\ot\ts\th_g\,\th_h\,\d_{c\bg}\,x_{d\bh}\,\bigr)\ \times
\\[4pt]
\label{cdgh}
\bigl(\ts\de_{ig}-
\sum_{e,f}\,
\Et_{ef}(v)\ot\th_i\,\th_g\,x_{\ts f\bg}\,\d_{e\bi}\,\bigr)\,
\bigl(\ts\de_{hj}-
\sum_{a,b}\,
\Et_{ab}(w)\ts\ot\ts x_{\ts bh}\,\d_{aj}\ts\bigr)\,.
\end{gather}

Let us perform the summation in \eqref{term410}
over the running indices $b\com e\ts$. Then let us replace
the running index $f$ by the index $b\ts$, which becomes free
after the summation. By adding the resulting sum to
the expession \eqref{term1} we get
$$
-\,(1+Z(v))\ts
\bigl(\ts W(v+l\ts)+1\ts\bigr)\,\Et_{ab}(w)
\ts\ot\ts x_{\ts bi}\,\d_{aj}
$$
due to \eqref{eet} and \eqref{can3}.
By performing the summation in \eqref{term24221}
over the running indices $c\com d$ and then adding the result
to the last displayed expression, we get
\begin{equation}
\label{term000}
-\,(1+Z(v))\,\Et_{ab}(w)
\ts\ot\ts x_{\ts bi}\,\d_{aj}\,.
\end{equation}

Now do the summation over all running indices
in the two expressions \eqref{term41},\eqref{term000} and then
add the two resulting sums to \eqref{term0}. By using the relation
\eqref{eet} once again, the total can be written as
the sum over the index $k$ of
\begin{gather}
\nonumber
\bigl(\ts(1+Z(v))\ot1\ts\bigr)\ \times
\\[4pt]
\label{k}
\bigl(\ts\de_{ik}-
\sum_{e,f}\,
\Et_{ef}(v)\ot\th_i\,\th_k\,x_{\ts f\bk}\,\d_{e\bi}\,\bigr)\,
\bigl(\ts\de_{kj}-
\sum_{a,b}\,
\Et_{ab}(w)\ts\ot\ts x_{\ts bk}\,\d_{aj}\ts\bigr)\,.
\end{gather}

By using the definition of the series $\Et_{ab}(v)$
as given before the relation \eqref{eet},
\begin{gather*}
\textstyle
\Et_{ef}(v)=(\ts v+l+\Ep)^{-1}_{fe}=-\,(-\ts u\mp\frac12+m-\Ep)^{-1}_{fe}\,,
\\[4pt]
\textstyle
\Et_{ab}(w)=(\ts w+l+\Ep)^{-1}_{ba}=-\,(\ts u\mp\frac12+m-\Ep)^{-1}_{ba}\,.
\end{gather*}
We also used the definitions \eqref{v}.
Hence the sum of the expressions \eqref{cdgh} over
the indices $c\com d$ and $g\com h$ plus the sum of the expessions
\eqref{k} over the index $k$ can be rewritten
as the sum over the indices $g\com h$ of the series in $u^{-1}$,
\begin{gather*}
\textstyle
\bigl(\ts(1+Z(\ts u\pm\frac12-m-l\ts))\ot1\ts\bigr)\ \times
\\[8pt]
\bigr(\,\de_{gh}\ \ts+
\sum_{c,d}\
\bigl(\,
\textstyle
F_{-c,-d\ts}(\ts u\pm\frac12-m\ts)
\ts\ot\ts
x_{cg}\,\d_{dh}
\,+\,
\textstyle
F_{-c,d\ts}(\ts u\pm\frac12-m\ts)
\ts\ot\ts
\th_h\,x_{cg}\,x_{d\bh}
\\
+\,\ts
\textstyle
F_{c,-d\ts}(\ts u\pm\frac12-m\ts)
\ts\ot\ts
\th_g\,\d_{c\bg}\,\d_{dh}
\,+\,
\textstyle
F_{cd\ts}(\ts u\pm\frac12-m\ts)
\ts\ot\ts
\th_g\,\th_h\,\d_{c\bg}\,x_{d\bh}\,\bigr)\bigr)\ \times
\\[10pt]
\bigl(\,\de_{ig}+
\sum_{e,f}\,
\textstyle
(-\ts u\mp\frac12+m-\Ep)^{-1}_{fe}
\displaystyle
\ot\th_i\,\th_g\,x_{\ts f\bg}\,\d_{e\bi}\,\bigr)\ \times
\\
\bigl(\,\de_{hj}+
\sum_{a,b}\,
\textstyle
(\ts u\mp\frac12+m-\Ep)^{-1}_{ba}
\displaystyle
\ts\ot\ts x_{\ts bh}\,\d_{aj}\ts\bigr)
\end{gather*}
with coefficients in the algebra
$\U(\ts\f_{m}\op\gl_{\ts l})\ot\ts\PD\ts(\CC^{\ts m+l}\ot\CC^{\ts n})\ts$.
By mapping the series $S_{ij}(u)$ to this sum
we describe the action of the extended twisted Yangian
$\X(\g_n)$ on the subspace \eqref{1vu} modulo $\q\cdot M$.
By comparing this sum with the product of the series
\eqref{baser} and \eqref{zuser} where $z=m\mp\frac12\,$,
we now prove that the map $\ph$  intertwines the actions of $\X(\g_n)\ts$;
here we use \eqref{ehom} and \eqref{fhom}.
This completes the proof of Theorem \ref{parind}.
\qed
\end{proof}


\section*{\normalsize 4. Zhelobenko operators}
\setcounter{section}{4}
\setcounter{equation}{0}
\setcounter{theorem*}{0}

Let us consider the \textit{hyperoctahedral group\/} $\H_m\ts$. This
is the semidirect product of the symmetric group $\Sym_m$ and
the Abelian group $\ZZ_2^m\ts$, where $\Sym_m$ acts by permutations
of the $m$ copies of $\ZZ_2\ts$.
In this section, we assume that $m>0\ts$.
The group $\H_m$ is generated by the elements $\si_a$ with $a=1\lcd m\ts$.
The elements $\si_a$ with the indices $a=1\lcd m-1$ are elementary
transpositions generating the symmetric group $\Sym_m\ts$, so that
$\si_a=(a\com a+1)\ts$. Then $\si_m$ is the generator of the
$m\ts$-th factor $\ZZ_2$ of $\ZZ_2^m\ts$.
The elements
$\si_1\lcd\si_m\in\H_m$ are involutions and satisfy the braid relations
\begin{align*}
\si_a\,\si_{a+1}\,\si_a
&\,=\,
\si_{a+1}\,\si_{a}\,\si_{a+1}
\!\!\!\quad\quad\textrm{for}\ \quad
a=1\lcd m-2\ts;
\\
\si_a\,\si_{\ts b}
&\,=\,
\si_{\ts b}\,\si_a
\hspace{34pt}
\ \quad\textrm{for}\ \quad
a=1\lcd b-2\ts;
\\
\si_{m-1}\,\si_m\,\si_{m-1}\,\si_m
&\,=\,
\si_{m}\,\si_{m-1}\,\si_{m}\,\si_{m-1}.
\end{align*}

Note that $\H_m$ is the Weyl group of the simple Lie algebra $\sp_{2m}\ts$.
Let $\Hh_m$ be the braid group corresponding to $\sp_{2m}\ts$.
It is generated by the elements $\sih_1\lcd\sih_m$
which by definition satisfy the above displayed relations, instead of
the involutions $\si_1\lcd\si_m$ respectively. For any reduced decomposition
$\si=\si_{a_1}\ldots\si_{a_K}$ in $\H_m$ put
\begin{equation}
\label{sih}
\sih=\sih_{a_1}\ldots\,\sih_{a_K}.
\end{equation}
The definition of $\,\sih$
is independent of the choice of a reduced decomposition~of~$\si\ts$.

The 
group $\H_m$ also contains
the Weyl group of the reductive
Lie algebra $\so_{2m}$ as a subgroup of index two.
Denote this subgroup by $\Hp_m\ts$,
it is generated by the elementary transpositions
$\si_1\lcd\si_{m-1}$ and by the involution
$\sip_m=\si_m\,\si_{m-1}\,\si_m\ts$. Along with the braid relations
between $\si_1\lcd\si_{m-1}\ts$, we also have the braid relations
involving $\sip_m\ts$,
\begin{align*}
\si_a\,\sip_{\ts m}
&\,=\,
\sip_{\ts m}\,\si_a
\ \quad\textrm{for}\ \quad
a=1\lcd m-3\com m-1\ts;
\\
\si_{m-2}\,\sip_m\,\si_{m-2}
&\,=\,
\sip_m\,\si_{m-2}\,\sip_m\ts.
\end{align*}
When $m>1\ts$, the braid group of $\so_{2m}$
is generated by $m$ elements
satisfying the same braid relations instead of the
$m$ involutions $\si_1\lcd\si_{m-1}\com\sip_m$ respectively.
When $m=1\ts$, the braid group corresponding to $\f_m=\so_2$ consists
of the identity element only.

Now let the indices
$c\com d$ run through $-m\lcd-1\com1\lcd m\ts$.
For $c>0$ we denote $\cb=m+1-c\ts$; for $c<0$
denote $\cb=-\,m-1-c\ts$.
Consider a representation $\si\mapsto\sib$ of the
group $\H_m$ by permutations of $-m\lcd-1\com1\lcd m$ such that
\begin{equation}
\label{sib}
\sib\ts(c)=\overline{\si(\ts\cb\ts)\!}\,
\quad\text{for}\quad\si\in\Sym_m
\end{equation}
and
$\sib_m\ts(c)=-c$ if $|c|=1\ts$,
while
$\sib_m\ts(c)=c$ if $|c|>1\ts$.
We can define an action of the braid group $\Hh_m$ by automorphisms of
the Lie algebra $\f_m\ts$, by the assignments
\begin{align}
\label{siact}
\sih:\,F_{cd}
&\,\mapsto\,
F_{\sib(c)\ts\sib(d)}
\quad\text{for}\quad\si\in\Sym_m\ts,
\\[2pt]
\label{simact}
\sih_m:\,F_{cd}
&\,\mapsto\,
(\ts\pm\ts1\ts)^{\ts\de_{c1}\ts+\,\de_{d1}}F_{\sib_m(c)\ts\sib_m(d)}\,;
\end{align}
cf.\ \cite{T}.
According to our convention on double signs,
the upper sign in $\pm$ corresponds to $\f_m=\so_{2m}\ts$,
while the lower sign corresponds to $\f_m=\sp_{2m}\ts$.
The automorphism property can be checked
by using the relations \eqref{ufmrel}, see
the proof of Part~(i) of Lemma~\ref{lemma41} below.
This action of the group $\Hh_m$ on $\f_m$
extends to an action of $\Hh_m$ by automorphisms of
the associative algebra $\U(\ts\f_m)\ts$.
Note that in the case $\f_m=\so_{2m}$
the action of $\Hh_m$ on $\U(\f_m)$
factors to an action of the group $\H_m\ts$.

Further, one can define an action
of the braid group $\Hh_m$ by automorphisms of the algebra
$\PD\ts(\CC^{\ts m}\ot\CC^{\ts n})$ in the following way. Put
\begin{align}
\nonumber
\sih\ts(x_{ai})&=x_{\ts\sib(a)\ts i}
\hspace{-16pt}
&
\textrm{and}
&
&
\hspace{-16pt}
\sih\ts(\d_{ai})&=\d_{\ts\sib(a)\ts i}
\hspace{-16pt}
&
\textrm{for}
&
&
&
\hspace{-16pt}
\si\in\Sym_m\,,
\\
\nonumber
\sih_m(x_{ai})&=x_{ai}
\hspace{-16pt}
&
\textrm{and}
&
&
\hspace{-16pt}
\sih_m(\d_{ai})&=\d_{ai}
\hspace{-16pt}
&
\textrm{for}
&
&
&
\hspace{-16pt}
a>1\ts,
\\
\label{fourier}
\sih_m\ts(x_{1i})&=\th_i\,\d_{\ts1\bi}
\hspace{-16pt}
&
\textrm{and}
&
&
\hspace{-16pt}
\sih_m\ts(\d_{\ts1i})&=\th_i\,x_{1\bi}
\hspace{-16pt}
&
\phantom{\textrm{for}}
&
&
&
\hspace{-16pt}
\end{align}
where $i=1\lcd n\ts$.
Note that in the case $\f_m=\so_{2m}$ the element
$\sih_m^{\,2}\in\Hh_m$
acts on $x_{1i}$ and on $\d_{\ts1i}$ as the identity, so that
the action of $\Hh_m$ on $\PD\ts(\CC^{\ts m}\ot\CC^{\ts n})$
factors to an action of the group $\H_m\ts$. But in the case
$\f_m=\sp_{2m}$ the element $\sih_m^{\,2}$
acts on $x_{1i}$ and on $\d_{\ts1i}$ as minus the identity,
because $\th_i\,\th_{\ts\bi}=-1$ in this case. This is
why we use the braid group, rather than the Weyl group $\H_m$
of the simple Lie algebra~$\sp_{2m}\ts$.
Taking the tensor product of the actions of $\Hh_m$
on the algebras $\U(\ts\f_m)$ and $\PD\ts(\CC^{\ts m}\ot\CC^{\ts n})\ts$,
we get an action of $\Hh_m$ by automorphisms
of the algebra $\B_m=\U(\ts\f_m)\ot\ts\PD\ts(\CC^{\ts m}\ot\CC^{\ts n})\ts$.

\begin{lemma*}
\label{lemma41}
{\rm\,\,(i)}
The map $\zeta_n:\U(\ts\f_m)\to\PD\ts(\CC^{\ts m}\ot\CC^{\ts n})$
is\/ $\Hh_m\ts$-equivariant.
\\
{\rm(ii)}
The action of\/ $\Hh_m$ on $\B_m$ leaves invariant
any element of the image of\/ $\X(\g_n)$
under the homomorphism $\be_m$.
\end{lemma*}

\begin{proof}
Let us employ the elements $p_{ci}$ and $q_{ci}$
of the algebra $\PD\ts(\CC^{\ts m}\ot\CC^{\ts n})\ts$,
introduced immediately after stating Proposition~\ref{xb}.
In terms of these elements, the action of $\Hh_m$ on the algebra
$\PD\ts(\CC^{\ts m}\ot\CC^{\ts n})$ can be described by setting
\begin{equation*}
\begin{aligned}
\sih\ts(\ts p_{ci})&=p_{\ts\sib(c)\ts i}
&
\quad\textrm{and}\quad&
&
\sih\ts(q_{ci})&=q_{\ts\sib(c)\ts i}
\quad\quad\textrm{for}\quad\quad
\si\in\Sym_m\,,
\\
\sih_m\ts(\ts p_{ci})&=(\ts\pm\ts1\ts)^{\ts\de_{c1}}p_{\ts\sib_m(c)i}
&
\quad\textrm{and}\quad&
&
\sih_m\ts(q_{ci})&=(\ts\pm\ts1\ts)^{\ts\de_{c1}}q_{\ts\sib_m(c)i}
\end{aligned}
\end{equation*}
where $c=-m\lcd-1\com1\lcd m\ts$.
Part (i) follows
by comparing our definition of the action of $\Hh_m$ on $\f_m$
with the description \eqref{ganpq} of the homomorphism~$\zeta_n\ts$.
Part~(ii) follows similarly, using
the description \eqref{fhompq} of 
$\be_m\ts$.
\qed
\end{proof}

Consider the Cartan subalgebra $\h$ from
the triangular decomposition \eqref{tridec}.
In the notation of this section,
our chosen basis of $\h$ is $(\ts F_{-\ab,-\ab}\,|\,a=1\lcd m\ts)\ts$.
Now let $(\ts\ep_a\,|\,a=1\lcd m\ts)\subset\h^\ast$ be the dual basis,
so that $\ep_b\ts(F_{-\ab,-\ab})=\de_{ab}\ts$.
For $c<0$ put $\ep_c=-\ts\ep_{-c}\ts$.
Thus the element $\ep_c\in\h^\ast$ is defined for every index
$c=-\ts m\lcd-1\com1\lcd m\ts$.

Consider the root system of the Lie algebra $\f_m$ in $\h^\ast\ts$. Put
\begin{equation*}
\eta_{\ts a}=\ep_a-\ep_{a+1}
\ \quad\textrm{for}\ \quad
a=1\lcd m-1\ts.
\end{equation*}
Also put $\eta_m=\ep_{m-1}+\ep_m$ in the case $\f_m=\so_{2m}\ts$,
and $\eta_m=2\ts\ep_m$ in the case $\f_m=\sp_{2m}\ts$.
Then $\eta_1\lcd\eta_m$ are the \textit{simple roots\/} of $\f_m\ts$.
Denote by $\De^+$
the set of \textit{positive roots\/} of $\f_m\ts$.
These are the weights $\ep_a-\ep_b$ and
$\ep_a+\ep_b$ where $1\le a<b\le m$ in the case $\f_m=\so_{2m}\ts$,
and the same weights together with $2\ts\ep_a$
where $1\le a\le m$ in the case $\f_m=\sp_{2m}\ts$.
We assume that in the case $\f_m=\so_2$ the root system of $\f_m$ is empty.
Let $\rho$ be halfsum of positive roots of $\f_m\ts$,
so that its sequence of labels
$(\ts\rho_1\lcd\rho_m\ts)$ is
$(\ts m-1\lcd0\ts)$ in the case $\f_m=\so_{2m}\ts$, and is
$(\ts m\lcd1\ts)$ in the case $\f_m=\sp_{2m}\ts$.
For each $a=1,...,m-1$ put
\begin{equation}
\label{Fc}
E_a=F_{-\ab,-\overline{a+1}}\,,\quad
F_a=F_{-\overline{a+1},-\ab}\,,\quad
H_a= F_{-\ab,-\ab}-F_{-\overline{a+1},-\overline{a+1}}\,\ts.
\end{equation}
Put
\begin{equation}
\label{Fmm}
E_m=F_{-\overline{m-1},\overline{m}}\,,\quad
F_m=F_{\overline{m},-\overline{m-1}}\,,\quad
H_m=F_{-\overline{m-1},-\overline{m-1}}+F_{-\overline{m},-\overline{m}}
\end{equation}
in the case $\f_m=\so_{2m}$ with $m>1$.
In the case when $\f_m=\sp_{2m}\ts$, put
\begin{equation}
\label{Fm}
E_m=F_{-\overline{m},\overline{m}}\,/\ts2\ts,\quad
F_m=F_{\overline{m},-\overline{m}}\,/\ts2\ts,\quad
H_m=F_{-\overline{m},-\overline{m}}\,.
\end{equation}
For every possible index
$a$ the three elements $E_a\ts,F_a\ts,H_a$ of the Lie algebra $\f_m$
span a subalgebra isomorphic to $\mathfrak{sl}_2\ts$.
They satisfy the commutation relations
\begin{equation}
\label{sltwo}
[\ts E_a\com F_a\ts]=H_a\ts,
\quad
[\ts H_a\com E_a\ts]=2\ts E_a\ts,
\quad
[\ts H_a\com F_a\ts]=-\ts2\ts F_a\ts.
\end{equation}

So far we denoted by $\B_m$ the associative algebra
$\U(\ts\f_m)\ot\ts\PD\ts(\CC^{\ts m}\ot\CC^{\ts n})\ts$.
Let us now use a different presentation
of the same algebra. Namely, from now until the end of the next section,
on we will regard $\B_m$
as the associative algebra generated by the algebras $\U(\ts\f_m)$
and $\PD\ts(\CC^{\ts m}\ot\CC^{\ts n})$ with the cross relations
\begin{equation}
\label{defar}
[\ts X\com Y\ts]=[\ts\zeta_n(X)\com Y\ts]
\end{equation}
for any $X\in\f_m$ and $Y\in\PD\ts(\CC^{\ts m}\ot\CC^{\ts n})\ts$.
The brackets at the left hand side of the relation
\eqref{defar} denote the commutator in $\B_m\ts$,
while the brackets at the right hand side denote the commutator in
the algebra $\PD\ts(\CC^{\ts m}\ot\CC^{\ts n})$ embedded to $\B_m\ts$.
In particular, we will regard
$\U(\ts\f_m)$ as a subalgebra of $\B_m\ts$.
An isomorphism of this $\B_m$ with the
tensor product $\U(\ts\f_m)\ot\ts\PD\ts(\CC^{\ts m}\ot\CC^{\ts n})$ can be
defined by mapping the elements
$X\in\f_m$ and $Y\in\PD\ts(\CC^{\ts m}\ot\CC^{\ts n})$ of $\B_m$
respectively to the elements
$$
X\ot1+1\ot\zeta_n(X)
\quad\textrm{and}\quad
1\ot Y
$$
of $\U(\ts\f_m)\ot\ts\PD\ts(\CC^{\ts m}\ot\CC^{\ts n})\ts$.
Here we use \eqref{gan}.
The action of the braid group $\Hh_m$ on $\B_m$ is defined
via its isomorphism of $\B_m$ with 
$\U(\ts\f_m)\ot\ts\PD\ts(\CC^{\ts m}\ot\CC^{\ts n})\ts$.
Since the map $\zeta_n$ is
$\Hh_m\ts$-equivariant, the same action of
$\Hh_m$ is obtained by extending
the actions of $\Hh_m$ from the subalgebras
$\U(\ts\f_m)$ and $\PD\ts(\CC^{\ts m}\ot\CC^{\ts n})$ to $\B_m\ts$.

Now consider the following two sets of elements
of the algebra $\U(\h)\subset\U(\ts\f_m)\ts$:
\begin{gather}
\label{denset1}
\{\,F_{aa}-F_{\ts bb}+z\ts,\ F_{aa}+F_{\ts bb}+z
\ |\
1\leqslant a<b\leqslant m\ts,\ z\in\ZZ\,\}\,,
\\[2pt]
\label{denset2}
\{\,F_{aa}+z
\ |\
1\leqslant a\leqslant m\ts,\ z\in\ZZ\,\}\,.
\end{gather}
In the case $\f_m=\so_{2m}\ts$,
denote by $\Uhb$ the ring of fractions of the commutative
algebra $\U(\h)$ relative to the set of denominators \eqref{denset1} .
In the case $\f_m=\sp_{2m}\ts$,
denote by $\Uhb$ the ring of fractions of $\U(\h)$
relative to the union of sets \eqref{denset1} and \eqref{denset2}.
The elements of the ring $\Uhb$ can also
be regarded as rational functions on the vector space
$\h^\ast\ts$. The elements of the subalgebra $\U(\h)\subset\Uhb$
are then regarded as polynomial functions~on~$\h^\ast$.

Denote by $\Bb_m$ the ring of fractions of $\B_m$
relative to the same set of denominators as was used to
define the ring of fractions $\Uhb$. But now we regard these
denominators as elements of $\B_m$ using the
embedding of $\h\subset\f_m$ into $\B_m\ts$. The ring $\Bb_m$
is defined due to the following relations
in $\B_m\ts$. For $c<0$ put $\ep_c=-\ts\ep_{-c}\ts$.
Thus the element $\ep_c\in\h^\ast$ is defined for every
$c=-m\lcd-1\com1\lcd m\ts$. Then for any element $H\in\h$ we have
\begin{gather*}
[\ts H\ts,F_{cd}\ts]=(\ts\ep_{\ts\db}-\ep_{\ts\cb}\ts)(H)\ts F_{cd}
\quad\text{for}\quad
c\com d=-m\lcd-1\com1\lcd m\ts;
\\[2pt]
[\ts H\ts,x_{ci}\ts]=
-\ts\ep_{\ts\cb}\ts(H)\,x_{ci}
\quad\text{and}\quad
[\ts H\ts,\d_{\ts ci}\ts]=
\ep_{\ts\cb}\ts(H)\,\d_{\ts ci}
\quad\text{for}\quad
c=1\lcd m\ts.
\end{gather*}
So the ring $\B_m$ obeys the Ore condition relative
to our set of denominators.
Using left multiplication by elements of
$\,\overline{\!\U(\h)\!\!\!}\,\,\,$,
the ring 
$\Bb_m$ becomes a module~of~$\Uhb\ts$.

The ring $\Bb_m$ is also an associative algebra over $\CC\ts$.
The action of the braid group $\Hh_m$ on $\B_m$
preserves the set of denominators, so that $\Hh_m$ also
acts by automorphisms of the algebra $\Bb_m\ts$.
Using the elements
\eqref{Fc} and \eqref{Fmm} when $\f_m=\so_{2m}\ts$,
or the elements
\eqref{Fc} and \eqref{Fm} when $\f_m=\sp_{2m}\ts$,
for every simple root $\eta_{\ts a}$ of $\f_m$
define a linear map
$$
\xi_{\ts a}:\B_m\to\Bb_m
$$
by setting
\begin{equation}
\label{q1}
\xi_{\ts a}(Y)=
Y+\,\sum_{s=1}^\infty\,\,
(\ts s\ts!\,H_a^{\ts(s)}\ts)^{-1}\ts E_a^{\ts s}\,
\widehat{F}_a^{\ts s}(\ts Y)
\end{equation}
where
$$
H_a^{\ts(s)}=(H_a)(H_a-1)\cdots(H_a-s+1)
$$
and $\widehat{F}_a$ is the operator of adjoint action
corresponding to the element $F_a\in\B_m\ts$,
$$
\widehat{F}_a(\ts Y)=[\ts F_a\ts\com Y\ts]\ts.
$$
For a given element $Y\in\B_m$ only finitely many terms of the sum
\eqref{q1} differ from zero. In the case $\f_m=\so_{2}$ there are
no roots of $\f_m\ts$, and no corresponding operators $\B_m\to\Bb_m\ts$.
On the other hand, in the case when $\f_m=\so_{2m}$ with $m>1\ts$,
by \eqref{simact}
\begin{equation*}
\xi_m\,\sih_m=\sih_m\,\xi_{m-1}\ts,
\end{equation*}
because
$$
\sih_m:\,
E_{m-1}\mapsto E_m\ts,\
F_{m-1}\mapsto F_m\ts,\
H_{m-1}\mapsto H_m\ts.
\vspace{2pt}
$$

Let $\J$ and $\Jb$ be the right ideals of algebras $\B_m$ and $\Bb_m$
respectively, generated by all elements of the subalgebra $\n\subset\f_m\ts$.
The following two properties of the linear operator $\xi_{\ts a}$
go back to \cite[Section~2]{Z1}.
For any elements $X\in\h$ and $Y\in\B_m\ts$, 
\begin{gather}
\label{q11}
\xi_{\ts a}(X\ts Y)
\,\in\,
(\ts X+\eta_{\ts a}\ts(X))\,\ts\xi_{\ts a}(\ts Y)\ts+\ts\Jb\ts,
\\[2pt]
\nonumber
\xi_{\ts a}(\ts Y X)
\,\in\,
\,\xi_{\ts a}(\ts Y)\ts(\ts X+\eta_{\ts a}\ts(X))\ts+\ts\Jb\ts.
\end{gather}
See \cite[Section 3]{KN1} for detailed proofs of these two properties.
The proofs use only the commutation relations \eqref{sltwo},
not the actual form of elements $E_a\ts,F_a\ts,H_a\ts$.

The property \eqref{q11} allows us to define a linear map
$
\bar\xi_{\ts a}:\Bb_m\to\Jb\,\ts\backslash\,\Bb_m
$
by 
\begin{equation*}
\bar\xi_{\ts a}(X\,Y)=Z\,\xi_{\ts a}(\ts Y)\ts+\ts\Jb
\quad\text{for}\quad
X\in\Uhb
\quad\text{and}\quad
Y\in\B_m\ts,
\end{equation*}
where the element $Z\in\Uhb$ is defined by the equality
\begin{equation*}
Z(\mu)=X(\ts\mu+\eta_{\ts a})
\quad\text{for}\quad
\mu\in\h^\ast
\end{equation*}
when both $X$ and $Z$ are regarded as rational functions on $\h^\ast\ts$.
The backslash in $\Jb\,\ts\backslash\,\Bb_m$ indicates that the quotient
is taken relative to a \textit{right\/} ideal of $\Bb_m\ts$.
For the proofs of the next two propositions see \cite[Section 4]{KN3}.

\begin{proposition*}
\label{p2}
For any simple root\/ $\eta_{\ts a}\!$ of\/ $\f_m$ we have
the inclusion\/ 
$
\sih\ts(\ts\Jb\ts)\ts\subset\ts\ker\ts\bar\xi_{\ts a}
$
where $\si=\si_a$ unless\/ $\f_m=\so_{2m}$ and\/ $a=m$,
in which case $\si=\sip_m\ts$.
\end{proposition*}

Recall that $\np$ denotes the nilpotent subalgebra of $\f_m$ spanned by
all the elements $F_{cd}$ with $c<d\ts$. Due to the relation
$F_{cd}=-\,\ep_{cd}\,F_{-d,-c}$ the subalgebra $\np$ is also
spanned by the elements $F_{cd}$ with $c<d$ and $c<0\ts$.
Now for any $a=1\lcd m$ denote by
$\np_{\ts a}$ the vector subspace of $\f_m$ spanned by all
the elements $F_{cd}$ with $c<d$ and $c<0\ts$, except the element $E_a\ts$.
Denote by $\Jp$ the left ideal of $\B_m\ts$, generated by the
elements $X-\zeta_n(X)$ where $X\in\np\ts$. Under the
isomorphism of $\B_m$ with
$\U(\ts\f_m)\ot\ts\PD\ts(\CC^{\ts m}\ot\CC^{\ts n})\ts$,
for any $X\in\f_m$
the difference $X-\zeta_n(X)\in\B_m$ is mapped to the element
\begin{equation}
\label{differ}
X\ot1\in\U(\ts\f_m)\ot1
\subset
\U(\ts\f_m)\ot\ts\PD\ts(\CC^{\ts m}\ot\CC^{\ts n})\ts.
\end{equation}
Let $\Jp_{\ts a}$ be
the left ideal of $\B_m\ts$, generated by the elements
$X-\zeta_n(X)$ with $X\in\np_{\ts a}\ts$, and by the element $E_a\in\B_m\ts$.
Denote $\Jpb=\Uhb\ts\Jp$ and
$\Jpb_{\ts a}=\Uhb\ts\Jp_{\ts a}\ts$.
Then both $\Jpb$ and $\Jpb_{\ts a}$ are left ideals of the algebra $\Bb_m\ts$.

\begin{proposition*}
\label{prop3N}
For any simple root\/ $\eta_{\ts a}\!$ of\/ $\f_m$ we have\/
$
\bar\xi_{\ts a}\ts(\,\sih\,(\ts\Jpb_{\ts a}\ts))
\ts\subset\ts\Jpb\ns+\Jb
$
where $\si=\si_a$ unless $\f_m=\so_{2m}$ and $a=m\ts$,
in which case $\si=\sip_m\ts$.
\end{proposition*}

Proposition \ref{p2}
allows us to define for any simple root $\eta_{\ts a}$ 
a linear map
\begin{equation*}
\xic_{\ts a}:\,\Jb\,\ts\backslash\,\Bb_m\to\Jb\,\ts\backslash\,\Bb_m
\end{equation*}
as the composition $\bar\xi_{\ts a}\,\sih$
applied to the elements of $\Bb_m$
taken modulo $\Jb\ts$. Here the simple reflection
$\si\in\H_m$ is chosen as in Proposition \ref{p2}.
In their present form, the operators $\xic_{\ts 1}\lcd\xic_{\ts m}$
on the vector space
$\Jb\,\ts\backslash\,\Bb_m$ have been defined in \cite{KO}.
We call them the \textit{Zhelobenko operators}.
For the proof of the next proposition see \cite[Sections 4 and 6]{KO}.

\begin{proposition*}
\label{p3}
The Zhelobenko operators
satisfy the braid relations corresponding to the Lie algebra\/ $\f_m\ts$.
Namely, in the case\/ $\f_m=\sp_{2m}$ we have
\begin{align}
\label{braid1}
\xic_{\ts a}\,\xic_{\ts a+1}\,\xic_{\ts a}
&\,=\,
\xic_{\ts a+1}\,\xic_{\ts a}\,\xic_{\ts a+1}
\!\!\!\quad\quad\textit{for}\ \quad
a=1\lcd m-2\ts;
\\
\label{braid2}
\xic_{\ts a}\,\xic_{\ts b}
&\,=\,
\xic_{\ts b}\,\xic_{\ts a}
\hspace{38pt}
\quad\textit{for}\ \quad
a=1\lcd b-2\ts;
\\[2pt]
\nonumber
\xic_{\ts m-1}\,\xic_{\ts m}\,\xic_{\ts m-1}\,\xic_{\ts m}
&\,=\,
\xic_{\ts m}\,\xic_{\ts m-1}\,\xic_{\ts m}\,\xic_{\ts m-1}\ts.
\end{align}
In the case when\/ $\f_m=\so_{2m}$ and\/ $m>1$, we have the same
relations \eqref{braid1} and\/ \eqref{braid2}
between\/ $\xic_{\ts1}\lcd\xic_{\ts m-1}$ as in the case\/ $\f_m=\sp_{2m}$
above, and also the relations 
\begin{align}
\label{xiam}
\xic_a\,\xic_{\ts m}
&\,=\,
\xic_{\ts m}\,\xic_{\ts a}
\ \quad\textit{for}\ \quad
a=1\lcd m-3\com m-1\ts;
\\[2pt]
\nonumber
\xic_{\ts m-2}\,\xic_{\ts m}\,\xic_{\ts m-2}
&\,=\,
\xic_{\ts m}\,\xic_{\ts m-2}\,\xic_{\ts m}\ts.
\end{align}
\end{proposition*}

For $\f_m=\sp_{2m}\ts$, by using
any reduced decomposition of an element $\si\in\H_m$ in terms
of the involutions $\si_1\lcd\si_m\ts$, we can
define a linear operator
\begin{equation}
\label{xis}
\xic_{\ts\si}:\,\Jb\,\ts\backslash\,\Bb_m\to\Jb\,\ts\backslash\,\Bb_m
\end{equation}
in the usual way, like in \eqref{sih}.
This definition of $\xic_{\ts\si}$
is independent of the choice of a reduced decomposition of $\si$
due to Proposition \ref{p3}.

When $\f_m=\sp_{2m}\ts$, the number
of the factors $\si_1\lcd\si_m$ in any reduced decomposition
$\si\in\H_m$ will be denoted $\ell\ts(\si)\ts$. This number
is also independent of the choice
of the decomposition, and is equal to
the number of elements in the set
\begin{equation}
\label{posneg}
\De_{\ts\si}=\{\ts\eta\in\De^+\,|\,\si\ts(\eta)\notin\De^+\,\}
\end{equation}
where $\De^+$ denotes the set of positive roots of
the Lie algebra $\sp_{2m}\ts$.

Now suppose that $\f_m=\so_{2m}\ts$.
Then by using any reduced decomposition in terms of
$\si_1\lcd\si_{m-1}\com\sip_m\ts$,
we can define a linear operator \eqref{xis} for every element
$\si\in\Hp_m\ts$. Again, this definition
is independent of the choice of a reduced decomposition of $\si$
due to Proposition \ref{p3}.
It turns out that in this case
we can extend the definition of the operator \eqref{xis}
to any element $\si\in\H_m\ts$, where $m\geqslant1$. Note that
in this case the action of the element
$\sih_m$ on $\Bb_m$ preserves the ideal $\Jb\ts$,
and therefore induces a linear operator
on the quotient vector space $\Jb\,\ts\backslash\,\Bb_m\ts$.
This operator will be again denoted by $\sih_m\ts$.
The extension of the definition of the operators \eqref{xis} to
$\si\in\H_m$ is based on the next lemma,
which has been proved in \cite[Section 4]{KN3}.

\begin{lemma*}
\label{lemma45}
When $\f_m=\so_{2m}$ and\/ $m>1\ts$, the operators\/
$\xic_{\ts1}\lcd\xic_{\ts m-1}\com\sih_m$ on $\Jb\,\ts\backslash\,\Bb_m$
satisfy the same relations, as the\/ $m$
generators of the braid group $\Hh_m$ respectively.
Then we also have the relation
\begin{equation}
\label{ximh}
\xic_{\ts m}=
\sih_m\,\xic_{\ts m-1}\,\sih_m\ts.
\end{equation}
\end{lemma*}

Now for $\f_m=\so_{2m}$ with any $m\geqslant1$,
take any decomposition of an element $\si\in\H_m$ in terms
of the involutions $\si_1\lcd\si_m$ such that the
number of occurencies of $\si_1\lcd\si_{m-1}$ in the decomposion
is minimal possible. For $\f_m=\so_{2m}$ the symbol
$\ell\ts(\si)$ will denote this minimal number.
Note that unlike for $\f_m=\sp_{2m}\ts$,
here we do not count the occurencies of $\si_m$ in the
decomposition. All the decompositions of $\si\in\H_m$ with
the minimal number of occurencies of $\si_1\lcd\si_{m-1}$
can be obtained from each other by using the braid
relations between $\si_1\lcd\si_m\in\H_m$ along with
the relation $\si_m^{\ts2}=1\ts$.

By substituting the operators
$\xic_{\ts1}\lcd\xic_{\ts m-1}\com\sih_m$ on
$\Jb\,\ts\backslash\,\Bb_m$ for
involutions $\si_1\lcd\si_m$ in such a decomposition of
$\si\in\H_m\ts$, we obtain another operator on
$\Jb\,\ts\backslash\,\Bb_m\ts$. The latter operator
does not depend on the choice of a decomposition
because of the first statement of Lemma \ref{lemma45},
and because the operator $\sih_m^{\,2}$ on
the vector space $\Jb\,\ts\backslash\,\Bb_m$ is the identity
in the case $\f_m=\so_{2m}$ considered here.
Moreover for $\si\in\Hp_m\subset\H_m\ts$,
the operator on $\Jb\,\ts\backslash\,\Bb_m$
obtained by the latter substitution coincides with the operator \eqref{xis}.
Indeed, for $\f_m=\so_{2m}$
the operator \eqref{xis} has been defined by
substituting the Zhelobenko operators
$\xic_{\ts1}\lcd\xic_{\ts m-1}\com\xic_{\ts m}$
for $\si_1\lcd\si_{m-1}\com\sip_m$
in any reduced decomposition of $\si\in\Hp_m\ts$.
The coincidence of the two operators for $\si\in\Hp_m$
now follows from the relation \eqref{ximh}.
Thus we have extended the definition of
the operator \eqref{xis} from $\si\in\Hp_m$ to all $\si\in\H_m\ts$.

Note that for $\f_m=\so_{2m}$ and
$\si\in\Hp_m\ts$, the number $\ell\ts(\si)$ is equal to
the length of a reduced decomposition of $\si$ in terms of
$\si_1\lcd\si_{m-1}\com\sip_m\ts$. Thus
we have also extended the standard length function
from the Weyl group $\Hp_m$ of $\so_{2m}$ to
the hyperoctahedral group $\H_m\ts$.
Moreover for any $\si\in\H_m\ts$,
not only for $\si\in\Hp_m\ts$,
the number $\ell\ts(\si)$ equals
the number of elements in the set \eqref{posneg}, where
$\De^+$ is the set of positive roots of $\so_{2m}\ts$.

>From now we shall consider
$\f_m=\so_{2m}$ and $\f_m=\sp_{2m}$ simultaneously,
and will work with the operators \eqref{xis} for all elements
$\si\in\H_m\ts$.
In particular, in the case $\f_m=\so_{2m}$
we will assume that the operator \eqref{xis} with $\si=\si_m$
acts as~$\sih_m\ts$.

The restriction of
the action \eqref{siact}\com\ts\eqref{simact} of the braid group $\Hh_m$
on $\f_m$ to the Cartan subalgebra $\h$ factors
to an action of the hyperoctahedral group $\H_m\ts$. This is the
standard action of the Weyl group of $\f_m=\sp_{2m}\ts$.
The resulting action of the subgroup $\Hp_m\subset\H_m$ on $\h$
is the standard action of the Weyl group of $\f_m=\so_{2m}\ts$.
The group $\H_m$ also acts on the dual vector space $\h^\ast\ts$, so that
$\si\ts(\ep_c)=\ep_{\si(c)}$ for any $\si\in\H_m$ and any
$c=-\ts m\lcd-1\com1\lcd m\ts$. Unlike in \eqref{sib}, here we use
the natural action of the group $\H_m$ by permutations
of $-\ts m\lcd-1\com1\lcd m\ts$. Thus $\si_a\in\H_m$
with $1\le a<m$ exchanges $a\com a+1$ and also exchanges
$-a\com-a-1\ts$ while $\si_m\in\H_m$ exchanges $m\com-m\ts$.
Note that we always have $\si\ts(-c)=-\si\ts(c)\ts$.
If we identify each weight $\mu\in\h^\ast$ with
the sequence $(\ts\mu_1\lcd\mu_m)$ of its labels, then
\begin{align*}
\si:(\ts\mu_1\lcd\mu_m)
&\,\mapsto\,
(\,\mu_{\ts\si^{-1}(1)}\lcd\mu_{\ts\si^{-1}(m)}\ts)
\quad\text{for}\quad
\si\in\Sym_m\ts,
\\
\si_m:(\ts\mu_1\lcd\mu_m)
&\,\mapsto\,
(\,\mu_{1}\lcd\mu_{m-1},-\mu_{m}\ts)\ts.
\end{align*}

The \textit{shifted\/} action of the
group $\H_m$ on the set $\h^\ast$ is defined by the assignment
\begin{equation*}
\mu\,\mapsto\,\si\circ\mu=\si\ts(\mu+\rho)-\rho
\ \quad\textrm{for}\quad\
\si\in\H_m\ts.
\end{equation*}
By regarding the elements of the commutative algebra $\Uhb$
as rational functions on the vector space $\h^\ast$
we can also define an action of the group $\H_m$ on this \textrm{algebra:}
\begin{equation}
(\si\circ X)(\mu)=X(\si^{-1}\ns\circ\mu)
\,\quad\textrm{for}\quad
X\in\Uhb\ts.
\end{equation}
The next proposition has been also proved in \cite[Section 4]{KN3}.

\begin{proposition*}
\label{saction}
For any $\si\in\H_m\ts$, $X\in\Uhb$ and
$Y\in\Jb\,\ts\backslash\,\Bb_m$ we have the relations
\begin{align}
\label{q121}
\xic_{\ts\si}(X\ts Y)
&\,=\,
(\ts\si\ts\circ X)\,\ts\xic_{\ts\si}(\ts Y)\ts,
\\
\nonumber
\xic_{\ts\si}(\ts YX)
&\,=\,
\,\xic_{\ts\si}(\ts Y)\ts(\ts\si\ts\circ X)\ts.
\end{align}
\end{proposition*}


\section*{\normalsize 5. Intertwining operators}
\setcounter{section}{5}
\setcounter{equation}{0}
\setcounter{theorem*}{0}

Let $\de=(\ts\de_1\lcd\de_m)$ be any sequence of $m$
elements from $\{1\com-1\}\ts$. The hyperoctahedral
group $\H_m$ acts on the set of all these sequences naturally,
so that the generator $\si_a\in\H_m$ with $a<m$ acts on $\de$ as
the transposition
of $\de_a$ and $\de_{a+1}\ts$, while the generator $\si_m\in\H_m$
changes the sign of $\de_m\ts$. Let
$\de_+=(1\lcd1)$ be
the sequence of $m$ elements $1$.
Given any sequence $\de\ts$, take the composition of the automorphisms
of the ring $\PD\ts(\CC^{\ts m}\ot\CC^{\ts n})\ts$,
\begin{equation}
\label{compfour}
x_{\ts\ab i}\mapsto\th_i\,\d_{\,\ab\ts\bi}
\quad\text{and}\quad
\d_{\,\ab i}\mapsto\th_i\,x_{\ts\ab\ts\bi}
\quad\text{whenever}\quad
\de_a=-1\ts.
\end{equation}
Here $a\ge1$ and $i=1\lcd n\ts$.
Let us denote by $\varpi$ this composition.
In particular, the automorphism $\varpi$
corresponding to $\de=(1\lcd1\com-1)$ coincides
with the action of $\sih_m$ on $\PD\ts(\CC^{\ts m}\ot\CC^{\ts n})\ts$,
see \eqref{fourier}. In the case $\f_m=\so_{2m}\ts$,
the automorphism $\varpi$ is involutive for any
$\de\ts$. But in the case $\f_m=\sp_{2m}$ the square $\varpi^2$ maps
$$
x_{\ts\ab i}\mapsto-\ts x_{\,\ab i}
\quad\text{and}\quad
\d_{\,\ab i}\mapsto-\ts\d_{\,\ab i}
\quad\text{whenever}\quad
\de_a=-1\ts.
$$

For any $\f_m\ts$-module $V$, the action of
$\X(\g_n)$ on 
$\F_m(V)=V\ot\P\ts(\CC^{\ts m}\ot\CC^{\ts n})$
is defined by the homomorphism
$\be_m:\X(\g_n)\to\U(\ts\f_m)\ot\ts\PD\ts(\CC^{\ts m}\ot\CC^{\ts n})\ts$,
see Proposition \ref{xb}. Further,
the action of Lie algebra $\f_m$ on the second tensor factor
$\P\ts(\CC^{\ts m}\ot\CC^{\ts n})$ of $\F_m(V)$
is defined by means of homomorphism
$\zeta_n:\U(\ts\f_m)\to\PD\ts(\CC^{\ts m}\ot\CC^{\ts n})\ts$, see
definition \eqref{gan}.
Here any element of the ring
$\PD\ts(\CC^{\ts m}\ot\CC^{\ts n})$ acts on the vector space
$\P\ts(\CC^{\ts m}\ot\CC^{\ts n})$ naturally.
We can modify the latter action,
by making any element $Y\in\PD\ts(\CC^{\ts m}\ot\CC^{\ts n})$
act on $\P\ts(\CC^{\ts m}\ot\CC^{\ts n})$ via
the natural action of $\varpi\ts(Y)\ts$.
Then we get another $\PD\ts(\CC^{\ts m}\ot\CC^{\ts n})\ts$-module,
with the same underlying vector space
$\P\ts(\CC^{\ts m}\ot\CC^{\ts n})$ for every~$\de\ts$.

For any $\f_m\ts$-module $V$, we can now define a bimodule
$\F_\de\ts(V)$ of $\f_m$ and $\X(\g_n)\ts$.
Its underlying vector space is the same
$V\ot\P\ts(\CC^{\ts m}\ot\CC^{\ts n})$ for every~$\de\ts$.
The action of $\X(\g_n)$ on $\F_\de\ts(V)$
is defined by pushing the homomorphism $\be_m$
forward through the automorphism $\varpi$,
applied to $\PD\ts(\CC^{\ts m}\ot\CC^{\ts n})$
as to the second tensor factor of the target of $\be_m\ts$.
The action of $\f_m$ on $\F_\de\ts(V)$
is also defined by pushing the homomorphism $\zeta_n$
forward through the automorphism $\varpi\ts$.
Thus the actions of $\X(\g_n)$ and $\f_m$ on the bimodule $\F_\de\ts(V)$
are respectively determined by the compositions
of the homomorphisms
\begin{align*}
\X(\g_n) \underset{\be_m}\longrightarrow
\U(\ts\f_m)\ot\ts\PD\ts(\CC^{\ts m}\ot\CC^{\ts n})
\underset{1\ot\tts\varpi\ts}\longrightarrow
\U(\ts\f_m)\ot\ts\PD\ts(\CC^{\ts m}\ot\CC^{\ts n})\,,&
\\
\U(\ts\f_m) \underset{1\ot\tts\zeta_n}\longrightarrow
\U(\ts\f_m)\ot\ts\PD\ts(\CC^{\ts m}\ot\CC^{\ts n})
\underset{1\ot\tts\varpi\ts}\longrightarrow
\U(\ts\f_m)\ot\ts\PD\ts(\CC^{\ts m}\ot\CC^{\ts n})\,.&
\end{align*}
Note that here we have
$$
\F_m(V)=\F_{\ts\de_+}(V)\,.
$$

Now let $\mu\in\h^\ast$ be any weight of $\f_m\ts$, such that
\begin{equation}
\label{notinso}
\mu_a-\mu_b\notin\ZZ
\quad\text{and}\quad
\mu_a+\mu_b\notin\ZZ
\quad\text{whenever}\quad
1\le a<b\le m\ts.
\end{equation}
In the case $\f_m=\sp_{2m}$ also suppose that,
in addition to \eqref{notinso},
\begin{equation}
\label{notinsp}
2\ts\mu_a\notin\ZZ
\quad\text{whenever}\quad
1\le a\le m\ts.
\end{equation}
We shall now proceed to show how for every element $\si\in\H_m\ts$,
the Zhelobenko operator \eqref{xis}
determines an $\X(\g_n)$-intertwining operator
\begin{equation}
\label{distoper}
\F_{m}\ts(\ts M_{\ts\mu})_{\ts\n}
\,\to\,
\F_{\ts\de}\ts(\ts M_{\ts\si\,\circ\ts\mu})_{\ts\n}
\ \quad\text{where}\quad\
\de=\si\ts(\de_+)\ts.
\end{equation}

In this section we keep regarding $\B_m$
as the associative algebra generated by
$\U(\ts\f_m)$ and $\PD\ts(\CC^{\ts m}\ot\CC^{\ts n})$ with the cross relations
\eqref{defar}.
Let $\I_{\ts\de}$ be the left ideal of algebra $\B_m$
generated by the elements $x_{\ts\ab\ts i}$ with $\de_a=-1\ts$,
and the elements $\d_{\ts\ab\ts i}$ with $\de_a=1\ts$.
Here $a=1\lcd m$ and $i=1\lcd m\ts$.
Note that in terms of the elements $q_{\ts ci}$
introduced immediately after stating Proposition~\ref{xb},
the left ideal $\I_{\ts\de}$ is generated by the elements
$q_{\ts-\de_a\ab,i}$ where again $a=1\lcd m$ and $i=1\lcd m\ts$.
In particular, the ideal $\I_{\ts\de_+}$ is generated by
all the left derivations $\d_{\ts ai}\ts$.
Let $\Ib_{\ts\de}$ be the left ideal of 
$\Bb_m$
generated by the same elements as the ideal of $\I_{\ts\de}$ of $\B_m\ts$.

Consider the image of the ideal $\Ib_{\ts\de}$
in the quotient space $\Jb\,\ts\backslash\,\Bb_m\ts$,
that is the subspace $\Jb\,\ts\backslash\,(\ts\Ib_{\ts\de}+\Jb\ts)$
in the quotient space $\Jb\,\ts\backslash\,\Bb_m\ts$.
The image will be occasionally denoted by the same symbol
$\Ib_{\ts\de}\ts$. In the context of the
next proposition, this 
should cause no confusion.

\begin{proposition*}
\label{propznak}
For any $\si\in\H_m$ the operator $\xic_{\ts\si}$ maps the subspace\/
$\I_{\ts\de_+}$ to\/ 
$\ts\I_{\ts\si(\de_+)}\ts$.
\end{proposition*}

\begin{proof}
For any $a=1\lcd m-1$ consider the operator $\widehat{F}_a$
corresponding to the element $F_a\in\B_m\ts$. By 
\eqref{Fc} and \eqref{gan},\eqref{defar} for any
$Y\in\PD\ts(\CC^{\ts m}\ot\CC^{\ts n})$ we have
$$
\widehat{F}_a(\ts Y)=\,
-\,\sum_{k=1}^n\,\,
[\ts x_{\ts\ab\ts k}\,\d_{\,\ts\overline{a+1}\ts k}\,,Y\ts]\ts.
$$
Similarly, in the case $\f_m=\sp_{2m}$ by 
\eqref{Fm} for any $Y\in\PD\ts(\CC^{\ts m}\ot\CC^{\ts n})$ we have
$$
\widehat{F}_m(\ts Y)=\,
\sum_{k=1}^n\,\,
[\ts x_{\ts\overline{m}\ts\bk}\,x_{\ts\overline{m}\ts k}\,,Y\ts]\ts/\ts2\ts.
$$
In the case $\f_m=\so_{2m}$ we do not need to consider
the operator $\widehat{F}_m\ts$, because in this case
the operator \eqref{xis} corresponding to $\si=\si_m$
acts on $\Jb\,\ts\backslash\,\Bb_m$ as $\sih_m$ by our definition.

The above description of the action of $\widehat{F}_a$ with $a<m$ on
$\PD\ts(\CC^{\ts m}\ot\CC^{\ts n})$ shows that this action
preserves each of the two $2n$ dimensional subspaces, spanned by the vectors
\begin{gather}
\label{xxideal}
x_{\ts\ab\ts i}
\quad\text{and}\quad
x_{\ts\overline{a+1}\ts i}
\quad\text{where}\quad
i=1\lcd n\ts;
\\[2pt]
\label{ddideal}
\d_{\,\ab\ts i}
\quad\text{and}\quad
\d_{\,\overline{a+1}\ts i}
\quad\text{where}\quad
i=1\lcd n\ts.
\end{gather}
This action also maps to zero the $2n$ dimensional subspace, spanned by
\begin{equation}
\label{xdideal}
x_{\ts\ab\ts i}
\quad\text{and}\quad
\d_{\,\overline{a+1}\ts i}
\quad\text{where}\quad
i=1\lcd n\ts.
\end{equation}
Therefore for any $\de\ts$, the operator $\bar\xi_a$ with $a<m$ maps
the left ideal $\Ib_{\ts\de}$ of $\Bb_m$ to the image of $\Ib_{\ts\de}$ in
$\Jb\,\ts\backslash\,\Bb_m\ts$, unless $\de_a=1$ and $\de_{a+1}=-1\ts$.
The operator $\xic_{\ts a}$ on
$\Jb\,\ts\backslash\,\Bb_m$ was defined by
taking the composition of 
$\bar\xi_a$ and
$\sih_a\ts$. Hence $\xic_{\ts a}$ with $a<m$ maps the image of
$\Ib_{\ts\de}$ to the image of $\Ib_{\ts\si_a(\de)}\ts$, unless
$\de_a=-1$ and $\de_{a+1}=1\ts$.

In the case $\f_m=\sp_{2m}\ts$,
the action of $\widehat{F}_m$
on the vector space $\PD\ts(\CC^{\ts m}\ot\CC^{\ts n})$
maps to zero the $n$ dimensional subspace spanned by the elements
\begin{equation}
\label{xideal}
x_{\ts\overline{m}\ts i}\ts=\ts x_{1i}
\quad\text{where}\quad
i=1\lcd n\ts.
\end{equation}
Therefore the operator $\bar\xi_m$ maps
the left ideal $\Ib_{\ts\de}$ of $\Bb_m$ to the image of $\Ib_{\ts\de}$ in
$\Jb\,\ts\backslash\,\Bb_m\ts$, unless $\de_m=1$.
Hence the operator $\xic_{\ts m}$ on
$\Jb\,\ts\backslash\,\Bb_m$ maps the image of
$\Ib_{\ts\de}$ 
to the image of $\Ib_{\ts\si_m(\de)}\ts$,
unless $\de_m=-1\ts$. In the case $\f_m=\so_{2m}\ts$,
we just note that
$\sih_m$ maps the image of
$\Ib_{\ts\de}$ in $\Jb\,\ts\backslash\,\Bb_m$ to the image of
$\Ib_{\ts\si_m(\de)}\ts$.

>From now on we will denote
the image of the ideal $\Ib_{\ts\de}$
in the quotient space $\Jb\,\ts\backslash\,\Bb_m\ts$
by the same symbol. Put
$$
\widehat{\de}=\sum_{a=1}^m\,\de_a\ts\ep_a\in\h^\ast\ts.
$$
Then for every $\si\in\H_m$ we have the equality
$\widehat{\si\ts(\de\ts)}=\si\ts(\,\widehat{\de}\,)$
where at the right hand side we use the action
of the group $\H_m$ on $\h^\ast$.
Let $(\ ,\,)$ be the standard bilinear
form on $\h^\ast$, so that the basis of weights $\ep_a$ with $a=1\lcd m$
is orthonormal. The above remarks on the action
of the Zhelobenko operators on $\Ib_{\ts\de}$ can now be restated
as follows:
\begin{align}
\label{rest1}
\textrm{if}\quad(\ts\widehat{\de}\com\ep_a-\ep_{a+1})&\ge0
&
\textrm{then}&
&
&\ts\xic_{\ts a}\ts(\,\Ib_\de\ts)\subset\ts\Ib_{\ts\si_a(\de)}
\,\quad\textrm{for}\quad
a=1\lcd m-1\,;
\\
\label{rest2}
\textrm{if}\quad(\ts\widehat{\de}\com\ep_m)&>0
&
\textrm{then}&
&
&\ns\xic_{\ts m}(\,\Ib_\de\ts)\subset\Ib_{\ts\si_m(\de)}
\quad\textrm{for}\quad
\f_m=\sp_{2m}\ts.
\end{align}

We shall prove Proposition \ref{propznak}
by induction on the length of a reduced decomposition
of $\si\in\H_m$ in terms of $\si_1\lcd\si_m\ts$.
This number has been denoted by $\ell(\si)$ in the case
$\f_m=\sp_{2m}\ts$, but may be different from the number
denoted by $\ell(\si)$ in the case $\f_m=\so_{2m}\ts$.
Recall that in both cases $\ell\ts(\si)$ equals
the number of elements in the set \eqref{posneg},
where $\De^+$ is the set of positive roots of $\f_m\ts$.

If $\si$ is the identity element of $\H_m\ts$,
Proposition \ref{propznak} is tautological.
Suppose that for some $\si\in\H_m\ts$,
$$
\xic_{\ts\si}(\,\Ib_{\ts\de_+}\ts)\subset\ts\Ib_{\ts\si(\de_+)}\ts.
$$
Take $\si_a\in\H_m$ with $1\le a\le m\ts$,
such that $\si_a\ts\si$ has a longer reduced decomposition
in terms of $\si_1\lcd\si_m$ than $\si\ts$.
If $\f_m=\so_{2m}$ and $a=m\ts$, then
$\xic_{\ts\si_m\ts\si}=\sih_m\,\xic_{\ts\si}$ and we need the inclusion
\begin{equation}
\label{obvinc}
\sih_m\ts(\,\Ib_{\ts\si(\de_+)}\ts)\subset\ts\Ib_{\ts\si_m\si(\de_+)}\ts,
\end{equation}
which holds by the definition of the action of $\H_m$
on $\Jb\,\ts\backslash\,\Bb_m\ts$.

We may exclude the case when $\f_m=\so_{2m}$ and $a=m\ts$, and assume that
\begin{equation}
\label{lengthup}
\ell\ts(\ts\si_a\ts\sigma)=\ell\ts(\sigma)+1\ts.
\end{equation}
Firstly, suppose that $a<m$ here. Let us then prove the inclusion
$$
\xic_{\ts a}(\,\Ib_{\ts\si(\de_+)}\ts)\subset\ts\Ib_{\ts\si_a\si(\de_+)}\ts.
$$
By \eqref{rest1}, the latter inclusion will have place if
\begin{equation*}
(\ts\widehat{\si\ts(\de_+)}\ts\com\ts\ep_a-\ep_{a+1})=
(\ts\si\ts(\ts\widehat{\de}_+)\ts\com\ts\ep_a-\ep_{a+1})\ge0\ts.
\end{equation*}
But the condition \eqref{lengthup} for $a<m$ implies that
$\ep_a-\ep_{a+1}\in\sigma\ts(\De^+)\ts$. Indeed, because the root
$\ep_a-\ep_{a+1}$ of $\f_m$ is simple,
$\si_a(\eta)\in\De^+$ for any
$\eta\in\De^+$ such that $\eta\neq\ep_a-\ep_{a+1}\ts$.
Since $\ell\ts(\sigma)$ and $\ell\ts(\ts\si_a\ts\sigma)$
are the numbers of elements in
$\De_{\ts\si}$ and $\De_{\ts\si_a\ts\si}$ respectively,
here $\ep_a-\ep_{a+1}\in\sigma\ts(\De^+)\ts$.
So $\ep_a-\ep_{a+1}=\si\ts(\ts\ep_b-\ep_c)$ where
$1\le b\le m$ and $1\le|c|\le m\ts$. Thus
$$
(\ts\si\ts(\ts\widehat{\de}_+)\ts\com\ts\ep_a-\ep_{a+1})
=
(\ts\si\ts(\ts\widehat{\de}_+)\ts\com\ts\si\ts(\ts\ep_b-\ep_c)\ts)
=
(\ts\widehat{\de}_+\ts\com\ts\ep_b-\ep_c\ts)\ge0\ts.
$$

Now suppose that $a=m\ts$. Here we assume that
$\f_m=\sp_{2m}\ts$. We need the inclusion
$$
\xic_{\ts m}\ts(\,\Ib_{\ts\si(\de_+)}\ts)
\subset\ts
\Ib_{\ts\si_m\si(\de_+)}\ts.
$$
It will have place if
\begin{equation*}
(\ts\widehat{\si\ts(\de_+)}\com\ep_m\ts)
=
(\ts\si\ts(\ts\widehat{\de}_+)\ts\com\ep_m\ts)>0\ts.
\end{equation*}
But the condition \eqref{lengthup} for $a=m$ implies
that $2\ts\ep_m\in\sigma\ts(\De^+)\ts$,
where $\De^+$ is the set of positive roots of $\sp_{2m}\ts$.
Indeed, because the root
$2\ts\ep_m$ of $\sp_{2m}$ is simple,
$\si_m(\eta)\in\De^+$ for any
$\eta\in\De^+$ such that $\eta\neq2\ts\ep_m\ts$.
Since $\ell\ts(\sigma)$ and $\ell\ts(\ts\si_m\ts\sigma)$
are the numbers of elements in
$\De_\si$ and $\De_{\ts\si_m\ts\si}$ respectively,
here $2\ts\ep_m\in\sigma\ts(\De^+)\ts$.
So $\ep_m=\si\ts(\ts\ep_b)$ where $1\le b\le m\ts$. Thus
$$
(\ts\si\ts(\ts\widehat{\de}_+)\ts\com\ep_m\ts)
=
(\ts\si\ts(\ts\widehat{\de}_+)\ts\com\si\ts(\ts\ep_b)\ts)
=
(\ts\widehat{\de}_+\ts\com\ts\ep_b\ts)>0\ts.
\eqno\qed
$$
\end{proof}

\begin{corollary*}
\label{corollary4.6}
For any $\si\in\H_m$ the operator $\xic_{\ts\si}$ on
$\Jb\,\ts\backslash\,\Bb_m$ maps 
$$
\Jb\,\ts\backslash\,(\Jpb+\Ib_{\ts\de_+}+\Jb\ts)
\,\to\ts
\Jb\,\ts\backslash\,(\Jpb+\Ib_{\ts\si(\de_+)}+\Jb\ts)\ts.
$$
\end{corollary*}

\begin{proof}
We will extend the arguments used in the proof of
Proposition \ref{propznak}. In particular, we will again use
the length of a reduced decomposition
of $\si$ in terms of $\si_1\lcd\si_m\ts$.
If $\si$ is the identity element of $\H_m\ts$,
then the required statement is tautological.
Now suppose that for some $\si\in\H_m$
the statement of Corollary \ref{corollary4.6} is true.
Take any simple reflection $\si_a\in\H_m$ with $1\le a\le m\ts$,
such that $\si_a\ts\si$ has a longer reduced decomposition
in terms of $\si_1\lcd\si_m$ than $\si\ts$.
In the case $\f_m=\so_{2m}$ we may assume
that $a<m\ts$, because in that case
the required statement for $\si_m\ts\si$ instead of $\si$
is provided by 
\eqref{obvinc}.

Thus we have the equality
\eqref{lengthup}.
With the above assumption on $a\ts$, we have proved that
\eqref{lengthup} implies 
\begin{equation}
\label{thischeck}
(\ts\widehat{\si\ts(\de_+)}\ts\com\ts\eta_{\ts a})\ge0\ts.
\end{equation}
Here $\eta_a$ is the simple root corresponding to 
$\si_a\ts$. But \eqref{thischeck} implies the equality
\begin{equation}
\label{N10}
\Jpb+\Ib_{\ts\si(\de_+)}=
\Jpb_{\ts a}+\Ib_{\ts\si(\de_+)}
\end{equation}
of left ideals of $\Bb_m\ts$.
Indeed, the left and right hand sides of \eqref{N10} differ by
the elements $Y\ts\zeta_n(E_a)$ where $Y$ ranges over $\Bb_m\ts$.
The condition \eqref{thischeck} implies that
$\zeta_n(E_a)\in\Ib_{\ts\si(\de_+)}\ts$, see
the definition \eqref{gan} and the arguments in
the beginning of proof of Proposition~\ref{propznak}.
Using Proposition \ref{prop3N} and the induction step
from our proof of Proposition~\ref{propznak}, 
$\xic_{\ts a}$ maps
$$
\Jb\,\ts\backslash\,(\Jpb+\Ib_{\ts\si(\de_+)}+\Jb\ts)
=
\Jb\,\ts\backslash\,(\Jpb_{\ts a}+\Ib_{\ts\si(\de_+)}+\Jb\ts)
\,\to\ts
\Jb\,\ts\backslash\,(\Jpb+\Ib_{\ts\si_a\si(\de_+)}+\Jb\ts)\ts.
$$
This makes the induction step of our proof of Corollary \ref{corollary4.6}.
\qed
\end{proof}

Let $\I_{\ts\mu,\de}$ be the left ideal of the algebra $\B_m$
generated by $\I_{\ts\de}+\Jp$ and by the elements
$$
F_{-\ab,-\ab}-\zeta_n\ts(F_{-\ab,-\ab})-\mu_a
\quad\text{where}\quad
a=1\lcd m\ts.
$$
Recall that under the isomorphism of the algebra $\B_m$ with
$\U(\ts\f_m)\ot\ts\PD\ts(\CC^{\ts m}\ot\CC^{\ts n})\ts$,
the difference $X-\zeta_n(X)\in\B_m$ for every $X\in\f_m$
is mapped to the element \eqref{differ}.
Denote by $\Ib_{\ts\mu,\de}$ the subspace
$\Uhb\,\I_{\ts\mu,\de}\ts$ of $\Bb_m\ts$, this
is also a left ideal of 
$\Bb_m\ts$.

\begin{theorem*}
\label{proposition4.5}
For any element $\si\in\H_m$ the operator $\xic_{\ts\si}$ on
$\Jb\,\ts\backslash\,\Bb_m$ maps 
$$
\Jb\,\ts\backslash\,
(\ts\Ib_{\ts\mu,\de_+}+\Jb\ts)
\,\to\ts
\Jb\,\ts\backslash\,
(\ts\Ib_{\ts\si\ts\circ\ts\mu\ts,\ts\si(\de_+)}+\Jb\ts)\ts.
$$
\end{theorem*}

\begin{proof}
Let $\ka$ be a weight of $\f_m$ with the sequence of
labels $(\ts\ka_1\lcd\ka_m)\ts$. Suppose that the weight $\ka$
satisfies the conditions
\eqref{notinso} instead of $\mu\ts$. In the case $\f_m=\sp_{2m}$
we also suppose that $\ka$ satisfies the conditions
\eqref{notinsp} instead of $\mu\ts$. Denote by
$\tilde\I_{\ts\ka,\de}$ be the left ideal of $\Bb_m$
generated by $\I_{\ts\de}+\Jp$ and by the elements
$$
F_{-\ab,-\ab}-\ka_a
\quad\text{where}\quad
a=1\lcd m\ts.
$$
Proposition \ref{saction} and Corollary \ref{corollary4.6}
imply that the operator $\xic_{\ts\si}$ on
$\Jb\,\ts\backslash\,\Bb_m$ maps
$$
\Jb\,\ts\backslash\,
(\ts\tilde\I_{\ts\ka,\de_+}+\Jb\ts)
\,\to\ts
\Jb\,\ts\backslash\,
(\ts\tilde\I_{\ts\si\ts\circ\ts\ka\ts,\ts\si(\de_+)}+\Jb\ts)\ts.
$$

Now choose
\begin{equation}
\label{kappa}
\ka_a=\mu_a+{n}/2
\quad\text{for}\quad
a=1\lcd m\ts.
\end{equation}
Then the conditions on $\ka$ stated in the beginning of this proof
are satisfied. For every $\si\in\H_m$ we shall prove the equality
of left ideals of $\Bb_m\ts$,
\begin{equation}
\label{lastin}
\tilde\I_{\ts\si\ts\circ\ts\ka\ts,\ts\si(\de_+)}=
\Ib_{\ts\si\ts\circ\ts\mu\ts,\ts\si(\de_+)}\ts.
\end{equation}
Theorem \ref{proposition4.5} will then follow.
Denote $\de=\si\ts(\de_+)\ts$. Then by our choice of $\ka$ we have
\begin{equation}
\label{sika}
\si\ts\circ\ts\ka=\si\ts\circ\ts\mu\ts+{n\ts\de}/2
\end{equation}
where the sequence $\de$ is regarded as a weight of $\f_m\ts$,
by identifying the weights with their sequences of labels.
Let $a$ run through $1\lcd m\ts$. If
$\de_a=1$ then by the definition \eqref{gan},
$$
\zeta_n\ts(F_{-\ab,-\ab})-{n}/2\,=\,
-\,\sum_{k=1}^n\,x_{\ts\ab\ts k}\,\d_{\,\ab\ts k}\in\I_{\ts\de}\ts.
$$
If $\de_a=-1$ then the same definition \eqref{gan} implies that
$$
\zeta_n\ts(F_{-\ab,-\ab})+{n}/2\,=\,
\sum_{k=1}^n\,\d_{\,\ab\ts k}\,x_{\ts\ab\ts k}\in\I_{\ts\de}\ts.
$$
Hence the relation \eqref{sika} implies the equality \eqref{lastin}.
\qed
\end{proof}

Consider the quotient vector space
$\B_m\ts/\,\I_{\ts\mu,\de}$
for any sequence $\de\ts$.
The algebra $\U(\ts\f_m)$
acts on this quotient via left multiplication,
being regarded as a subalgebra of $\B_m\ts$. The
algebra $\X(\g_n)$ also acts on this quotient via left multiplication,
using the homomorphism $\be_m:\X(\g_n)\to\B_m\ts$. Recall that
in Section 2, the target algebra $\B_m$ of the homomorphism $\be_m$
was defined as 
$\U(\ts\f_m)\ot\ts\PD\ts(\CC^{\ts m}\ot\CC^{\ts n})\ts$.
Here we use a different presentation of
the same algebra, by means of the cross relations \eqref{defar}.
In particular, here the image of $\be_m$
commutes with the subalgebra $\U(\ts\f_m)$ of $\B_m\ts$; see
Proposition \ref{xb}, Part (ii).
Thus here the vector space
$\B_m\ts/\,\I_{\ts\mu,\de}$ becomes a bimodule
over $\f_m$ and $\X(\g_n)\ts$.

Consider the bimodule $\F_\de\ts(M_\mu)$
over $\f_m$ and $\X(\g_n)\ts$,
defined in the beginning of this section.
This bimodule is equivalent to
$\B_m\ts/\,\I_{\ts\mu,\de}\ts$. Indeed, let $Z$ run through
$\P\ts(\CC^{\ts m}\ot\CC^{\ts n})\ts$. Then a bijective linear map
$$
\F_\de\ts(M_\mu)\to\B_m\ts/\,\I_{\ts\mu,\de}
$$
intertwining the actions of $\f_m$ and $\X(\g_n)$ can be defined
by mapping the element
$$
1_\mu\ot Z\in M_\mu\ot\P\ts(\CC^{\ts m}\ot\CC^{\ts n})
$$
to the image of 
$$
\varpi^{-1}(Z)\in\PD\ts(\CC^{\ts m}\ot\CC^{\ts n})\subset\B_m
$$
in the quotient $\B_m\ts/\,\I_{\ts\mu,\de}\ts$.
The intertwining property here follows from the definitions
of $\F_\de\ts(M_\mu)$ and $\I_{\ts\mu,\de}\ts$.
The same mapping determines a bijective linear map
\begin{equation}
\label{fbi}
\F_\de\ts(M_\mu)\ts\to\ts\Bb_m\ts/\,\Ib_{\ts\mu,\de}\ts.
\end{equation}

In particular, the space $\F_\de\ts(M_\mu)_{\ts\n}$
of $\n\ts$-coinvariants of $\F_\de\ts(M_\mu)$
is equivalent to the quotient
$\Jb\,\ts\backslash\,\Bb_m\ts/\,\Ib_{\ts\mu,\de}$
as a bimodule over the Cartan subalgebra
$\h\subset\f_m$ and over $\X(\g_n)\ts$.
But Theorem \ref{proposition4.5} implies
that the operator $\ts\xic_{\ts\si}$ on
$\Jb\,\ts\backslash\,\Bb_m$ determines a linear map
\begin{equation}
\label{bbjioper}
\Jb\,\ts\backslash\,\Bb_m\ts/\,\Ib_{\ts\mu,\de_+}
\to\,
\Jb\,\ts\backslash\,\Bb_m\ts/\,\Ib_{\ts\si\ts\circ\ts\mu\ts,\ts\si(\de_+)}
\ts.
\end{equation}
The latter map intertwines the actions
of $\X(\g_n)$ on the source and the target vector spaces,
because the image of $\X(\g_n)$ in
$\B_m$ relative to $\be_m$ commutes with the subalgebra
$\U(\ts\f_m)\subset\B_m\ts$; see the definition \eqref{q1}.
We also use Lemma \ref{lemma41}, Part (ii).
Recall that $\F_m(V)=\F_{\ts\de_+}(V)\ts$.
Hence by using the equivalences \eqref{fbi}
for the sequences $\de=\de_+$ and
$\de=\si\ts(\de_+)\ts$, the operator \eqref{bbjioper} becomes
the desired $\X(\g_n)\ts$-intertwining operator \eqref{distoper}.

As usual, for any $\f_m$-module $V$ and any element $\la\in\h^\ast$
let $V^\la\subset V$ be the subspace of vectors \textit{of weight\/}
$\la\ts$, so that any $X\in\h$ acts on $V^\la$
via multiplication by $\la\ts(X)\in\CC\ts$.
It now follows from the property \eqref{q121} of $\ts\xic_{\ts\si}$
that the restriction of
our operator \eqref{distoper} to the subspace of weight $\la$
is an $\X(\g_n)\ts$-intertwining operator
\begin{equation}
\label{distoperla}
\F_{m}\ts(\ts M_{\ts\mu})_{\ts\n}^{\ts\la}
\,\to\,
\F_{\ts\de}\ts(\ts M_{\ts\si\,\circ\ts\mu})_{\ts\n}^{\,\si\,\circ\ts\la}
\ \quad\text{where}\quad\
\de=\si\ts(\de_+)\ts.
\end{equation}

At the end of Section 2, we defined the modules $P_z$ and
$\Pp_z$ over the Yangian $\Y(\gl_n)\ts$. The underlying vector space
of these modules is the Grassmann algebra $\P\,(\CC^{\ts n})\ts$.
This algebra is graded by $0\com1\lcd n\ts$.
The actions of $\Y(\gl_n)$ on $P_z$ and $\Pp_z$
preserve the degree. Now for any
$N=1\lcd n$ denote respectively by $P_z^{\ts N}$ and
$P_z^{\ts-N}$ the submodules in $P_z$ and $\Pp_z$ which
consist of the elements of degree $N\ts$.
Note that $\Y(\gl_n)$ acts on the subspace of $P_z$
of degree zero trivially, that is via the
counit homomorphism $\Y(\gl_n)\to\CC\ts$.
That action of $\Y(\gl_n)$ does not depend on $z\ts$.
It will be convenient to denote by $P_z^{\ts0}$ the
vector space $\CC$ with the trivial action of $\Y(\gl_n)\ts$.

Denote
\begin{equation}
\label{nua}
\nu_a={n}/2+\mu_a-\la_a
\quad\textrm{for}\quad
a=1\lcd m\ts.
\end{equation}
Suppose that $\nu_1\lcd\nu_m\in\{\tts0\com1\lcd n\}\ts$,
otherwise the source $\X(\g_n)\ts$-module in
\eqref{distoperla} would be zero by Corollary \ref{verma}.
Under our assumption, Corollary \ref{verma} implies that
the the source $\X(\g_n)\ts$-module in \eqref{distoperla} is equivalent to
\begin{equation}
\label{munup}
P_{\mu_m+z}^{\,\nu_m}
\ot
P_{\mu_{m-1}+z+1}^{\,\nu_{m-1}}
\ot\ldots\ot
P_{\mu_1+z+m-1}^{\,\nu_1}
\end{equation}
pulled back through the automorphism \eqref{fus} of $\X(\g_n)\ts$,
where $f(u)$ is given by \eqref{fuprod} and $z=\mp\tts\frac12\ts$.
A more general results is stated as Proposition \ref{siverma} below.
The tensor product in \eqref{munup} is that of
$\Y(\gl_n)\ts$-modules. Then we employ the embedding
$\Y(\g_n)\subset\Y(\gl_n)$ and the homomorphism
$\X(\g_n)\to\Y(\g_n)$ defined by \eqref{xy}.
By using the labels $\rho_1\lcd\rho_m$
of the halfsum $\rho$ of the positive roots of $\f_m\ts$,
the tensor product \eqref{munup} can be rewritten as
\begin{equation}
\label{pp}
P_{\mu_m-\frac12+\rho_m}^{\,\nu_m}
\ot\ldots\ot
P_{\mu_1-\frac12+\rho_1}^{\,\nu_1}\ts.
\end{equation}
By using the labels $\rho_1\lcd\rho_m$
we can also rewrite the product \eqref{fuprod} as
\begin{equation}
\label{muprod}
\prod_{a=1}^m\,\ts
\frac{u-\mu_a+\frac12-\rho_a}{u-\mu_a-\frac12-\rho_a}\ .
\end{equation}

Let us now consider the target $\X(\g_n)\ts$-module in \eqref{distoperla}.
For each $a=1\lcd m$ denote
$$
\widetilde\mu_a=\mu_{\ts|\si^{-1}(a)|}\ts,
\quad
\widetilde\nu_a=\nu_{\ts|\si^{-1}(a)|}\ts,
\quad
\widetilde\rho_a=\rho_{\ts|\si^{-1}(a)|}\ts.
$$
The above description of the source $\X(\g_n)\ts$-module in
\eqref{distoperla} can now be generalized to the target $\X(\g_n)\ts$-module,
which depends on an arbitrary element $\si\in\H_m\ts$.

\begin{proposition*}
\label{siverma}
For $\de=\si\ts(\de_+)$ the $\X(\g_n)\ts$-module
$\F_{\ts\de}\ts(\ts M_{\ts\si\,\circ\ts\mu})_{\ts\n}^{\,\si\,\circ\ts\la}$
is equivalent to the tensor product
\begin{equation}
\label{ppd}
P_{\widetilde\mu_m-\frac12+\widetilde\rho_m}^{\,\ts\de_m\ts\widetilde\nu_m}
\ot\ldots\ot
P_{\widetilde\mu_1-\frac12+\widetilde\rho_1}^{\,\ts\de_1\ts\widetilde\nu_1}
\end{equation}
pulled back through the automorphism \eqref{fus} of\/ $\X(\g_n)$
where $f(u)$ equals\/ \eqref{muprod}.
\end{proposition*}

\begin{proof}
First consider the bimodule
$\F_m(\ts  M_{\ts\si\,\circ\ts\mu})_{\ts\n}$ of $\h$ and $\X(\g_n)\ts$.
By Corollary~\ref{verma}, this bimodule is equivalent to the tensor product
\begin{equation}
\label{bim1}
P_{\ts\de_m\ts\widetilde\mu_m-\frac12+\ts\de_m\ts\widetilde\rho_m}
\ot\ldots\ot
P_{\ts\de_1\ts\widetilde\mu_1-\frac12+\ts\de_1\ts\widetilde\rho_1}
\end{equation}
pulled back through the automorphism \eqref{fus} of $\X(\g_n)$ where
$f(u)$ equals
\begin{equation}
\label{deprod}
\prod_{a=1}^m\,\ts\frac
{u-\ts\de_a\ts\widetilde\mu_a+\frac12-\ts\de_a\ts\widetilde\rho_a}
{u-\ts\de_a\ts\widetilde\mu_a-\frac12-\ts\de_a\ts\widetilde\rho_a}\ .
\end{equation}
For any $a=1\lcd m$ the element $F_{-\ab,-\ab}\in\h$ acts on
the tensor product \eqref{bim1} as
$$
{n}/2\tts-\tts\deg{\nns}_a+(\si\circ\mu)_a
$$
where $\deg{\nns}_a$ is the degree operator
on the $a\ts$-th tensor factor, counting the factors from right
to left. It acts on the vector space $\P\ts(\CC^{\ts n})$ of that
tensor factor as the Euler operator
\begin{equation}
\label{euler}
\sum_{k=1}^n\,x_k\ts\d_k\in\PD\ts(\CC^{\ts n})\,.
\end{equation}

A bimodule equivalent to
$\F_\de\ts(\ts  M_{\ts\si\,\circ\ts\mu})_{\ts\n}$ can be obtained
by pushing forward actions of $\h$ and $\X(\g_n)$ on \eqref{bim1}
through the composition of automorphisms \eqref{onefour},
for every tensor factor with number $a\ts$ such that $\de_a=-1$.
Here we number the $m$ tensor factors of \eqref{bim1} by $1\lcd m$
from right to left. Then we also have to pull the resulting
$\X(\g_n)\ts$-module back through the automorphism
\eqref{fus}, where the series $f(u)$ equals \eqref{deprod}.
The automorphism \eqref{onefour} maps the element
\eqref{euler}~to
$$
\sum_{k=1}^n\,\d_{\ts\bk}\ts x_{\ts\bk}\,=\,
n-\,\sum_{k=1}^n\,x_k\ts\d_k\ts.
$$
Hence if $\de_a=-1\ts$, the element
$F_{-\ab,-\ab}\in\h$ acts on the modified tensor product as
$$
-\,{n}/2+(\si\circ\mu)_a+\deg{\nns}_a\ts.
$$
By equating the last displayed expression to $(\si\circ\la)_a$
and by using \eqref{nua} together with the condition $\de_a=-1\ts$,
we get the equation $\deg{\nns}_a=\widetilde\nu_a\ts$.
But by Lemma~\ref{ppl}, pushing forward the $\Y(\gl_n)\ts$-module
$$
P_{-\widetilde\mu_a-\frac12-\widetilde\rho_a}^{\,\ts\widetilde\nu_a}
$$
through the automorphism \eqref{onefour} of $\PD\ts(\CC^{\ts n})$
yields the same $\Y(\gl_n)\ts$-module as pulling
$$
P_{\widetilde\mu_a-\frac12+\widetilde\rho_a}^{\ts-\ts\widetilde\nu_a}
$$
back through the automorphism \eqref{fut} of $\Y(\gl_n)$ where
$$
g(u)=
\frac{u-\widetilde\mu_a+\frac12-\widetilde\rho_a}
{u-\widetilde\mu_a-\frac12-\widetilde\rho_a}\ .
$$
Thus the $\X(\g_n)\ts$-module
$\F_{\ts\de}\ts(\ts M_{\ts\si\,\circ\ts\mu})_{\ts\n}^{\,\si\,\circ\ts\la}$
is equivalent to the tensor product \eqref{ppd}
pulled back through the automorphism \eqref{fus}
where the series $f(u)$ is obtained by multiplying \eqref{deprod}
by $g(-u)\ts g(u)$ for each index $a\ts$ such that $\de_a=-1\ts$;
see the definition \eqref{xy}. But for any 
the element $\si\in\H_m$ the product \eqref{muprod} equals
\begin{equation}
\label{tildaprod}
\prod_{a=1}^m\,\ts
\frac{u-\widetilde\mu_a+\frac12-\widetilde\rho_a}
{u-\widetilde\mu_a-\frac12-\widetilde\rho_a}\ .
\end{equation}
If $\de_a=-1$ then the factors of \eqref{deprod} and \eqref{tildaprod}
indexed by $a$ are equal to $g(-u)^{-1}$ and $g(u)$ respectively.
If $\de_a=1$ then the factors of \eqref{deprod} and \eqref{tildaprod}
indexed by $a$ coincide.
This comparison of \eqref{deprod} and \eqref{tildaprod}
completes the proof. 
\qed
\end{proof}

The vector spaces of two equivalent
$\X(\g_n)\ts$-modules in Proposition \ref{siverma} are
$$
(M_{\ts\si\,\circ\ts\mu}\ot
\P\ts(\CC^{\ts m}\ot\CC^{\ts n}))_{\ts\n}^{\,\si\,\circ\ts\la}
\ \quad\text{and}\ \quad
\P^{\,\widetilde\nu_m}\tts(\CC^{\ts n})
\ot\ldots\ot
\P^{\,\widetilde\nu_1}\tts(\CC^{\ts n})
$$
respectively.
We can define a linear map from the latter vector space
to the former,
by mapping $f_1\ot\ldots\ot f_m$ to
the class of 
$1_{\ts\si\,\circ\ts\mu}\ot f$ in
the space of $\n\ts$-coinvariants. Here
$$
f_1\in\P^{\,\widetilde\nu_m}\tts(\CC^{\ts n})
\,\lcd\,
f_m\in\P^{\,\widetilde\nu_1}\tts(\CC^{\ts n})
$$
and
$f\in\P\ts(\CC^{\ts m}\ot\CC^{\ts n})$ is defined by \eqref{fff}.
This linear map is an equivalence of the $\X(\g_n)\ts$-modules
in Proposition \ref{siverma},
see the remarks made after our proof of Corollary~\ref{verma}.

Thus for any $\nu_1\lcd\nu_m\in\{\tts0\com1\lcd n\}$
we have demonstrated how the Zhelobenko operator $\ts\xic_{\ts\si}$
on $\Jb\,\ts\backslash\,\Bb_m$ determines an intertwining operator
between the $\X(\g_n)\ts$-modules \eqref{pp} and \eqref{ppd}
pulled back via the automorphism \eqref{fus} of $\X(\g_n)\ts$,
where 
$f(u)$ is the same \eqref{muprod} for both modules.
Hence this operator also intertwines the $\X(\g_n)\ts$-modules
\begin{equation}
\label{ppi}
P_{\mu_m-\frac12+\rho_m}^{\,\nu_m}
\ot\ldots\ot
P_{\mu_1-\frac12+\rho_1}^{\,\nu_1}
\to\,
P_{\widetilde\mu_m-\frac12+\widetilde\rho_m}^{\,\ts\de_m\ts\widetilde\nu_m}
\ot\ldots\ot
P_{\widetilde\mu_1-\frac12+\widetilde\rho_1}^{\,\ts\de_1\ts\widetilde\nu_1}
\ts,
\end{equation}
neither of them being pulled back via the authomorphism \eqref{fus}.
It was proved in \cite{MN} that
both $\X(\g_n)\ts$-modules in \eqref{ppi} are irreducible
under our assumptions on $\mu\ts$.
Hence an intertwining operator between them
is unique up to a multiplier from $\CC\ts$.
For our intertwining operator,
this multiplier is determined by Proposition~\ref{isis} below.
Another expression for an intertwining operator
of the $\X(\g_n)\ts$-modules \eqref{ppi}
was given in \cite{N2}.

For any $a=1\lcd m$ and $s=1\lcd n$ let us define the elements
$f_{as}$ and $g_{\tts as}$ of the ring
$\PD\ts(\CC^{\ts m}\ot\CC^{\ts n})$ as follows.
Let us arrange the indices $1\lcd n$ into the sequence
\begin{equation}
\label{corder}
1\com3\lcd n-1\com n\lcd4\com 2
\ \quad\text{or}\ \quad
1\com3\lcd n-2\com n\com n-1\lcd4\com 2
\end{equation}
when 
$n$ is even or odd respectively.
The mapping $k\mapsto\bk$ reverses the sequence \eqref{corder}.
We will write $i\prec j$ when $i$ precedes $j\/$
in this sequence. Note that then the elements
$
E_{\ts ij}-\th_i\ts\th_j\tts E_{\ts\bj\ts\bi}\in\gl_n
$
with $i\prec j$ or $i=j$ span a Borel subalgebra of $\g_n\subset\gl_n\ts$,
while the elements $E_{\ts ii}-E_{\ts\bi\ts\bi}$
span the corresponding Cartan subalgebra of $\g_n\ts$.
Then $f_{as}$ and $g_{\tts as}$ are
the products of the elements $x_{ak}$ and $\d_{a\tts\bk}$
of $\ts\PD\ts(\CC^{\ts m}\ot\CC^{\ts n})$ respectively,
taken over the first $s$ indices $k$ in the sequence \eqref{corder}.
For example, if $n\ge4$ then
$f_{a\tts2}=x_{a1}x_{a3}$ and
$g_{\tts a\tts2}=\d_{a\tts2}\ts\d_{a\tts4}\ts$.
If $n=3$ then
$f_{a\tts2}=x_{a1}x_{a3}$ but
$g_{\tts a\tts2}=\d_{a\tts2}\ts\d_{a\tts3}\ts$.
We also set $f_{a\tts0}=g_{\tts a\tts0}=1\ts$.

Our proof of Proposition \ref{isis} will be based on the
following four lemmas. The proof of the first lemma
is very similar to the proof of the second one, and will be omitted.

\begin{lemma*}
\label{ddnorma}
For any\/ $a=1\lcd m-1$ and\/ $s\com t=0\com1\lcd n$
the operator\/ $\xic_{\ts a}$ on\/ $\Jb\,\ts\backslash\,\Bb_m$
maps the image in\/ $\Jb\,\ts\backslash\,\Bb_m$ of\/
$g_{\ts\overline{a}\ts s}\,g_{\ts\overline{a+1}\ts t}\in\Bb_m$
to the image in\/ $\Jb\,\ts\backslash\,\Bb_m$ of the product
$$
\sih_a\ts(\ts g_{\ts\overline{a}\ts s}\,g_{\ts\overline{a+1}\ts t}\ts)
\ts\cdot
\left\{\begin{array}{cc}
\dfrac{H_a-s+t+1}{H_a+1}
&\quad\text{if}\quad s<t\ts,
\\[10pt]
1
&\quad\text{if}\quad s\ge t\ts,
\end{array}\right.
$$
plus the images in\/ $\Jb\,\ts\backslash\,\Bb_m$
of elements of the left ideal in\/ $\Bb_m$ generated by\/ $\Jb'$ and
\eqref{xxideal}.
\end{lemma*}

\begin{lemma*}
\label{xxnorma}
For any\/ $a=1\lcd m-1$ and\/ $s\com t=0\com1\lcd n$
the operator\/ $\xic_{\ts a}$ on\/ $\Jb\,\ts\backslash\,\Bb_m$
maps the image in\/ $\Jb\,\ts\backslash\,\Bb_m$ of\/
$f_{\ts\overline{a}\ts s}\,f_{\ts\overline{a+1}\ts t}\in\Bb_m$
to the image in\/ $\Jb\,\ts\backslash\,\Bb_m$ of the product
$$
\sih_a\ts(\ts f_{\ts\overline{a}\ts s}\,f_{\ts\overline{a+1}\ts t}\ts)
\ts\cdot
\left\{\begin{array}{cc}
\dfrac{H_a+s-t+1}{H_a+1}
&\quad\text{if}\quad s>t\ts,
\\[10pt]
1
&\quad\text{if}\quad s\le t\ts,
\end{array}\right.
$$
plus the images in\/ $\Jb\,\ts\backslash\,\Bb_m$
of elements of the left ideal in\/ $\Bb_m$ generated by\/ $\Jb'$ and
\eqref{ddideal}.
\end{lemma*}

\begin{proof}
By the definitions \eqref{gan} and \eqref{Fc}, we have
\begin{equation}
\label{efahat}
\zeta_n\ts(E_a)=-\,\sum_{k=1}^n\,
x_{\ts\overline{a+1}\ts k}\,\d_{\ts\overline{a}\ts k}
\ \quad\text{and}\ \quad
\zeta_n\ts(F_a)=-\,\sum_{k=1}^n\,
x_{\ts\overline{a}\ts k}\,\d_{\ts\overline{a+1}\ts k}
\,.
\end{equation}
By \eqref{fourier}, we also have
$$
\sih_a\ts(\ts f_{\ts\overline{a}\ts s}\,f_{\ts\overline{a+1}\ts t}\ts)
\,=\,
f_{\ts\overline{a+1}\ts s}\,f_{\ts\overline{a}\ts t}\,.
$$
Let us now use the symbol $\,\equiv\,$ to indicate equalities
in the vector space $\Jb\,\ts\backslash\,\Bb_m$ modulo the subspace,
which is the image of the left ideal in $\Bb_m$ generated by $\Jb'$ and
by the elements \eqref{ddideal}. The element $E_a\in\Bb_m$
belongs to this left ideal. Therefore the operator $\xic_{\ts a}$ maps
the image in $\Jb\,\ts\backslash\,\Bb_m$ of
$f_{\ts\overline{a}\ts s}\,f_{\ts\overline{a+1}\ts t}\in\Bb_m$
to the image in $\Jb\,\ts\backslash\,\Bb_m$~of
\begin{align*}
\xi_{\ts a}\ts(\ts
f_{\ts\overline{a+1}\ts s}\,
f_{\ts\overline{a}\ts t}\ts)
&\,\,=\,\,
\sum_{r=0}^\infty\,\,
(\ts r\ts!\,H_a^{\ts(r)}\ts)^{-1}\ts E_a^{\ts r}\,
\widehat{F}_a^{\ts r}(\ts
f_{\ts\overline{a+1}\ts s}\,
f_{\ts\overline{a}\ts t}\ts)
\\
&\,\,\equiv\,\,
\sum_{r=0}^\infty\,\,
(\ts r\ts!\,H_a^{\ts(r)}\ts)^{-1}\ts
\widehat{E}_a^{\ts r}\,
\widehat{F}_a^{\ts r}(\ts
f_{\ts\overline{a+1}\ts s}\,
f_{\ts\overline{a}\ts t}\ts)\ts.
\end{align*}

Let us now use \eqref{defar} along with \eqref{efahat}.
By the definitions of
$f_{\ts\overline{a+1}\ts s}$ and $f_{\ts\overline{a}\ts t}$ we have
$$
\widehat{F}_a\ts(\ts
f_{\ts\overline{a+1}\ts s}\,
f_{\ts\overline{a}\ts t}\ts)
\,=\,
-\,\sum_{k=1}^n\,\,
[\ts x_{\ts\ab\ts k}\,\d_{\,\ts\overline{a+1}\ts k}\,,
f_{\ts\overline{a+1}\ts s}\,
f_{\ts\overline{a}\ts t}\ts]\,.
$$
If $s\le t\ts$, then every summand above is zero,
which proves the lemma in this case.
Now suppose that $s>t\ts$. Then by using the proof of
\cite[Proposition 3.7]{KN2},
$$
\xi_{\ts a}\ts(\ts
f_{\ts\overline{a+1}\ts s}\,
f_{\ts\overline{a}\ts t}\ts)
\,\,\equiv\,\,
\sum_{r=0}^{s-t}\,
\frac{\ts (s-t)\ldots(s-t-r+1)}{H_a\ldots(H_a-r+1)}\,\,
f_{\ts\overline{a+1}\ts s}\,
f_{\ts\overline{a}\ts t}\,.
$$
In the last line,
the sum of the fractions corresponding to $r=0\lcd s-t$ equals
$$
\frac{H_a+1}{H_a-s+t+1}\,;
$$
this equality can be easily proved by induction on the difference
$s-t\ts$.
Therefore
$$
\xi_{\ts a}\ts(\ts
f_{\ts\overline{a+1}\ts s}\,
f_{\ts\overline{a}\ts t}\ts)
\,\,\equiv\,\,
\frac{H_a+1}{H_a-s+t+1}
\,\,
f_{\ts\overline{a+1}\ts s}\,
f_{\ts\overline{a}\ts t}
\,\,=\,\,
f_{\ts\overline{a+1}\ts s}\,
f_{\ts\overline{a}\ts t}
\,\,
\frac{H_a+s-t+1}{\ts H_a+1}
$$
as required in the case when $s>t\ts$.
Here we also used the equality in the ring $\B_m\ts$,
$$
H_a\,
f_{\ts\overline{a+1}\ts s}\,
f_{\ts\overline{a}\ts t}
\,=\,
f_{\ts\overline{a+1}\ts s}\,
f_{\ts\overline{a}\ts t}\,
(H_a+s-t)
$$
which follows from \eqref{defar}, since
$$
\zeta_n(H_a)
\,=\,
\zeta_n\ts(\ts
F_{\ts\overline{a+1}\ts,\ts\overline{a+1}}\,-F_{\ts\ab\ts\ab}\ts)
\,=\,
\sum_{k=1}^n\,
(\ts
x_{\ts\overline{a+1}\ts k}\,
\d_{\ts\overline{a+1}\ts k}\,
-
x_{\ts\overline{a}\ts k}\,
\d_{\ts\overline{a}\ts k}\ts
)\ts.\qquad
\eqno\qed
$$
\end{proof}

\begin{lemma*}
\label{xdnorma}
For any\/ $a=1\lcd m-1$ and\/ $s\com t=0\com1\lcd n$
the operator\/ $\xic_{\ts a}$ on\/ $\Jb\,\ts\backslash\,\Bb_m$
maps the image in\/ $\Jb\,\ts\backslash\,\Bb_m$ of\/
$f_{\ts\overline{a}\ts s}\,g_{\ts\overline{a+1}\ts t}\in\Bb_m$
to the image in\/ $\Jb\,\ts\backslash\,\Bb_m$ of the product
$$
\sih_a\ts(\ts f_{\ts\overline{a}\ts s}\,g_{\ts\overline{a+1}\ts t}\ts)
\ts\cdot
\left\{\begin{array}{cc}
\dfrac{H_a+s+t+1}{H_a+n+1}
&\quad\text{if}\quad s+t>n\ts,
\\[10pt]
1
&\quad\text{if}\quad s+t\le n\ts,
\end{array}\right.
$$
plus the images in\/ $\Jb\,\ts\backslash\,\Bb_m$
of elements of the left ideal in\/ $\Bb_m$ generated by\/ $\Jb'$ and
\eqref{xdideal}.
\end{lemma*}

\begin{proof}
By \eqref{fourier},
$$
\sih_a\ts(\ts f_{\ts\overline{a}\ts s}\,g_{\ts\overline{a+1}\ts t}\ts)
\,=\,
f_{\ts\overline{a+1}\ts s}\,g_{\ts\overline{a}\ts t}\,.
$$
Let now us the symbol $\,\equiv\,$ to indicate equalities
in $\Jb\,\ts\backslash\,\Bb_m$ modulo the subspace,
which is the image of the left ideal in $\Bb_m$ generated by $\Jb'$ and
by the elements \eqref{xdideal}.
The elements $E_a-\zeta_n(E_a)$ and $\zeta_n(F_a)$ of $\Bb_m$
belong to this left ideal, see \eqref{efahat}.
Using \eqref{defar}, the operator $\xic_{\ts a}$ maps
the image in $\Jb\,\ts\backslash\,\Bb_m$ of
$f_{\ts\overline{a}\ts s}\,g_{\ts\overline{a+1}\ts t}\in\Bb_m$
to the image in $\Jb\,\ts\backslash\,\Bb_m$~of
\begin{align*}
\xi_{\ts a}\ts(\ts
f_{\ts\overline{a+1}\ts s}\,
g_{\ts\overline{a}\ts t}\ts)
&\,\,=\,\,
\sum_{r=0}^\infty\,\,
(\ts r\ts!\,H_a^{\ts(r)}\ts)^{-1}\ts E_a^{\ts r}\,
\widehat{F}_a^{\ts r}\ts(\ts
f_{\ts\overline{a+1}\ts s}\,
g_{\ts\overline{a}\ts t}\ts)
\\
&\,\,\equiv\,\,
\sum_{r=0}^\infty\,\,
(\ts r\ts!\,H_a^{\ts(r)}\ts)^{-1}\ts
\zeta_n(E_a)^{\ts r}\ts
\zeta_n(F_a)^{\ts r}\ts
f_{\ts\overline{a+1}\ts s}\,
g_{\ts\overline{a}\ts t}\,.
\end{align*}

We have
$$
\zeta_n(F_a)\ts
f_{\ts\overline{a+1}\ts s}\,
g_{\ts\overline{a}\ts t}\,=\,
-\,\sum_{k=1}^n\,
x_{\ts\overline{a}\ts k}\,
\d_{\ts\overline{a+1}\ts k}\,
f_{\ts\overline{a+1}\ts s}\,
g_{\ts\overline{a}\ts t}
$$
by \eqref{efahat}.
If $s+t\le n\ts$, then
every summand in the above displayed
sum is zero modulo the left ideal of $\Bb_m$ generated by
the elements \eqref{xdideal}, because then
there is no factors
$x_{\ts\overline{a+1}\ts i}$ of $f_{\ts\overline{a+1}\ts s}$
and $\d_{\ts\overline{a}\ts i}$ of $g_{\ts\overline{a}\ts t}$
with the same index $i\ts$. This proves the lemma in this case.
Now suppose that $s+t>n\ts$. Then, by using the proof of
\cite[Proposition 3.7]{KN2},
\begin{gather*}
\xi_{\ts a}\ts(\ts
f_{\ts\overline{a+1}\ts s}\,
g_{\ts\overline{a}\ts t}\ts)
\,\,\equiv
\sum_{r=0}^{s+t-n}\,
\frac{\ts (s+t-n)\ldots(s+t-n-r+1)}{H_a\ldots(H_a-r+1)}\,\,
f_{\ts\overline{a+1}\ts s}\,
g_{\ts\overline{a}\ts t}\,\,=
\\[8pt]
\frac{H_a+1}{H_a-s-t+n+1}
\,\,
f_{\ts\overline{a+1}\ts s}\,
g_{\ts\overline{a}\ts t}
\,\,=\,\,
f_{\ts\overline{a+1}\ts s}\,
g_{\ts\overline{a}\ts t}
\,\,
\frac{H_a+s+t+1}{\ts H_a+n+1}
\end{gather*}
as required. Here we have also used an equality in the ring $\B_m$
which follows from \eqref{defar},
$$
H_a\,
f_{\ts\overline{a+1}\ts s}\,
g_{\ts\overline{a}\ts t}
\,=\,
f_{\ts\overline{a+1}\ts s}\,
f_{\ts\overline{a}\ts t}\,
(H_a+s+t)\ts.
\eqno\qed
$$
\end{proof}

\begin{lemma*}
\label{xnorma}
If\/ $\f_m=\sp_{2m}\ts$, then
for any\/ $s=0\com1\lcd n$
the operator\/ $\xic_{\ts m}$ on\/ $\Jb\,\ts\backslash\,\Bb_m$
maps the image in\/ $\Jb\,\ts\backslash\,\Bb_m$ of\/
$f_{\ts\overline{m}\ts s}\in\Bb_m$
to the image in\/ $\Jb\,\ts\backslash\,\Bb_m$ of the product
$$
\sih_m\ts(\ts f_{\ts\overline{m}\ts s})
\ts\cdot
\left\{\begin{array}{cc}
\dfrac{H_a+s+1}{H_a+n/2+1}
&\quad\text{if}\quad s>n/2\ts,
\\[10pt]
1
&\quad\text{if}\quad s\le n/2\ts,
\end{array}\right.
$$
plus the images in\/ $\Jb\,\ts\backslash\,\Bb_m$
of elements of the left ideal in\/ $\Bb_m$ generated by\/ $\Jb'$ and
\eqref{xideal}.
\end{lemma*}

\begin{proof}
Let $\f_m=\sp_{2m}\ts$. Then $\g_n=\sp_n\ts$,
so that the number $n$ is even.
By \eqref{fourier}, we have
$$
\sih_m\ts(\ts f_{\ts\overline{m}\ts s})=
g_{\ts\overline{m}\ts s}
\ \quad\text{or}\ \quad
\sih_m\ts(\ts f_{\ts\overline{m}\ts s})=
(-1)^{\ts s-n/2}\,g_{\ts\overline{m}\ts s}
$$
when $s\le n/2$ or $s>n/2$ respectively.
Hence it suffices to consider for any $s=0\com1\lcd n$
the image in $\Jb\,\ts\backslash\,\Bb_m$ of the element
$\xi_{\ts m}\ts(\ts g_{\ts\overline{m}\ts s})\in\Bb_m\ts$.
By the definitions \eqref{gan} and \eqref{Fm},
\begin{equation*}
\label{efmhat}
\zeta_n\ts(E_m)\,=\,\sum_{k=1}^n\,
\theta_k\,
\d_{\ts\overline{m}\ts k}\,
\d_{\ts\overline{m}\ts\bk}\,/\ts2
\ \quad\text{and}\ \quad
\zeta_n\ts(F_m)\,=\,\sum_{k=1}^n\,
\theta_k\,
x_{\ts\overline{m}\ts\bk}\,
x_{\ts\overline{m}\ts k}\,/\ts2\,.
\end{equation*}
Now let the symbol $\,\equiv\,$ indicate equalities
in $\Jb\,\ts\backslash\,\Bb_m$ modulo the subspace,
which is the image of the left ideal in $\Bb_m$ generated by $\Jb'$ and
by the elements \eqref{xideal}.
The elements $E_m-\zeta_n(E_m)$ and $\zeta_n(F_m)$ of $\Bb_m$
belong to this left ideal. Therefore by using \eqref{defar},
\begin{align*}
\xi_{\ts m}\ts(\ts
g_{\ts\overline{m}\ts s})
&\,\,=\,\,
\sum_{r=0}^\infty\,\,
(\ts r\ts!\,H_m^{\ts(r)}\ts)^{-1}\ts E_m^{\ts r}\,
\widehat{F}_m^{\ts r}\ts(\ts
g_{\ts\overline{m}\ts s})
\\
&\,\,\equiv\,\,
\sum_{r=0}^\infty\,\,
(\ts r\ts!\,H_m^{\ts(r)}\ts)^{-1}\ts
\zeta_n(E_m)^{\ts r}\ts
\zeta_n(F_m)^{\ts r}\ts
g_{\ts\overline{m}\ts s}\,.
\end{align*}

We have
$$
\zeta_n(F_m)\,
g_{\ts\overline{m}\ts s}
\,=\,\sum_{k=1}^n\,
\theta_k\,
x_{\ts\overline{m}\ts\bk}\,
x_{\ts\overline{m}\ts k}\,
g_{\ts\overline{m}\ts s}\,/\ts2\,.
$$
If $s\le n/2\ts$, then
every summand in the above
sum is zero modulo the left ideal of $\Bb_m$ generated by
the elements \eqref{xideal}, because then for any index $k\ts$
there is no pair of factors
$\d_{\ts\overline{m}\ts k}$ and $\d_{\ts\overline{m}\ts\bk}$
in the product $g_{\ts\overline{m}\ts s}\ts$.
This proves the lemma in this case.
Now suppose that $s>n/2\ts$. Then, by using the proof of
\cite[Proposition 3.7]{KN2} once again, we have
\begin{gather*}
\xi_{\ts m}\ts(\ts g_{\ts\overline{m}\ts s})
\,\,\equiv
\sum_{r=0}^{s-n/2}\,
\frac{\ts (s-n/2)\ldots(s-n/2-r+1)}{H_m\ldots(H_m-r+1)}\,\,
g_{\ts\overline{m}\ts s}\,\,=
\\[8pt]
\frac{H_m+1}{H_m-s+n/2+1}
\,\,
g_{\ts\overline{m}\ts s}
\,\,=\,\,
g_{\ts\overline{m}\ts s}
\,\,
\frac{H_m+s+1}{\ts H_a+n/2+1}
\end{gather*}
as required. Here we also used the equality
$
H_m\,g_{\ts\overline{m}\ts s}=
g_{\ts\overline{m}\ts s}\,(H_m+s)
$
in the ring $\B_m\ts$, which follows from \eqref{defar},
because $\overline{m}=1$ and for $\f_m=\sp_{2m}$ by
\eqref{gan} and \eqref{Fm} we have
$$
\zeta_n(H_m)
\,=\,
-\,\zeta_n\ts(F_{11})
\,=\,
n/2\,-\,
\sum_{k=1}^n\,
x_{1k}\,\d_{1k}\,.
\eqno\qed
$$
\end{proof}

Let us now state Proposition \ref{isis}.
We assume that the weight $\mu$ satisfies the
conditions \eqref{notinso}, and also satisfies the conditions \eqref{notinsp}
in the case $\f_m=\sp_{2m}\ts$. We also assume that
$\nu_1\lcd\nu_m\in\{\tts0\com1\lcd n\}\ts$, see the definition \eqref{nua}.
Let $(\ts\mu_1^\ast\lcd\mu_m^\ast\ts)$ be the sequence of labels of the
weight $\mu+\rho\ts$. Thus for each $a=1\lcd m$ we have
$\mu_a^\ast=\mu_a+m-a$ in the case $\f_m=\so_{2m}\ts$, and
$\mu_a^\ast=\mu_a+m-a+1$ in the case $\f_m=\sp_{2m}\ts$.
Let $(\ts\la_1^\ast\lcd\la_m^\ast\ts)$ be the sequence of labels
of $\la+\rho\ts$.
For each positive root $\eta\in\De^+$ define a number $z_\eta\in\CC\ts$,
$$
z_\eta\,=\,
\left\{
\begin{array}{cll}
\,\dfrac{\la_b^\ast-\la_c^\ast}{\mu_b^\ast-\mu_c^\ast}
&\quad\text{if}\quad
\eta=\ep_{b}-\ep_c
&\quad\text{and}\quad\nu_b>\nu_c\,,
\\[12pt]
\,\dfrac{\la_b^\ast+\la_c^\ast}{\mu_b^\ast+\mu_c^\ast}
&\quad\text{if}\quad
\eta=\ep_{b}+\ep_c
&\quad\text{and}\quad\nu_b+\nu_c>n\,,
\\[12pt]
\,\dfrac{\la_b^\ast}{\mu_b^\ast}
&\quad\text{if}\quad
\eta=2\ts\ep_{b}
&\quad\text{and}\quad2\ts\nu_b>n\,,
\\[12pt]
\,1
&\quad\text{otherwise\ts.}&
\end{array}
\right.
$$
Note that in the first two cases above $1\le b<c\le m\ts$,
while in the third case $1\le b\le m$ and $\f_m=\sp_{2m}\ts$.
Let $v_{\mu}^{\ts\la}$ be the image
of the product
$f_{\ts\overline{1}\ts\nu_1}\ns\ldots\ts
f_{\ts\overline{m}\ts\nu_m}\in\Bb_m$
in the quotient vector space
$\Jb\,\ts\backslash\,\Bb_m\ts/\,\Ib_{\ts\mu,\de_+}\ts$.
This image is a highest vector relative to the action of
the Lie algebra $\g_n$ on this space\ts: it is annihilated by elements
$
E_{\ts ij}-\th_i\ts\th_j\tts E_{\ts\bj\ts\bi}\in\g_n
$
with~$i\prec j\ts$.

\begin{proposition*}
\label{isis}
{\rm\,\,(i)}
The vector $v_{\mu}^{\ts\la}$
is not in the zero coset of
$\Jb\,\ts\backslash\,\Bb_m\ts/\,\Ib_{\ts\mu,\de_+}\ts$.
\\
{\rm\,\,(ii)}
Under the action of\/ $\h$ on\/ 
$\Jb\,\ts\backslash\,\Bb_m\ts/\,\Ib_{\ts\mu,\de_+}$
the vector $v_{\mu}^{\ts\la}$ is of weight $\lambda\ts$.
\\
{\rm\,\,(iii)} For any\/ $\si\in\H_m$ the intertwining operator
\eqref{bbjioper} determined by\/ $\xic_{\ts\si}$ maps the vector\/
$v_{\mu}^{\ts\la}$ to the image in\/
$\Jb\,\ts\backslash\,\Bb_m\ts/\,\Ib_{\ts\si\ts\circ\ts\mu\ts,\ts\si(\de_+)}$
of\/ $ \sih\ts(\ts f_{\ts\overline{1}\ts\nu_1} \ns\ldots\ts
f_{\ts\overline{m}\ts\nu_m} ) \in\Bb_m $
multiplied by the product 
\begin{equation}
\label{isim}
\prod_{\eta\ts\in\De_\si}
\!z_\eta\,.
\end{equation}
\end{proposition*}

\begin{proof}
Part (i) of the proposition follows directly from the definition of
the ideal $\Ib_{\ts\mu,\de_+}\ts$. Let us prove Part (ii). The
elements of $\h$ act on
$\Jb\,\ts\backslash\,\Bb_m\ts/\,\Ib_{\ts\mu,\de_+}$ via their left
multiplication on $\Bb_m\ts$. Let us indicate by $\,\equiv\,$ the
equalities in $\Bb_m$  modulo the left ideal
$\Ib_{\ts\mu,\de_+}\ts$. Then by the definition \eqref{gan} for each
$a=1\lcd m$ we have the relations in the algebra $\Bb_m\ts$,
\begin{gather*}
F_{-\ab,-\ab}\,\ts
f_{\ts\overline{1}\ts\nu_1}
\ns\ldots\ts
f_{\ts\overline{m}\ts\nu_m}
\,=\,
f_{\ts\overline{1}\ts\nu_1}
\ns\ldots\ts
f_{\ts\overline{m}\ts\nu_m}\,F_{-\ab,-\ab}
\,-\sum_{k=1}^n\,
[\,
x_{\ts\overline{a}\ts k}\,\d_{\ts\overline{a}\ts k}\ts-n/2
\,\com\ts
f_{\ts\overline{1}\ts\nu_1}
\ns\ldots\ts
f_{\ts\overline{m}\ts\nu_m}\,]
\\[4pt]
=\,
f_{\ts\overline{1}\ts\nu_1}
\ns\ldots\ts
f_{\ts\overline{m}\ts\nu_m}
\,(\ts F_{-\ab,-\ab}-\nu_a\ts)
\,\equiv\,
f_{\ts\overline{1}\ts\nu_1}
\ns\ldots\ts
f_{\ts\overline{m}\ts\nu_m}
\,(\,\zeta_n(F_{-\ab,-\ab})+\mu_a-\nu_a\ts)
\\[12pt]
\equiv\, x_{\ts\overline{1}\ts
k}^{\,\nu_1}\,\ldots\,x_{\ts\overline{m}\ts k}^{\,\nu_m}
\,(\ts n/2+\mu_a-\nu_a\ts)
\,=\,
\la_a\,
f_{\ts\overline{1}\ts\nu_1}
\ns\ldots\ts
f_{\ts\overline{m}\ts\nu_m}\,.
\end{gather*}
Thus
$$
F_{-\ab,-\ab}\,\ts v_{\mu}^{\ts\la}\,=\,\la_a\,v_{\mu}^{\ts\la}\, \
\quad\text{for}\ \quad a=1\lcd m\ts.
$$

We will prove Part (iii) by induction on
the length of a reduced decomposition
of $\si$ in terms of $\si_1\lcd\si_m\ts$.
If $\si$ is the identity element of $\H_m\ts$,
then the required statement is tautological.
Now suppose that for some $\si\in\H_m$
the statement of (iii) is true.
Take any simple reflection $\si_a\in\H_m$ with $1\le a\le m\ts$,
such that $\si_a\ts\si$ has a longer reduced decomposition
in terms of $\si_1\lcd\si_m$ than $\si\ts$.
If $\f_m=\so_{2m}$ and $a=m\ts$, then we have
$\xic_{\ts\si_m\ts\si}=\sih_m\,\xic_{\ts\si}$
and $\De_{\ts\si_m\ts\si}=\De_{\ts\si}\,$,
so that the induction step is immediate.
We may now assume that $a<m$ in the case $\f_m=\so_{2m}\ts$.

Take the simple root $\eta_{\ts a}$ corresponding to the reflection
$\si_a\ts$. Let $\eta=\si^{-1}(\eta_{\ts a})\ts$. Then
$\eta\in\Delta^+$ and
$$
\si_a\ts\si\ts(\eta)=\si_a(\eta_{\ts a})=-\ts\eta_{\ts a}\notin\De^+.
$$
Hence
$$
\De_{\ts\si_a\ts\si}=\De_{\ts\si}\sqcup\ts\{\eta\tts\}\ts.
$$
Let $\ka\in\h^\ast$ be the weight with labels \eqref{kappa}.
Using the proof of Theorem~\ref{proposition4.5},
we get the equality of two left ideals of the algebra $\Bb_m\ts$,
$$
\Ib_{\ts(\ts\si_a\ts\si\ts)\ts\circ\ts\mu
\ts,\ts
(\ts\si_a\ts\si\ts)\ts(\de_+)}\,=\,
\tilde\I_{\ts(\ts\si_a\ts\si\ts)\ts\circ\ts\ka
\ts,\ts
(\ts\si_a\ts\si\ts)\ts(\de_+)}\ts.
$$
But modulo the second of these two ideals,
the element $H_a$ equals
\begin{gather}
\nonumber
((\ts\si_a\ts\si\ts)\circ\ka\ts)\ts(H_a)\,=\,
(\ts\si_a\ts\si\ts(\ka+\rho)-\rho\ts)\ts(H_a)\,=\,
(\ka+\rho)\ts(\ts\si^{-1}\si_a(H_a))-\rho\ts(H_a)\,=\,
\\[4pt]
\label{rhs}
-\,(\ka+\rho)\ts(\si^{-1}(H_a))-1\,=\,
-\,(\ka+\rho)\ts(H_\eta)-1\,=\,
-\,\frac{\,2\ts(\ts\ka+\rho\com\eta\ts)\,}{(\ts\eta\com\eta\ts)}\ts-1\ts.
\end{gather}
Here $H_\eta=\si^{-1}(H_a)$ is the coroot
corresponding to the root $\eta\,$, and we
use the standard bilinear form on $\h^\ast\ts$.
Using only the definition \eqref{kappa},
the right hand side of \eqref{rhs} equals
\begin{align*}
-\mu_b^\ast+\mu_c^\ast-1
\ \quad\textrm{if}\ \quad
\eta&=\ep_{b}-\ep_c\ts,
\\[5pt]
-\mu_b^\ast-\mu_c^\ast-n-1
\ \quad\textrm{if}\ \quad
\eta&=\ep_{b}+\ep_c\ts,
\\[5pt]
-\mu_b^\ast-n/2-1
\ \quad\textrm{if}\ \quad
\eta&=2\ts\ep_{b}\ts.
\end{align*}
We will now use 
(iii) as the induction assumption.
Denote $\de=\si(\de_+)\ts$. Consider five~cases.

I. Suppose $\eta=\ep_b-\ep_c$ where $1\le b<c\le m\ts$, while
$\si(\ep_b)=\ep_a$ and $\si(\ep_c)=\ep_{a+1}\ts$.
Then $\si_a=\ep_a-\ep_{a+1}$ and
$\de_a=\de_{a+1}=1\ts$. Hence
$$
\sih\ts(\ts
f_{\ts\overline{1}\ts\nu_1}
\ns\ldots\ts
f_{\ts\overline{m}\ts\nu_m})
=
f_{\ts\overline{a}\ts\nu_b}\,
f_{\ts\overline{a+1}\ts\nu_c}Y
$$
where $Y$ is an element of the subalgebra
of $\PD\ts(\CC^{\ts m}\ot\CC^{\ts n})$
generated by all $x_{dk}$ and $\d_{dk}$
with $d\neq\overline{a}\com\overline{a+1}\ts$.
Here Lemma \ref{xxnorma} with $s=\nu_b$ and
$t=\nu_c$ applies. With these $s$~and~$t\ts$, by substituting
$-\mu_b^\ast+\mu_c^\ast-1$ for $H_a$
in the fraction displayed
in that lemma, the fraction becomes 
\begin{equation}
\label{sampro}
\frac{-\mu_b^\ast+\mu_c^\ast-1+\nu_b-\nu_c+1}
{-\mu_b^\ast+\mu_c^\ast-1+1}\,=\,\frac{\la_b^\ast-\la_c^\ast}
{\mu_b^\ast-\mu_c^\ast}\,.
\end{equation}
The condition $s>t$ from Lemma \ref{xxnorma} means here that $\nu_b>\nu_c\ts$.

II. Suppose $\eta=\ep_b-\ep_c$ where $1\le b<c\le m\ts$, but
$\si(\ep_b)=-\ts\ep_{a+1}$ and $\si(\ep_c)=-\ts\ep_{a}\ts$.
Then $\si_a=\ep_a-\ep_{a+1}$ again, but
$\de_a=\de_{a+1}=-\ts1\ts$. Hence
$$
\sih\ts(\ts
f_{\ts\overline{1}\ts\nu_1}
\ns\ldots\ts
f_{\ts\overline{m}\ts\nu_m})
=
g_{\ts\overline{a}\ts\nu_c}\,
g_{\ts\overline{a+1}\ts\nu_b}Y
$$
where $Y$ is another element of the subalgebra
of $\PD\ts(\CC^{\ts m}\ot\CC^{\ts n})$
generated by all $x_{dk}$ and $\d_{dk}$
with $d\neq\overline{a}\com\overline{a+1}\ts$.
Now Lemma \ref{ddnorma} with $s=\nu_{\ts c}$ and
$t=\nu_b$ applies. With these $s$~and~$t\ts$,
by substituting $-\mu_b^\ast+\mu_c^\ast-1$ for $H_a$
in the fraction displayed
in Lemma \ref{ddnorma}, the fraction becomes
the same number \eqref{sampro} as in the previous case,
under the same condition $\nu_b>\nu_c\ts$.

III. Suppose $\eta=\ep_b+\ep_c$ and $1\le b<c\le m\ts$, while
$\si(\ep_b)=\ep_a$ and $\si(\ep_c)=-\ts\ep_{a+1}\ts$.
Then $\si_a=\ep_a-\ep_{a+1}$ again, but
$\de_a=1$ and $\de_{a+1}=-1\ts$. Hence
$$
\sih\ts(\ts
f_{\ts\overline{1}\ts\nu_1}
\ns\ldots\ts
f_{\ts\overline{m}\ts\nu_m})
=
f_{\ts\overline{a}\ts\nu_b}\,
g_{\ts\overline{a+1}\ts\nu_c}Y
$$
where $Y$ is another element of the subalgebra
of $\PD\ts(\CC^{\ts m}\ot\CC^{\ts n})$
generated by $x_{dk}$ and $\d_{dk}$
with $d\neq\overline{a}\com\overline{a+1}\ts$.
Here Lemma \ref{xdnorma} with $s=\nu_b$ and
$t=\nu_c$ applies. With these $s$~and~$t\ts$,
by substituting $-\mu_b^\ast-\mu_c^\ast-n-1$ for $H_a$
in the fraction displayed in that lemma,
the fraction becomes the number
\begin{equation}
\label{samsam}
\frac{-\mu_b^\ast-\mu_c^\ast-n-1+\nu_b+\nu_c+1}
{-\mu_b^\ast-\mu_c^\ast-n-1+n+1}
\,=\,
\frac
{\lambda_{\ts b}^\ast+\lambda_c^\ast}
{\,\mu_{\ts b}^\ast+\mu_c^\ast}\,.
\end{equation}
The condition $s+t>n$ from Lemma \ref{xdnorma} means 
here that $\nu_{\ts b}+\nu_{\ts c}>n\ts$.

IV. Suppose $\eta=\ep_b+\ep_c$ where $1\le b<c\le m\ts$, but
$\si(\ep_b)=-\ts\ep_{a+1}$ and $\si(\ep_c)=\ep_{a}\ts$.
Then $\si_a=\ep_a-\ep_{a+1}$ again, but
$\de_a=1$ and $\de_{a+1}=-\ts1\ts$. Hence
$$
\sih\ts(\ts
f_{\ts\overline{1}\ts\nu_1}
\ns\ldots\ts
f_{\ts\overline{m}\ts\nu_m})
=
f_{\ts\overline{a}\ts\nu_c}\,
g_{\ts\overline{a+1}\ts\nu_b}Y
$$
where $Y$ is another element of the subalgebra
of $\PD\ts(\CC^{\ts m}\ot\CC^{\ts n})$
generated by $x_{dk}$ and $\d_{dk}$
with $d\neq\overline{a}\com\overline{a+1}\ts$.
Now Lemma \ref{xdnorma} with $s=\nu_c$ and
$t=\nu_b$ applies. With these $s$~and~$t\ts$,
by substituting
$-\mu_b^\ast-\mu_c^\ast-n-1$ for $H_a$
in the fraction displayed in that lemma,
the fraction becomes the same number \eqref{samsam} as in the previous case,
under the same condition $\nu_{\ts b}+\nu_{\ts c}>n\ts$.

V. Suppose $\f_m=\sp_{2m}$ and $\eta=2\ts\ep_{\ts b}$
with $1\le b\le m\ts$. Then $\si(\ep_b)=\ep_{\ts m}$
and $\si_a=\si_m\ts$, while $\de_m=1\ts$. Hence
$$
\sih\ts(\ts
f_{\ts\overline{1}\ts\nu_1}
\ns\ldots\ts
f_{\ts\overline{m}\ts\nu_m})
=
f_{\ts\overline{m}\ts\nu_b}\,Y
$$
where $Y$ is now an element of the subalgebra
of $\PD\ts(\CC^{\ts m}\ot\CC^{\ts n})$
generated by $x_{dk}$~and~$\d_{dk}$
with $d\neq\overline{m}=1\ts$.
Here Lemma \ref{xnorma} with $s=\nu_b$ applies.
With this $s\ts$,
by substituting
$-\mu_{\ts b}^\ast-{n}/{2}-1$ for $H_m$ in
in the fraction displayed in that lemma, the fraction becomes 
$$
\frac{-\mu_b^\ast-n/2-1+\nu_b+1}
{-\mu_b^\ast-n/2-1+{n}/{2}+1}
\,=\,
\frac{\lambda_{\ts b}^\ast}{\mu_{\ts b}^\ast}
\,.
$$
The condition $s>n/2$ from Lemma \ref{xnorma}
means here that $2\ts\nu_{\ts b}>n\ts$.

Thus in all the five cases above, by using the induction assumption, the
intertwining operator
$$
\Jb\,\ts\backslash\,\Bb_m\ts/\,\Ib_{\ts\mu,\de_+}
\to\,
\Ib_{\ts(\ts\si_a\ts\si\ts)\ts\circ\ts\mu
\ts,\ts
(\ts\si_a\ts\si\ts)\ts(\de_+)}
$$
determined by\/ $\xic_{\ts\si_a\ts\si}$ maps the vector
$v_{\mu}^{\ts\la}$ to the image in
$
\,\Jb\,\ts\backslash\,\Bb_m\ts/\,
\Ib_{\ts(\ts\si_a\ts\si\ts)\ts\circ\ts\mu
\ts,\ts
(\ts\si_a\ts\si\ts)\ts(\de_+)}
$
of 
$$
\sih_a\ts\sih\ts(\ts
f_{\ts\overline{1}\ts\nu_1}
\ns\ldots\ts
f_{\ts\overline{m}\ts\nu_m})
\in\Bb_m
$$
multiplied by the product \eqref{isim} over the set
$\De_\si\ts$, and by an extra factor $z_\eta$
corresponding to the positive root
$\eta=\si^{-1}(\eta_{\ts a})\ts$.
This makes the induction step.
\qed
\end{proof}

The product \eqref{isim} in Proposition \ref{isis}
does not depend on the choice of a reduced decomposition of $\si\in\H_m$
in terms of $\si_1\lcd\si_m\ts$.
The uniqueness of the
intertwining operator \eqref{ppi} thus provides another proof
of the independence of our operator \eqref{distoperla}
on the decomposition of $\si\ts$,
not involving Proposition~\ref{p3}.
Proposition~\ref{isis} also shows that our
intertwining operator \eqref{distoperla} is not zero.


\section*{\normalsize 6. Olshanski homomorphism}
\setcounter{section}{6}
\setcounter{equation}{0}
\setcounter{theorem*}{0}

For a positive integer $l\ts$, take the vector space
$\CC^{\ts n+l}$.
In the case of an alternating form on
$\CC^{\ts n}$ choose $l$ to be even.
Let $e_1\lcd e_{\ts n+l}$
be the vectors of the standard basis in $\CC^{\ts n+l}$.
Consider the decomposition
$\CC^{\ts n+l}=\CC^{\ts n}\op\CC^{\ts l}$
where the direct summands $\CC^{\ts n}$ and $\CC^{\ts l}$
are spanned by the vectors $e_1\lcd e_n$ and $e_{\ts n+1}\lcd e_{\ts n+l}$
respectively. This defines an embedding of the direct sum
$\gl_n\op\ts\gl_{\ts l}$ of Lie algebras to $\gl_{\ts n+l}\ts$.
As a subalgebra of $\gl_{\ts n+l}\ts$,
the summand $\gl_n$ is spanned by the matrix units
$E_{ij}\in\gl_{\ts n+l}$ where
$i\com j=1\lcd n\ts$. The summand $\gl_{\ts l}$ is spanned by
the matrix units $E_{ij}$ where $i\com j=n+1\lcd n+l\ts$.

The subspace $\CC^{\ts n}\subset\CC^{\ts n+l}$ comes
with a bilinear form chosen in Section 1. Now
choose a bilinear form on the subspace $\CC^{\ts l}\subset\CC^{\ts n+l}$
in a similar way. Namely, let
$i$ be any of the indices $n+1\lcd n+l\ts$. If $i-n$ is even, then put
$\bi=i-1\ts$. If $i-n$ is odd and $i<n+l$,  then put $\bi=i+1\ts$.
If $i=n+l$ and $l$ is odd, then put $\bi=i\ts$.
Further, put $\th_i=1$ or $\th_i=(-1)^{\ts i-n-1}$
in the case of the symmetric or alternating form on $\CC^{\ts n}$.
For any basis vectors $e_i$ and $e_j$ of the subspace $\CC^{\ts l}$ put
$\langle\ts e_i\com e_j\ts\rangle=\th_i\,\de_{\ts\bi j}\ts$.
Equip the vector space $\CC^{\ts n+l}$ with the
bilinear form which is the sum of the forms on the direct summands.
The forms on $\CC^{\ts l}$ and $\CC^{\ts n+l}$ are of the same type
(symmetric or alternating) as the form on $\CC^n\ts$.

Now consider the subalgebras $\g_n\com\ts\g_{\ts l}$ and $\g_{n+l}$ of
the Lie algebras $\gl_n\com\ts\gl_{\ts l}$ and $\gl_{\ts n+l}$ respectively.
We have an embedding of the direct sum $\g_n\op\g_l$
to the Lie algebra $\g_{n+l}\ts$,
according to our choice of the bilinear forms made above.
We also have an embedding
of the direct product of Lie groups $G_n\times G_l$ to $G_{n+l}\ts$.
Let $\C_l$ denote the subalgebra of $G_l\ts $-invariants
in the universal enveloping algebra $\U(\g_{n+l})\ts$.
Then $\C_l$ contains the subalgebra $\U(\g_n)\subset\U(\g_{n+l})\ts$.
If $\g_n=\sp_n$ then $\C_l$ coincides with the
centralizer of the subalgebra $\U(\sp_l)\subset\U(\sp_{n+l})\ts$.
If $\g_n=\so_n$ then $\C_l$ is contained in the
centralizer of $\U(\so_l)\subset\U(\so_{n+l})\ts$,
but may not coincide with the centralizer.

Take the extended twisted Yangian $\X(\g_{n+l})\ts$.
The subalgebra of $\X(\gl_{n+l})$ generated~by
$$
S_{ij}^{\ts(1)},S_{ij}^{\ts(2)},\ts\ldots
\quad\text{where}\quad
i\com j=1\lcd n
$$
is isomorphic to $\X(\g_n)\ts$ as an associative algebra, see
\cite[Section 3.14]{MNO}.  Thus we have a natural embedding
$\X(\g_n)\to\X(\g_{n+l})\ts$, let us denote it by $\io_{\ts l}\ts$.
We also have a surjective homomorphism
$$
\pi_{n+l}:\ts\X(\g_{n+l})\to\U(\g_{n+l})\ts,
$$
see \eqref{pin}. Note that the composition
$\pi_{n+l}\,\io_{\ts l}$ coincides with the homomorphism~$\pi_n\ts$.

Further, consider the involutive automorphism $\om_{n+l}$
of the algebra $\X(\g_{n+l})\ts$, see the definition \eqref{sin}.
The image of the composition of homomorphisms
\begin{equation}
\label{ohom}
\pi_{n+l}\,\ts\om_{n+l}\,\ts\io_{\ts l}:\ts
\X(\g_n)\to
\U(\g_{n+l})
\end{equation}
belongs to subalgebra $\C_l\subset\U(\g_{n+l})\ts$.
Moreover, together with the subalgebra of
$G_{n+l}\ts$-invariants in $\U(\g_{n+l})\ts$,
this image generates
$\C_l\ts$. These two results are due to G.\,Olshanski \cite{O2},
for their detailed proofs see \cite[Section 4]{MO}.
We will use the composition of the homomorphisms
$$
\ga_{\tts l}=
\pi_{n+l}\,\ts
\om_{n+l}\,\ts
\io_{\ts l}\,\ts
\om_n\ts.
$$
We will call it the \textit{Olshanski homomorphism\/}.
The images of the homomorphisms $\ga_{\tts l}$ and \eqref{ohom}
in $\U(\g_{n+l})$ coincide.
The reason for using the homomorphism $\ga_{\tts l}$ rather than
the homomorphism \eqref{ohom} will become apparent when we state
Theorem~\ref{5.1}.

An irreducible representation of the group $G_n$ is called
\textit{polynomial\/} if it appears as a subrepresentation
of some tensor power of the defining representation $\CC^{\ts n}$.
According to \cite[Sections V.7 and VI.3]{W}
the irreducible polynomial representations of the group $G_n$
are parameterized
by all the partitions $\nu$ of $N=0\com1\com2\com\,\ldots$
such that $2\ts\nus_1\le n$ in the case $G_n=Sp_n\ts$, and
$\nus_1+\nus_2\le n$ in the case $G_n=O_n\ts$.
Here $\nus$ is the partition conjugate to $\nu$
while $\nus_1\ts\com\nus_2\ts\com\,\ldots$ are the parts of $\nus$.
Note that in the case $G_n=O_n$ we still have
$2\ts\nus_2\le n\ts$.
Denote by $W_\nu$ the irreducible polynomial representation
of the group $G_n$ corresponding to $\nu\ts$.
Let $\nu_1\ts\com\nu_2\ts\com\,\ldots$ be the parts of $\nu\ts$.

\vbox{
Let $\nut$ be the weight of the Lie algebra
$\f_m$ with the sequence of labels
$$
(\ts{n}/2-\nus_m\ts\lcd\ts{n}/2-\nus_1\ts)\ts.
$$
Due to conditions on $\nu$,
the labels $\nut_1\lcd\nut_m$ of $\nut$ in the case $\f_m=\sp_{2m}$
are integers such that $\nut_1\ge\ldots\ge\nut_m\ge0\ts$.
In the case $\f_m=\so_{2m}$ either all labels of $\nut$ are integers,
or all if them are half-integers. In the case $\f_m=\so_{2m}$ we have
$\nut_1\ge\ldots\ge\nut_{m-1}\ge|\ts\nut_m|\ts$.
}

Consider $\P\ts(\CC^{\ts m}\ot\CC^{\ts n})$ as a bimodule
over $\f_m$ and $G_n\ts$. Then by \cite[Subsection 3.8.9]{H}
when $G_n=Sp_n\ts$, or by \cite[Subsection 4.3.5]{H}
when $G_n=O_n\ts$, we have a 
decomposition
\begin{equation}
\label{dirsum}
\P\ts(\CC^{\ts m}\ot\CC^{\ts n})\,=\,\,
\mathop{\op}\limits_\nu\,L_{\ts\nut}\ot W_\nu
\end{equation}
where $\nu$ ranges over all parameters of
irreducible polynomial representations of $G_n$ such that $\nu_1\le m\ts$.
Here $L_{\ts\nut}$ is the irreducible $\f_m\ts$-module
of the highest weight $\nut\ts$.

Let $\la$ and $\mu$ be parameters of any
irreducible polynomial representations of the groups $G_{n+l}$ and $G_l$
respectively. Suppose that $\la_1\com\ts\mu_1\le m\ts$.
Using the action of the group $G_l$ on $W_\la$
via its embedding to $G_{n+l}$
as the second direct factor of the subgroup
$G_n\times G_l$ consider the vector space
\begin{equation}
\label{hll}
\Hom_{\,G_l}(\ts W_\mu\ts,W_\la)\ts.
\end{equation}
The subalgebra $\C_l\subset\U(\g_{n+l}\ts)$ acts on this vector
space through the action of $\U(\g_{n+l})$ on $W_\la\ts$.
In the case $G_n=Sp_n\ts$,
the vector space \eqref{hll} is irreducible under the action of the
algebra $\C_l\ts$; see \cite[Theorem 9.1.12]{D}.
In the case $G_n=O_n\ts$, the $\C_l\ts$-module \eqref{hll}
is either irreducible or splits to a direct
sum of two irreducible $\C_l\ts$-modules.
It is irreducible if $W_\la$ is irreducible
as a $\so_{n+l}\ts$-module, that is if $2\ts\las_1\neq n+l$
by \cite[Section V.9]{W}.
Note that in the case $G_n=O_n\ts$, the condition $2\ts\las_1\neq n+l$
is sufficient but not necessary for the irreducibility
of the $\C_l\ts$-module \eqref{hll}\ts; see \cite[Section 1.7]{N2}.

In any case, the vector space
\eqref{hll} is irreducible under joint action of
the subalgebra $\C_l\subset\U(\g_{n+l})$ and of
the subgroup $G_n\subset G_{n+l}\ts$; see again \cite[Section 1.7]{N2}.
Hence the following identifications of bimodules over
$C_l$ and $G_n$ are unique up to rescaling of their
vector spaces\ts:
$$
\Hom_{\,G_l}(\ts W_\mu\ts,W_\la)\ts=
$$
$$
\Hom_{\,G_l}(\ts W_\mu\ts,
\Hom_{\,\f_m}(\ts L_{\ts\lat}\,,\P\ts(\ts\CC^{\ts m}\ot\CC^{\ts n+l}\ts)))\ts=
\vspace{2pt}
$$
$$
\Hom_{\,G_l}(\ts W_\mu\ts,
\Hom_{\,\f_m}(\ts L_{\ts\lat}\,,
\P\ts(\ts\CC^{\ts m}\ot\CC^{\ts l}\ts)\ot
\P\ts(\ts\CC^{\ts m}\ot\CC^{\ts n}\ts)))\ts=
\vspace{2pt}
$$
\begin{equation}
\label{hlls}
\,\Hom_{\,\f_m}(\ts L_{\ts\lat}\,,
L_{\ts\mut}\ot\P\ts(\ts\CC^{\ts m}\ot\CC^{\ts n}\ts))\ts.
\vspace{4pt}
\end{equation}
We use the decompositions \eqref{dirsum} for $n+l$
and $l$ instead of $n\ts$, and the identification
\begin{equation}
\label{pmnl}
\P\ts(\ts\CC^{\ts m}\ot\CC^{\ts n+l}\ts)=
\P\ts(\ts\CC^{\ts m}\ot\CC^{\ts l}\ts)\ot
\P\ts(\ts\CC^{\ts m}\ot\CC^{\ts n}\ts)
\end{equation}
of vector spaces. Thus in \eqref{hlls},
the labels of the weights $\lat$ and $\mut$ of $\ts\f_m$
are respectively
$$
(\ts n/2+l\ts/2-\las_m\ts\lcd\ts n/2+l\ts/2-\las_1\ts)
\ \quad\text{and}\ \quad
(\ts{l}\ts/2-\mus_m\ts\lcd\ts{l}\ts/2-\mus_1\ts)\ts.
$$

By pulling back via the Olshanski
homomorphism $\ga_{\tts l}:\X(\g_n)\to\C_l\ts$,
the vector space \eqref{hll} becomes a module over the
extended twisted Yangian $\X(\g_n)\ts$. Using the above identifications,
the vector space \eqref{hlls} than also becomes a module over $\X(\g_n)\ts$.
But the target $\f_m\ts$-module
$L_{\ts\mut}\ot\P\ts(\ts\CC^{\ts m}\ot\CC^{\ts n}\ts)$
in \eqref{hlls} coincides with the $\f_m\ts$-module
$\F_m(L_{\ts\mut})\ts$.

\begin{theorem*}
\label{5.1}
The action of\/ $\X(\g_n)$ on the vector space \eqref{hlls}
via the homomorphism\/ $\ga_{\tts l}$ coincides with the action, obtained by
pulling the action of\/ $\X(\g_n)$ on the bimodule\/
$\F_m(L_{\ts\mut})$ back through the homomorphism \eqref{fus} where
\begin{equation}
\label{fuss}
f(u)=1\,-\,m\ts(\ts u-l\ts/2\pm1/2\ts)^{-1}\,.
\end{equation}
\end{theorem*}

\begin{proof}
Take the action of the subalgebra $\C_{\ts l}\subset\U(\gl_{n+l})$
on the space $\P\ts(\ts\CC^{\ts m}\ot\CC^{\ts n+l}\ts)\ts$.
The extended twisted Yangian $\X(\g_n)$ acts on this vector space via
the homomorphism
$\ga_{\tts l}:\X(\g_n)\to\C_{\ts l}\ts$. Using the decomposition \eqref{pmnl}
we will show that for $i\com j=1\lcd n$ the generators
$
S_{ij}^{\ts(1)},S_{ij}^{\ts(2)},\ts\ldots
$
of $\X(\g_n)$ act on this
vector space respectively as the coefficients at
$u^{-1}\com u^{-2}\com\ts\ldots\,$
of the series \eqref{fhom} multiplied by the series
\eqref{fuss}.

For any $i\com j=1\lcd n+l$
the element $F_{ij}\in\U(\g_{n+l})$ acts on
$\P(\ts\CC^{\ts m}\ot\CC^{\ts n+l}\ts)$ as the operator
$$
\sum_{c=1}^m\,\,
(\ts x_{ci}\,\d_{\ts cj}-
\ts\th_i\,\th_j\,x_{c\ts\bj}\,\d_{\ts c\ts\bi}\ts)\,.
$$
Here we use the standard coordinate functions
$x_{ci}$ on $\CC^{\ts m}\ot\CC^{\ts n+l}$
with $c=1\lcd m$ and $i=1\lcd n+l\ts$.
Then $\d_{ci}$ is the left derivation on the Grassmann algebra
$\P\ts(\CC^{\ts m}\ot\CC^{\ts n+l})$ relative to $x_{ci}\ts$.
The functions $x_{ci}$ with $c\leqslant n$ and $c>n$
correspond to the direct summands $\CC^{\ts n}$ and
$\CC^{\ts l}$ of $\CC^{\ts n+l}$.
Consider the $(n+l)\times(n+l)$ matrix whose $i\com j$ entry is
$$
\de_{ij}\,+\,(\ts u-l\ts/2\pm1/2\ts)^{-1}\,
\sum_{c=1}^m\,\,
(\ts x_{ci}\,\d_{\ts cj}-
\ts\th_i\,\th_j\,x_{c\ts\bj}\,\d_{\ts c\ts\bi}\ts)\,.
$$
Write this matrix and its inverse as the block matrices
$$
\begin{bmatrix}
\,A\,&B\,\\\,C\,&D\,
\end{bmatrix}
\quad\textrm{and}\quad
\begin{bmatrix}\ \At\,&\Bt\ \\ \ \Ct\,&\Dt\
\end{bmatrix}
$$
where the blocks $A,B,C,D$ and $\At,\Bt,\Ct,\Dt$ are matrices of sizes
$n\times n$, $n\times l$, $l\times n$, $l\times l$ respectively.
The action of the algebra $\X(\g_n)$ on the vector space
$\P(\ts\CC^{\ts m}\ot\CC^{\ts n+l}\ts)$ via the homomorphism
$\ga_{\tts l}:\X(\gl_n)\to\C_{\ts l}$ can now be described by assigning
to the series $S_{ij}(u)$ with $i\com j=1\lcd n$ the
$i\com j$ entry of the matrix $\At^{\,-1}\ts$.

Introduce the $(n+l)\times 2\ts m$ matrix whose $i\com c$ entry
for $c=-\ts m\lcd-1$ is the operator of the left multiplication by
$x_{ci}$ on $\P(\ts\CC^{\ts m}\ot\CC^{\ts n+l}\ts)\ts$.
For $c=1\lcd m$ let the $i\com c$ entry of this matrix be
the operator $\th_i\,\d_{c\ts\bi}\,$. Write this matrix as
$$
\begin{bmatrix}
\,P\,
\\
\,\Pb\,
\end{bmatrix}
$$
where the blocks $P$ and $\Pb$ are matrices of sizes
$n\times2\ts m$ and $l\times2\ts m$ respectively. Further,
introduce the $2\ts m\times(n+l)$ matrix whose $c\com j$ entry
for $c=-\ts m\lcd-1$ is the operator $\d_{cj}\ts$.
For $c=1\lcd m$ let the $c\com j$ entry of this matrix be
the operator of left multiplication by $\th_j\,x_{c\ts\bj}\,$.
Write this matrix as
$$
\begin{bmatrix}
\,Q\,\,\Qb\,\ts
\end{bmatrix}
$$
where $Q$ and $\Qb$ are matrices of sizes
$2\ts m\times n$ and $2\ts m\times l$ respectively. Then
$$
\begin{bmatrix}
\,A\,&B\,\\\,C\,&D\,
\end{bmatrix}
\,=\,
1\,+\,(\ts u-l\ts/2\pm1/2\ts)^{-1}
\begin{bmatrix}\,
\,P\ts Q-m\,&\,P\ts\Qb
\\
\Pb\ts Q\,&\,\Pb\ts \Qb-m\,
\end{bmatrix}
$$
which can be also written as the matrix
$$
1\,+\,(\ts u-l\ts/2\pm1/2-m\ts)^{-1}
\begin{bmatrix}\,
\,P\ts Q&\ P\ts\Qb\
\\
\,\Pb\ts Q&\ \Pb\ts \Qb\
\end{bmatrix}
$$
multiplied by the series $f(u)$ determined by \eqref{fuss}.
Using a well known formula for $\At^{\,-1}\ts$,
$$
\At^{\,-1}=A-B\ts D^{-1}\ts C
\,=\,f(u)\,\bigl(\ts
1+(\ts u-l\ts/2\pm1/2-m\ts)^{-1}\ts P\ts Q
$$
$$
-\,\,
(\ts u-l\ts/2\pm1/2-m\ts)^{-2}\,P\ts\Qb\,\ts\bigl(\ts1+
(\ts u-l\ts/2\pm1/2-m\ts)^{-1}\ts\Pb\ts\Qb\,\bigr)^{-1}\ts\Pb\ts Q
\,\bigr)
\vspace{2pt}
$$
\begin{equation}
\label{pulq}
=\,f(u)\,\bigl(\ts1+P\ts
(\ts u-l\ts/2\pm1/2-m+\Qb\ts\Pb\,\bigr)^{-1}\ts Q\,\bigr)\,.
\vspace{4pt}
\end{equation}
Consider the $2m\times2m$ matrix $\Qb\ts\Pb$ appearing in the last line.
For any indices
$a\com b=-\ts m\lcd-1\com1\lcd m$ the $a\com b$ entry of this matrix
is the operator
$$
\de_{ab}\,l\ts/2\ts+\,\bar\zeta_{\ts l}\ts(F_{ab})
$$
where $\bar\zeta_{\ts l}:\U(\ts\f_m)\to\PD\ts(\CC^{\ts m}\ot\CC^{\ts n+l})$
is the homomorphism corresponding to the action of the
Lie algebra $\f_m$ on $\P\ts(\CC^{\ts m}\ot\CC^{\ts n+l})$
via the tensor factor $\P\ts(\CC^{\ts m}\ot\CC^{\ts l})$ in \eqref{pmnl},
similar to the homomorphism \eqref{gan}.
Namely for $a\com b=1\lcd m$ we have
$$
\bar\zeta_{\ts l}\ts(F_{ab})\,=\,
-\ts\de_{ab}\,l\ts/2\ +
\sum_{k=n+1}^{n+l}\,x_{ak}\,\d_{\ts bk}\,,
$$
$$
\bar\zeta_{\ts l}\ts(F_{a,-b})\,=\,\sum_{k=n+1}^{n+l}\,
\th_k\,x_{a\bk}\,x_{\ts bk}\,,
\qquad
\bar\zeta_{\ts l}\ts(F_{-a,b})\,=\sum_{k=n+1}^{n+l}\,
\th_k\,\d_{ak}\,\d_{\ts b\bk}\,.
$$
Hence any entry of the $2m\times2m$ matrix
$$
(\ts u-l\ts/2\pm1/2-m+\Qb\ts\Pb\,)^{-1}
$$
can be obtained by applying the homomorphism $\bar\zeta_{\ts l}$
to the respective entry
of the matrix $F\ts(u\pm\frac12-m)\ts$;
the latter entries are series in $u^{-1}$
with coefficients in $\U(\ts\f_m)\ts$.
We now complete the proof by
comparing the $i\com j$ entry of the $n\times n$ matrix \eqref{pulq}
with the series, obtained from \eqref{fhom}
by replacing $F_{ab}(u\pm\frac12-m)$ there by
$\bar\zeta_{\ts l}\ts(F_{ab}(u\pm\frac12-m))$ for all indices
$a\com b=-\ts m\lcd-1\com1\lcd m\ts$.
\qed
\end{proof}

Set $\C_0=\U(\g_n)$ and $\ga_{\ts0}=\pi_n\ts$.
Then Theorem \ref{5.1} remains
valid in the case $l=0\ts$. In this case we assume that $\g_{\ts l}=\{0\}\ts$.
Note that our proof of Theorem \ref{5.1} also implies Proposition~\ref{xb},
because the kernels of homomorphisms $\bar\zeta_{\ts l}$ with
$l=0\com1\com2\com\,\ldots$ have only zero intersection.
For $\f_m=\so_2$ the latter follows directly from
the definition \eqref{gan}.
For $\f_m\neq\so_2$
all irreducible finite-dimensional $\f_m\ts$-modules
arise from the skew Howe~duality.

Let $\la$ and $\mu$ be the parameters of any
irreducible polynomial representations of $G_{n+l}$ and $G_l$
respectively. The vector space \eqref{hll} is not zero if and only if
\begin{equation}
\label{llcon}
\la_k\ge\mu_k
\quad\textrm{and}\quad
\las_k-\mus_k\le n
\quad\textrm{for every}\quad
k=1\com2\com\ts\ldots\ts;
\end{equation}
see \cite[Section 1.3]{N2}.
Suppose that $\la_1\com\mu_1\le m\ts$.
Then we can identify
the vector spaces \eqref{hll} and \eqref{hlls}.
Then the algebra $\C_l$ acts on \eqref{hlls} irreducibly, if $G_n=Sp_n\ts$.
If $G_n=O_n\ts$, then \eqref{hlls} is irreducible
under the joint action of the algebra $\C_l$ and the group $O_n\ts$.
In both cases, the $G_{n+l}\ts$-invariant elements of $\U(\gl_{n+l})$
act on \eqref{hlls} via multiplication by scalars.
Then Theorem \ref{5.1} has a corollary,
which refers to the
action of $\X(\g_n)$ on the vector space \eqref{hlls}
inherited from the bimodule\/ $\F_m(L_{\ts\mut})\ts$.

\begin{corollary*}
\label{xirred}
The algebra\/ $\X(\g_n)$ acts on space \eqref{hlls} irreducibly,
if\/ $G_n=Sp_n\ts$.
If $G_n=O_n\ts$, the space \eqref{hlls} is irreducible
under the joint action of\/ $\X(\g_n)$ and\/ $O_n\ts$.
\end{corollary*}

Now suppose that $\f_m\neq\so_2\ts$.
Then any irreducible finite-dimensional module $V$ of $\f_m$
is equivalent to $L_{\ts\mut}$ for some
non-negative integer $l$ and the label $\mu$
of some irreducible polynomial representation of the group $G_l$
with $\mu_1\le m\ts$. If $V^{\ts\prime}$ is another
irreducible finite-dimensional $\f_m\ts$-module,
such that the vector space \eqref{vvp} is non zero, then
$V^{\ts\prime}$ has to be equivalent to $L_{\ts\lat}$ for the label $\la$
of some irreducible polynomial representation of $G_{n+l}$ with
$\la_1\le m\ts$. Thus any non-zero vector space \eqref{vvp}
has to be of the form \eqref{hlls}.


\section*{Acknowledgments}
We are gratefil to P.\,Kulish for amiable attention to this work.
The first author has been supported by
the RFBR grant 08-01-02934,
the grant for Support of Scientific Schools 8065-2006-2,
by the Atomic Energy Agency of the Russian Federation,
and by the ANR grant 05-BLAN-0029-01.
The second author has been supported by the EPSRC grant C511166,
and by the EC grant MRTN-CT2003-505078.
This work began when both authors visited
the Max Planck Institute for Mathematics in Bonn.
We are grateful to the staff of the institute for
their kind help and generous hospitality.




\begin{thebibliography}{[EHW]}

\bibitem[A]{A}
{T.\,Arakawa},
\textit{Drinfeld functor and finite-dimensional representations of the
Yangian},
{Commun. Math. Phys.}
\textbf{205}
(1999),
1--18.

\bibitem[AS]{AS}
{T.\,Arakawa and T.\,Suzuki},
\textit{Lie algebras and degenerate affine Hecke algebras of type $A$\/},
{J. Algebra}
\textbf{209}
(1998),
288--304.

\bibitem[AST]{AST}
{T.\,Arakawa, T.\,Suzuki and A.\,Tsuchiya},
\textit{Degenerate double affine Hecke algebras and conformal field theory},
{Progress Math.}
\textbf{160}
(1998),
1--34.

\bibitem[C]{C}
{I.\,Cherednik},
\textit{Lectures on Knizhnik-Zamolodchikov
equations and Hecke algebras},
{Math. Soc. Japan Memoirs}
{\bf 1}
(1998),
1--96.

\bibitem[D]{D}
J.\,Dixmier,
\textit{Alg\`ebres enveloppantes},
Gauthier-Villars, Paris, 1974.

\bibitem[D1]{D1}
{V.\,Drinfeld},
\textit{Hopf algebras and the quantum Yang-Baxter equation},
{Soviet Math.\,Dokl.}
\textbf{32}
(1985),
254--258.

\bibitem[D2]{D2}
{V.\,Drinfeld},
\textit{Degenerate affine Hecke algebras and Yangians},
{Funct. Anal. Appl.}
\textbf{20}
(1986),
56--58.

\bibitem[H]{H}
{R.\,Howe},
\textit{Perspectives on invariant theory:
Schur duality, multiplicity-free actions and beyond},
Israel Math. Conf. Proc.
\textbf{8}
(1995),
1--182.

\bibitem[KN1]{KN1}
S.\,Khoroshkin and M.\,Nazarov,
\textit{Yangians and Mickelsson algebras I\/},
Transformation Groups
\textbf{11}
(2006),
625--658.

\bibitem[KN2]{KN2}
S.\,Khoroshkin and M.\,Nazarov,
\textit{Yangians and Mickelsson algebras II\/},
Moscow Math. J.
\textbf{6}
(2006),
477--504.

\bibitem[KN3]{KN3}
S.\,Khoroshkin and M.\,Nazarov,
\textit{Twisted Yangians and Mickelsson algebras I\/},
{Selecta Math.}
\textbf{13}
(2007),
69--136.

\bibitem[KO]{KO}
{S.\,Khoroshkin and O.\,Ogievetsky},
\textit{Mickelsson algebras and Zhelobenko operators},
{J. Algebra}
\textbf{319}
(2008),
2113--2165.

\bibitem[KS]{KS}
{P. Kulish and E. Sklyanin},
\textit{Algebraic structures related to reflection equations},
{J. Phys.}
\textbf{A25}
(1992),
5963--5975.

\bibitem[M1]{M1}
{J.\,Mickelsson},
\textit{Step algebras of semi-simple subalgebras of Lie algebras},
Reports Math. Phys.
\textbf{4}
(1973),
307--318.

\bibitem[M2]{M2}
{J.\,Mickelsson},
\textit{On irreducible modules of a Lie algebra which are composed of
finite-dimensional modules of a subalgebra},
{Ann. Acad. Sci. Fenn. Ser. A I}
\textbf{598}
(1975),
1--16.

\bibitem[M]{M}
{A.\,Molev},
\textit{Skew representations of twisted Yangians},
Selecta Math.
\textbf{12}
(2006),
1--38.

\bibitem[MNO]{MNO}
{A.\,Molev, M.\,Nazarov and G.\,Olshanski},
\textit{Yangians and classical Lie algebras},
{Russian Math.\ Surveys}
\textbf{51}
(1996),
205--282.

\bibitem[MO]{MO}
{A.\,Molev and G.\,Olshanski},
\textit{Centralizer construction for twisted Yangians},
{Selecta Math.}
\textbf{6}
(2000),
269--317.

\bibitem[MN]{MN}
{A.\,Mudrov and M.\,Nazarov},
\textit{On irreducibility of modules over twisted Yangians},
in preparation.

\bibitem[N]{N2}
{M.\,Nazarov},
\textit{Representations of twisted
Yangians associated with skew Young diagrams},
{Selecta Math.}
\textbf{10}
(2004),
71--129.

\bibitem[NT]{NT}
M.\,Nazarov and V.\,Tarasov,
\textit{On irreducibility of tensor products of Yangian modules
associated with skew Young diagrams},
Duke Math. J.
{\bf 112}
(2002),
343–-378.

\bibitem[O1]{O1}
{G.\,Olshanski},
\textit{Extension of the algebra $U(g)$ for infinite-dimensional classical
Lie algebras $g$, and the Yangians $Y(gl(m))$},
{Soviet Math.\ Dokl.}
\textbf{36}
(1988),
569--573.

\bibitem[O2]{O2}
{G.\,Olshanski},
\textit{Twisted Yangians and infinite-dimensional classical Lie algebras},
{Lecture Notes Math.}
\textbf{1510}
(1992),
103--120.

\bibitem[PP]{PP}
{A.\,Perelomov and V.\,Popov},
\textit{Casimir operators for semi-simple Lie groups},
{Math. USSR Izv.}
\textbf{2}
(1968),
1313--1335.

\bibitem[TV]{TV}
{V.\,Tarasov and A.\,Varchenko},
\textit{Duality for Knizhnik-Zamolodchikov and dynamical equations},
{Acta Appl. Math.}
\textbf{73}
(2002),
141--154.

\bibitem[T]{T}
{J.\, Tits},
\textit{Normalisateurs de tores. I. Groupes de Coxeter \'etendu},
J. Algebra
\textbf{4}
(1966),
96--116.

\bibitem[W]{W}
{H.\,Weyl},
\textit{Classical Groups, their Invariants and Representations},
Princeton University Press, Princeton, 1946.

\bibitem[Z]{Z1}
D.\,Zhelobenko,
\textit{Extremal cocycles on Weyl groups},
{Funct. Anal. Appl.}
{\bf 21}
(1987),
183--192.

\end{thebibliography}
\end{document}